\documentclass[11pt,reqno, final]{amsart}

\usepackage{amsmath,amssymb,amsthm,bbm,bm}
\usepackage{color}
\usepackage[colorlinks=true, allcolors=blue]{hyperref}
\usepackage{mathrsfs}
\usepackage{mathtools}

\usepackage[noabbrev,capitalize,nameinlink]{cleveref}
\crefname{equation}{}{}
\usepackage[noadjust]{cite}

\usepackage{fullpage}

\usepackage{graphics}
\usepackage{pifont}
\usepackage[T1]{fontenc}
\usepackage{tikz}
\usetikzlibrary{arrows.meta}
\usepackage{environ}
\usepackage{framed}
\usepackage{url}
\usepackage[linesnumbered,ruled,vlined]{algorithm2e}
\usepackage[noend]{algpseudocode}
\usepackage[labelfont=bf]{caption}
\usepackage{cite}
\usepackage{framed}
\usepackage[framemethod=tikz]{mdframed}
\usepackage{appendix}
\usepackage{graphicx}
\usepackage[textsize=tiny]{todonotes}
\usepackage{tcolorbox}
\allowdisplaybreaks[1]
\usepackage{enumerate}
\usepackage[shortlabels]{enumitem}
\crefformat{enumi}{#2#1#3}
\crefrangeformat{enumi}{#3#1#4 to~#5#2#6}
\crefmultiformat{enumi}{#2#1#3}
{ and~#2#1#3}{, #2#1#3}{, and~#2#1#3}

\crefname{algocf}{Algorithm}{Algorithms}

\crefname{equation}{}{} %remove ``Equation''
\AtBeginEnvironment{appendices}{\crefalias{section}{appendix}} %appendices

\usepackage[color]{showkeys} %add in 'final' into parameter to remove showkeys
\colorlet{refkey}{orange!20}
\colorlet{labelkey}{blue!30}

\crefname{algocf}{Algorithm}{Algorithms}

\numberwithin{equation}{section}
\newtheorem{theorem}{Theorem}[section]

\newtheorem{lemma}[theorem]{Lemma}
\newtheorem{claim}[theorem]{Claim}
\crefname{claim}{Claim}{Claims}

\newtheorem{corollary}[theorem]{Corollary}

\newtheorem*{question*}{Question}

\theoremstyle{definition}
\newtheorem{definition}[theorem]{Definition}

\newtheorem*{definition*}{Definition}

\theoremstyle{remark}
\newtheorem{remark}[theorem]{Remark}

\newcommand{\mb}{\mathbb}

\newcommand{\mc}{\mathcal}
\newcommand{\mf}{\mathfrak}
\newcommand{\mr}{\mathrm}

\newcommand{\on}{\operatorname}

\newcommand{\wt}{\widetilde}

\newcommand{\eps}{\varepsilon}
\renewcommand{\epsilon}{\varepsilon}

\allowdisplaybreaks

\title{A Toolkit for Robust Thresholds}

\author[Pham]{Huy Tuan Pham}
\address{Department of Mathematics, Stanford University}
\email{huypham@stanford.edu}
\author[Sah]{Ashwin Sah}
\author[Sawhney]{Mehtaab Sawhney}
\address{Department of Mathematics, Massachusetts Institute of Technology}
\email{\{asah,msawhney\}@mit.edu}

\author[Simkin]{Michael Simkin}
\address{Center of Mathematical Sciences and Applications, Harvard University}
\email{msimkin@cmsa.fas.harvard.edu}

\thanks{Pham is supported by a Two Sigma Fellowship. Sah is supported by the PD Soros Fellowship. Sawhney is supported by the Churchill foundation. Sah and Sawhney are supported by NSF Graduate Research Fellowship Program DGE-2141064. Simkin is supported by the Center of Mathematical Sciences and Applications at Harvard University.}

\begin{document}

\begin{abstract}
Consider a host hypergraph $G$ which contains a spanning structure due to minimum degree considerations. We collect three results proving that if the edges of $G$ are sampled at the appropriate rate then the spanning structure still appears with high probability in the sampled hypergraph. We prove such results for perfect matchings in hypergraphs above Dirac thresholds, for $K_r$-factors in graphs satisfying the Hajnal--Szemer\'edi minimum degree condition, and for bounded-degree spanning trees. In each case our proof is based on constructing a spread measure and then applying recent results on the (fractional) Kahn--Kalai conjecture connecting the existence of such measures with probabilistic thresholds.

For our second result we give a shorter and more general proof of a recent theorem of Allen, B\"ottcher, Corsten, Davies, Jenssen, Morris, Roberts, and Skokan which handles the $r=3$ case with different techniques. In particular, we answer a question of theirs with regards to the number of $K_r$-factors in graphs satisfying the Hajnal--Szemer\'edi minimum degree condition.
\end{abstract}

\maketitle

\section{Introduction}\label{sec:intro}

A central pursuit in both random and extremal (hyper)graph theory is to determine thresholds for various properties. In the random setting the natural question is to determine for which $p=p(n)$ the random graph $\mathbb{G}(n,p)$ satisfies a particular property with high probability.\footnote{A sequence of events, indexed by $n$, holds \textit{with high probability (w.h.p.)} if the probabilities of their occurrence tend to $1$ as $n \to \infty$.} In the extremal setting it is natural to consider the \textit{minimum-degree}, or \textit{Dirac}, threshold, i.e., for which $d=d(n)$ every $n$-vertex graph $G$ with $\delta(G) \geq d$ satisfies a particular property.

Prominent examples include thresholds for connectivity~\cite{ER59}, for containing a fixed size subgraph \cite{ER60}, for containing a perfect matching \cite{ER64,ER66,JKV08,KO06,HPS09,Khan16,RRS09}, for containing a Hamilton cycle \cite{Pos76,Dir52}, for Ramsey properties \cite{RR95,FL07}, for containing a clique (or subgraph) factor \cite{JKV08,CH63,HS70}, and for containing a given bounded-degree spanning tree \cite{Mon19,KSS95}.

One interpretation of probabilistic thresholds, suggested by Krivelevich, Lee, and Sudakov \cite{KLS14}, is as a measure of \textit{robustness}. For example, $K_n$ is extremely robust with respect to containing a perfect matching since the threshold for this property in $\mb{G}(n,p)$ is $\log n / n$. Below this density w.h.p.~isolated vertices, which are very simple local obstructions, begin to appear.

In this work we collect several results that exemplify the philosophy that if the minimum degree of a hypergraph $G$ is above the minimum-degree threshold for a particular property then it not only satisfies the property but satisfies it robustly. That is, up to multiplicative constants, the probabilistic threshold for the property does not depend on whether one considers subgraphs of $G$ or subgraphs of the complete hypergraph. We refer the reader to a survey of Sudakov \cite{Sud17} where a number of previous results in this direction are collected.

In order to prove our results we take advantage of the recently established connection between so-called \textit{spread measures} and thresholds \cite{FKNP21}. In fact, a second purpose of this paper is to demonstrate methods to construct such measures. The techniques we use are related to regularity, robust perfect matchings, random greedy algorithms, and iterative absorption.

We prove robustness results in three settings: (i) the Dirac threshold for containing a perfect matching in a $k$-uniform hypergraph, (ii) the Hajnal--Szemer\'edi minimum degree condition for containing a $K_r$-factor, and (iii) the Dirac threshold for containing bounded-degree spanning trees. It is worth mentioning that finding the probabilistic thresholds for these properties in random subgraphs of the complete (hyper)graph were major open problems in the field. The thresholds for (i) and (ii) were found by Johansson, Kahn, and Vu \cite{JKV08} and the threshold for (iii) was found by Montgomery \cite{Mon19}. The fact that we are able to consider the minimum-degree setting for all three in just one paper speaks to the strength of the spread measure results in \cite{FKNP21}.

We remark that very recently Allen, B\"ottcher, Corsten, Davies, Jenssen, Morris, Roberts, and Skokan \cite{ABCDJMRS22} proved (ii) for $r=3$, but without spread techniques and using techniques related to Johansson--Kahn--Vu \cite{JKV08}. Results falling under (i) have been independently proved by Kang, Kelly, K{\"u}hn, Osthus, and Pfenninger \cite{KKKOP22} using a similar strategy but with differences in implementation. Additionally, for the Dirac threshold for $(k-1)$-degrees, they prove a sharper ``stability'' version; we refer to the discussion after \cref{thm:pm} for more details.

The idea of using spread measures to bound thresholds comes from the \textit{fractional expectation threshold vs.~threshold} conjecture of Talagrand \cite{Tal10} which was only recently established in a breakthrough due to Frankston, Kahn, Narayanan, and Park \cite{FKNP21}. We note that this result is a fractional version of the Kahn--Kalai \textit{expectation threshold vs.~threshold} conjecture \cite{KK07} and that the full conjecture, which implies the fractional version, was resolved in very recent work of Park and the first author \cite{PP22}. For our purposes, the crucial corollary is a connection between so-called \emph{spread measures} and thresholds.

\begin{definition}\label{def:spread-measure}
Consider a finite ground set $Z$ and fix a nonempty collection of subsets $\mc{H} \subseteq 2^Z$. Let $\mu$ be a probability measure on $\mc{H}$. For $q > 0$ we say that $\mu$ is \textit{$q$-spread} if for every set $S \subseteq Z$:
\[\mu \left( \{ A \in \mc{H}\colon S \subseteq A \} \right) \le q^{|S|}.\]
\end{definition}

For a finite set $Z$ and $p \in [0,1]$ we denote by $Z(p)$ the binomial distribution on subsets of $Z$ where each vertex is present with probability $p$, independently of the other vertices. For a hypergraph $H$ we abuse notation and write $H(p)$ instead of $E(H)(p)$, i.e., the binomial distribution on the edge set of $H$.

\begin{theorem}[{From \cite[Theorem~1.6]{FKNP21}}]\label{thm:FKNP}
There exists a constant $C = C_{\ref{thm:FKNP}} > 0$ such that the following holds. Consider a non-empty ground set $Z$ and fix a nonempty collection of subsets $\mc{H}\subseteq 2^Z$. Suppose that there exists a $q$-spread probability measure on $\mc{H}$. Then $Z(\min(Cq\log|Z|,1))$ contains an element of $\mc{H}$ as a subset with probability $1-o_{|Z|\to\infty}(1)$.
\end{theorem}

We note that given $\mc{H} \subseteq 2^Z$, determining the expectation-threshold or fractional expectation-threshold is a difficult task in general. In particular, the fractional expectation-threshold is the solution to a linear program with variables corresponding to the collection of subsets in $\mc{H}$, which for many applications is of exponential size. (The expectation-threshold is in general an \emph{integer} linear program.) Nevertheless, the results of \cite{FKNP21,PP22} immediately imply several previously difficult results including the threshold for containing a perfect matching in a random hypergraph (due to Johansson, Kahn, and Vu \cite{JKV08}) or containing a given bounded-degree spanning tree in a random graph (due to Montgomery \cite{Mon19}). A crucial factor in these applications is that the uniform distribution has optimal spread. This follows, for example, from the fact that vertex permutations act transitively on the objects in question (see \cite[Section~7]{FKNP21}). In contrast, in the minimum-degree setting using vertex permutations is a non-starter, and it is not obvious that the uniform distribution has sufficiently small spread to find the thresholds. The focus of this work is on providing methods for constructing spread measures in situations where one cannot rely on the ``spread from vertex permutations''. A similar difficulty is present in work on the threshold for containing a Steiner triple system \cite{SSS22,YKKMO22,Kee22,JP22}, where the connection between spread and thresholds was exploited alongside other tools for understanding Steiner systems and Latin squares.

We will now state each of our results (i-iii), outline their proofs, and discuss relations to existing literature. Before doing so we remark that the proofs of the three result categories can be read independently of each other. We describe the paper's structure in detail below, in \cref{sec:organization}.

\subsection{Perfect matchings in Dirac hypergraphs}\label{sub:dirac}

For a $k$-uniform hypergraph $H$ and $S \subseteq V(H)$, let $\deg_H(S)$ be the number of hyperedges containing $S$. For $1 \le \ell \le k$ let $\delta_\ell(H)$, the minimum $\ell$-degree of $H$, be the smallest degree of an $\ell$-set in $V(H)$.

\begin{definition}\label{def:dirac-threshold}
For integers $1 \le \ell \le k$ and $n$, with $n \in k\mb N$, let $t(n,k,\ell)$ be the smallest $d$ such that every $n$-vertex $k$-uniform hypergraph with $\delta_\ell(H) \ge d$ contains a perfect matching. The \textit{$\ell$-degree (Dirac) threshold} for perfect matchings in $k$-uniform hypergraphs is
\[\delta_{\ell,k}^+ \coloneqq \lim_{\substack{n\to \infty\\k|n }} \frac{t(n,k,\ell)}{\binom{n}{k-\ell}}.\]
\end{definition}

\begin{remark}\label{rmk:large-degree-implies-matching}
The existence of the limit in \cref{def:dirac-threshold} is not \emph{a priori} clear but was proven as \cite[Theorem~1.2]{FK22}. The actual value of $\delta_{\ell,k}^+$ is known in only a few cases. For instance, it is an open problem to find $\delta_{1,k}^+$ for $k \geq 5$.

An immediate consequence of this definition is that for every $1 \le \ell < k$ and $\eps>0$ there exists some $n_{\ell,k,\eps}$ such that every $k$-uniform hypergraph $H$ on $n \ge n_{\ell,k,\eps}$ vertices with $n \in k \mb N$ and $\delta_\ell(H) \ge (\delta_{\ell,k}^+ + \eps) \binom{n}{k-\ell}$ contains a perfect matching.
\end{remark}

Johansson, Kahn, and Vu \cite{JKV08} proved that the threshold for the appearance of perfect matchings in $\mathbb{G}^{(k)}(n,p)$ has order $\log n / n^{k-1}$. We show that this holds more generally for binomial random subgraphs of any hypergraph satisfying the minimum-degree condition above.

Recall that for a hypergraph $\mc{H}$ and $p \in [0,1]$ we write $\mc{H}(p)$ for the random hypergraph where each hyperedge of $\mc{H}$ is retained with probability $p$, independently of all other choices.

\begin{theorem}\label{thm:pm}
Let $\eps > 0$ and $k,\ell \in \mb{N}$ be fixed. There exists $C = C_{\ref{thm:pm}}(\ell,k,\eps)>0$ such that the following holds. Let $\mc{H}$ be an $n$-vertex hypergraph with $k | n$ and $\delta_\ell(\mc{H}) \ge (\delta_{\ell,k}^+ + \eps) \binom{n}{k-\ell}$. Then $\mc{H}(C\log n/n^{k-1})$ contains a perfect matching with high probability.
\end{theorem}

\begin{remark}\label{rmk:pm}
An earlier version of this paper \cite{PSSS22v2} contained an error in the proof of the $\ell>1$ case of \cref{thm:pm} (see \cref{rmk:ell=1 regular complex proof}, below). This is given a simple correction in the current version. Concurrently with posting the previous version, Kang, Kelly, K{\"u}hn, Osthus, and Pfenninger \cite{KKKOP22} released a proof of \cref{thm:pm} (valid for all $\ell\ge 1$). Additionally, they proved that when $\ell = k-1$ the condition $\delta_\ell(\mc H) \geq (\delta_{\ell,k}^+ + \varepsilon) n$ can be replaced by the condition $\delta_\ell(\mc H) \geq t(n,k,k-1)$. Hence, the $\epsilon$ in the minimum degree (and thus also the $\epsilon$-dependence in the threshold) can be removed, thereby providing a tight result in this case. In particular, they give a robust version of the minimum degree result of R\"{o}dl, Ruci\'{n}ski, and Szemer\'{e}di~\cite{RRS09}. Kang, Kelly, K{\"u}hn, Osthus, and Pfenninger also proved that \cref{thm:pm} can be extended to ``optimal matchings'' in cases when $n \notin k\mb{N}$; we note that our techniques also extend to such ``optimal matchings''.

In terms of the differences between the two proofs, \cite{KKKOP22} also use the strategy of exhibiting a sufficiently spread distribution on perfect matchings. However, the approaches differ in the specifics of constructing this distribution. In particular, \cite{KKKOP22} relies on hypergraph regularity whereas we do not.
\end{remark}

In terms of robustness for the Dirac degree threshold for graphs, a result of Sudakov and Vu~\cite{SV08} shows that for $p=\omega(\log n/n)$ any subgraph of $\mb{G}(n,p)$ with minimum degree $(1/2+o(1))np$ has a perfect matching. This \textit{local resilience} property is strictly stronger than the robustness of \cref{thm:pm} in the case $(k,\ell)=(2,1)$ since such a ``minimum-degree subgraph'' can be specified by intersecting the minimum degree host with the random graph $\mb{G}(n,p)$. It remains a problem of substantial interest to derive any such resilient threshold (or closely related universality result) from spread related techniques.

For hypergraphs, analogues of results of Sudakov--Vu \cite{SV08} were proven by Ferber and Kwan~\cite{FK22} for $k$-uniform hypergraphs and $p\gtrsim\max\{n^{-k/2+o(1)},n^{-k+2}\log n\}$. In particular, for $G\sim\mb{G}^{(k)}(n,p)$ with $p$ above the specified threshold, with high probability any $G'\subseteq G$ satisfying $\delta_{\ell}(G')\ge {(\delta_{\ell,k}^+ + \eps)\binom{n-\ell}{k-\ell}p}$ contains a perfect matching. (The case when $\ell = k-1$ was handled in earlier work of Ferber and Hirschfeld \cite{FH20}.) Note that this result is nontrivial only when $p = \Omega(n^{\ell-k}\log n)$ since otherwise one has $\delta_{\ell}(G) = 0$ and there are no subgraphs satisfying the hypothesis of the statement. This indicates a difference between the notion of robustness considered in our work and the notion of local resilence considered in previous works as the probabilistic ranges of interest in \cref{thm:pm} and \cite[Conjecture~1.3]{FK22} only coincide when $\ell = 1$.

We now briefly outline the proof of \cref{thm:pm}; the short proof is given in \cref{sec:pm-hyper}. Fix integers $k > \ell \ge 1$, $\eps>0$, and a hypergraph $\mc{H}$ satisfying the conditions of \cref{thm:pm}. By \cref{thm:FKNP}, it suffices to construct an $O(1/n^{k-1})$-spread distribution on perfect matchings in $\mc{H}$.

We construct such a spread distribution on matchings via \textit{iterative absorption}, which was first introduced by K\"uhn and Osthus \cite{KO13} and Knox, K\"uhn, and Osthus \cite{KKO15} to prove results on Hamilton decompositions. This powerful method has played a prominent role in several recent breakthroughs. We mention only a proof that combinatorial designs exist \cite{GKLO16} and the proof of the Erd\H{o}s--Faber--Lov\'asz conjecture \cite{KKKMO21}. The method is substantially simplified in the context of perfect matchings.

Fix $\eta \ll \eps$ (independent of $n$) and uniformly at random choose a \textit{vortex} $V(\mc{H}) = V_0 \supseteq V_1 \supseteq \cdots \supseteq V_N = X$ with $|X| \le n^{1/(k+1)}$ and $|V_{i+1}| \approx \eta|V_i|$ for every $i$.

The heart of the construction is a \textit{cover-down} procedure: Suppose that for $0 \le i <\ell$ we have constructed a matching $M_i\subseteq\mc{H}$ that covers $V_0 \setminus V_i$ and only a small number of vertices in $V_i$. We then construct an $O(1/|V_i|^{k-1})$-spread distribution on matchings $M_{i+1} \supseteq M_i$, contained in $\mc{H}$ that cover all vertices of $V_0 \setminus V_{i+1}$ and only a small number of vertices in $V_{i+1}$.

Inductively applying the cover-down procedure we obtain a matching $M_{\ell}$ that covers $V_0 \setminus X$, and only a small number of vertices in $X$. In fact, we will ensure that the number of vertices in $V(M_\ell) \cap X$ is small enough that $\delta_\ell(\mc{H}[X \setminus V(M)]) \ge (\delta_{\ell,k}^+ + \eps/2) \binom{|X|}{k-\ell}$, so by \cref{rmk:large-degree-implies-matching} there exists a perfect matching $\wt{M} \subseteq \mc{H}[X \setminus V(M)]$. Therefore $M_\ell \cup \wt{M} \subseteq H$ is a perfect matching.

To see why this procedure has a spread of $O(1/n^{k-1})$, we heuristically analyze the probability that any specific hyperedge $T \in \mc{H}$ is used in $M$. (Since the goal is to apply \cref{thm:FKNP} in the actual analysis one needs to make the analogous calculation for any \textit{set} of hyperedges.) For any $i$, the probability that $T$ is spanned by $V_i$ is $(|V_i|/n)^k$. Since the matching $M_{i+1}$ is $O(1/|V_i|^{k-1})$-spread, the probability that $T \in M_{i+1}\setminus M_i$ is $O((|V_i|/n)^k\! / |V_i|^{k-1}) = O(|V_i|/n^k) = O(1/n^{k-1})$. Finally, the probability that $T$ is in $\wt{M}$ is at most the probability that $T$ is spanned by $X$, which is $(|X|/n)^k \le 1/n^{k-1}$.

\subsection{\texorpdfstring{$K_r$}{Kr}-factors in edge-percolated graphs}\label{sub:K_r-factor}

The Hajnal--Szemer\'edi theorem \cite{HS70} states that every $n$-vertex graph with minumum degree at least $(1-1/r)n$ and $n \in r\mb N$ contains a $K_r$-factor, and that this is tight. (The $r=3$ case is the earlier Corr\'adi--Hajnal \cite{CH63} theorem.) We consider the robustness of this minimum degree condition. We also give a lower bound on the number of $K_r$-factors in such graphs.

\begin{theorem}\label{thm:robust-corradi-hajnal}
Let $r \in \mb N$. There exists a constant $C=C_{\ref{thm:robust-corradi-hajnal}}(r) > 0$ such that for every graph $G$ on $n \in r \mb N$ vertices with $\delta(G) \ge (1-1/r)n$ the following hold.
\begin{enumerate}[(\arabic{enumi})]
    \item\label{itm:corradi hajnal counting} $G$ contains at least $(n/C)^{(r-1)n/r}$  $K_r$-factors.

    \item\label{itm:corradi hajnal threshold} For $p=C(\log n)^{2/(r(r-1))}n^{-2/r}$ with high probability $G(p)$ contains a $K_r$-factor.
\end{enumerate}
\end{theorem}

We note that $p=((r-1)! \log n)^{2/(r(r-1))}\! / n^{2/r}$ is the (sharp) threshold for the property that every vertex in $\mathbb{G}(n,p)$ is contained in an $r$-clique. Hence, \cref{thm:robust-corradi-hajnal} \ref{itm:corradi hajnal threshold} is optimal up to a constant factor. That this is indeed the threshold for a $K_r$-factor in $\mb{G}(n,p)$ was first proved by Johansson, Kahn, and Vu \cite{JKV08}. Additionally, the $r = 3$ case of \cref{thm:robust-corradi-hajnal} \ref{itm:corradi hajnal threshold} was proved very recently by Allen, B\"ottcher, Corsten, Davies, Jenssen, Morris, Roberts, and Skokan \cite{ABCDJMRS22}. The counting statement in \cref{thm:robust-corradi-hajnal} \ref{itm:corradi hajnal counting} (which is new even for $r = 3$) proves \cite[Conjecture~1]{ABCDJMRS22}.

One difficulty in proving \cref{thm:robust-corradi-hajnal} is that there does not exist a $q$-spread probability measure on $K_r$-factors in $G$ with $q < \binom{n-1}{r-1}^{-2/(r(r-1))} = \Omega(n^{-2/r})$ (indeed, each vertex is contained in at least one $K_r$ that is chosen with probability at least $\binom{n-1}{r-1}^{-1}$). Hence, the best bound one can hope for from applying \cref{thm:FKNP} directly is $O( n^{-2/r} \log n)$. This exceeds the threshold in \cref{thm:robust-corradi-hajnal} by a fractional power of $\log n$. To circumvent this difficulty we prove \cref{thm:robust-ch-clique}, which asserts that an appropriately spread measure exists on perfect matchings in the \textit{$r$-clique complex} of $G$. We then apply \cref{thm:FKNP} to obtain an optimal bound on the threshold for perfect matchings in the clique complex. Finally, we apply results of Riordan \cite{Rio18} to couple percolations of the clique complex with percolations of $G$. This gives the correct logarithmic power in \cref{thm:robust-corradi-hajnal}.

\begin{theorem}\label{thm:robust-ch-clique}
Let $r \in \mb N$. There exists a constant $C=C_r$ such that for every graph $G$ on $n \in r \mb N$ vertices with $\delta(G)\ge (1-1/r)n$ if $\mc{H}$ is the set of $r$-cliques in $G$ then there is a $C/n^{r-1}$-spread probability measure (with respect to ground set $\mc{H}$) on the $K_r$-factors of $G$.
\end{theorem}

The proof of \cref{thm:robust-ch-clique} essentially breaks into two steps. The first is proving a version of \cref{thm:robust-ch-clique} in which the minimum degree condition is replaced by the assumption that one has an $r$-partite graph $G$ where the graph between each pair of parts is super-regular; this is \cref{thm:super-regular}. In the second step we derive \cref{thm:robust-ch-clique} from \cref{thm:super-regular}.

\begin{definition}\label{def:super-regular}
Given a pair of vertex sets $X_1$, $X_2$ in a graph $G$ define $d_G(X_1,X_2) = \frac{e_G(X_1,X_2)}{|X_1||X_2|}$ (when the graph $G$ is clear from context we may omit the subscript). A pair $(A_1,A_2)$ is $\eps$-regular if for all $X_i\subseteq A_i$ with $|X_i|\ge \eps |A_i|$ we have that $|d(A_1,A_2)-d(X_1,X_2)|\le \eps$. We say a pair $(A_1,A_2)$ is $(d,\eps)$-regular if $d(A_1,A_2) = d$. 

Furthermore we say $(A_1,A_2)$ is $(d,\eps,\delta)$-super-regular if it is $(d,\eps)$-regular and for all $v\in A_i$ we have $\on{deg}(v,A_{3-i})\ge \delta |A_{3-i}|$. We say a pair is $(d,\eps)$-super-regular if it is $(d,\eps, d-\eps)$-super-regular. Finally we say a pair $(A_1,A_2)$ is $(d^{+},\eps)$-super-regular if it is $(d',\eps)$-super-regular for some $d'\ge d$.
\end{definition}

\begin{theorem}\label{thm:super-regular}
Fix $r\ge 2$ and suppose $1/n\ll\eps\ll d$. Let $G=(V,E)$ be an $r$-partite graph on $V=\bigcup_{i=1}^rA_i$ where $|A_i|=n$ for all $i\in[r]$. Suppose $G[A_i,A_j]$ is $(d_{i,j},\eps)$-super-regular with $d_{i,j}\ge d$ for all $i\neq j$. Let $\mc{H}$ be the $r$-uniform hypergraph where edges in $\mc{H}$ correspond to $r$-partite cliques of $G$. There exists an $O_{d,\eps}(1/n^{r-1})$-spread distribution on the set of perfect matchings in $\mc{H}$. 
\end{theorem}

We note that the $r=3$ version of \cref{thm:super-regular} immediately implies (along with work of Riordan \cite{Rio18} which provides a coupling to the $K_r$-factor version) the crucial technical statement in \cite{ABCDJMRS22} and one of the main contributions of this work is providing a pair of short proofs of \cref{thm:super-regular}. 

The first proof of \cref{thm:super-regular} proceeds via induction on $r$. The base case $r = 2$ is the key step; for the inductive step one chooses a spread matching between a pair of parts and then constructs an auxiliary $(r-1)$-partite graph where an edge chosen between the initial pair of parts is connected to a vertex in the remaining parts if and only if both endpoints connect to the vertex. A careful application of the union bound proves that the associated graph is super-regular and therefore the key step is the case $r = 2$. For this, the crucial idea is to consider a subgraph $G'$ of the underlying bipartite graph $G$ where one chooses to keep a uniformly random large constant number of neighbors of each vertex (an edge is kept if either of its vertices wish to keep it). This subgraph is trivially seen to be appropriately spread by construction and the desired result follows by using Hall's theorem to verify that $G'$ has a perfect matching w.h.p. We note that this idea of constructing a ``spread'' matching via considering a random subgraph where each vertex has a constant number of out-neighbors also plays a role in forthcoming work of the last three authors on the planar assignment problem \cite{SSSAlgo} and will also be used in the proof of \cref{thm:spanning-trees}. The proof here is given in the short \cref{sec:Kr-short}.

The second proof follows closely along the lines of the proof of \cref{thm:pm}; in particular one considers the set of $r$-cliques as hyperedges in an $r$-uniform hypergraph and the heuristic derivation of spread in \cref{sub:dirac} is unchanged. The changes between the proof of \cref{thm:pm} and \cref{thm:super-regular} are due to a ``regularity boosting'' procedure in the cover-down step and finding the perfect matching in the final vortex set. The second issue is handled immediately by the influential blow-up lemma of Koml\'os, S\'ark\"ozy, and Szemer\'edi \cite{KSS01} while the first issue can be handled by a straightforward regularity boosting procedure based on the counting lemma. This proof is given in \cref{sec:Kr-IA}.

The basic strategy to deduce \cref{thm:robust-ch-clique} from \cref{thm:super-regular} is as follows: one applies Szemer\'edi's regularity lemma, finds a $K_r$-factor between the reduced graph of the regularity partition, and then uses \cref{thm:super-regular} for each $K_r$ in the factor on the reduced graph. However, this sketch is a gross oversimplification of the necessary stability analysis as one is forced to deal with various exceptional vertices which occur in the partition. Our stability analysis is combination of those given in \cite{KSS01} and \cite{ABCDJMRS22}. In particular we handle the case where the underlying graph $G$ is far from being $r$-partite using essentially an identical argument to that given in \cite{ABCDJMRS22} (for $r=2$ there is an additional special case where the graph is near the union of two complete graphs). The extremal case, however, where there is a subset of size $|V(G)|/r$ which is nearly empty, is handled via an argument closely related to that given in \cite{KSS01}; in particular we rely on the devices of special stars (which are somewhat simpler in the $K_r$-factor case) used throughout \cite{KSS01} to rebalance parts of a vertex partition so that one can convert various density constraints into minimum degree constraints. Furthermore the various necessary modifications can be made in a ``spread'' manner, which completes our analysis. The stability analysis is carried out in \cref{sec:robust-counting}.

Finally, we note that \cref{thm:robust-corradi-hajnal} is not the first result establishing a robust version of the Hajnal--Szemer\'{e}di theorem. The counting statement in \cref{thm:robust-corradi-hajnal} generalizes a result of S\'{a}rk\"{o}zy, Selkow, and Szemer\'{e}di \cite{SSS03}, where they proved that a graph with minimum degree $d \geq n/2$ has at least $(cn)^{n/2}$ perfect matchings for a small constant $c$. We note that Cuckler and Kahn \cite{CK09} proved that the optimal constant is $c=d/(en)+o(1)$ (and this is witnessed by $d$-regular graphs); the extension of such a result to $K_r$-factors would be of substantial interest. Second, a resilience version of the Corr\'{a}di--Hajnal theorem \cite{CH63} (which is the $r=3$ case of the Hajnal--Szemer\'{e}di theorem) was proven by Balogh, Lee, and Samotij \cite{BLS12}: with high probability $\mb{G}(n,p)$ with $p=\omega((\log n/n)^{1/2})$ is such that every subgraph with minimum degree $(2/3+o(1))np$ contains a triangle factor covering all but $O(p^{-2})$ vertices. This theorem is optimal both with the range of $p$ and the size of the exceptional set of vertices.

\subsection{Bounded-degree spanning trees}\label{sub:spanning}
Using a precursor of the influential blow-up lemma, Koml\'os, S\'ark\"ozy, and Szemer\'edi \cite{KSS95} proved that for every $\Delta \in \mathbb N$ and $\delta > 0$ every sufficiently large $n$-vertex graph with minimum degree at least $(1/2+\delta)n$ contains a copy of every spanning tree of maximum degree at most $\Delta$. Since there exist disconnected $n$-vertex graphs with minimum degree $\lceil n/2 \rceil - 1$, the ``$1/2$'' cannot be improved.

Regarding the corresponding threshold in random graphs, Montgomery \cite{Mon19} proved that w.h.p.~for any given spanning tree with degrees bounded by $\Delta$, the random graph $\mathbb{G}(n,O_\Delta(\log n/n))$ contains a copy. (In fact, Montgomery proved the stronger \textit{universality} result that w.h.p.~$\mathbb{G}(n,O_\Delta(\log n /n))$ contains a copy of \textit{every} spanning tree with maximal degree at most $\Delta$.) Note that a lower bound of $\log n /n$ for this threshold follows from considering the presence of isolated vertices. Therefore Montgomery's result is tight up to a multiplicative constant.

Our third result is that graphs satisfying the K\'omlos--S\'ark\"ozy--Szemer\'edi \cite{KSS95} minimum degree condition are robust with respect to containing a given bounded-degree spanning tree.

\begin{theorem}\label{thm:spanning-trees}
For every $\Delta \in \mb{N}$, $\delta > 0$, there exists $C = C_{\ref{thm:spanning-trees}}(\Delta,\delta)$ such that the following holds. Suppose that $G$ is an $n$-vertex graph satisfying $\delta(G) \ge (1/2 + \delta)n$ and $T$ is an $n$-vertex tree with $\Delta(T) \le \Delta$. Then w.h.p.~$G(C\log n /n)$ contains a copy of $T$.
\end{theorem}

Our proof of \cref{thm:spanning-trees} closely follows that of \cite{KSS95} and in fact we quote the main preprocessing statement from their proof. Roughly, the proof in \cite{KSS95} proceeds by embedding most vertices of the tree via a greedy process into a regularity decomposition of the graph $G$. The remaining vertices, which either form paths of length $4$ or stars of various sizes, are embedded via specialized lemmas which are essentially special cases of the blow-up lemma. The initial embedding of the large portion of the tree can easily be made $O(1/n)$-spread by using a \textit{random} greedy algorithm, and observing that there are always $\Omega(n)$ choices for where to place each vertex. The paths of length $4$ lemma can be made spread by essentially quoting the $r=3$ case of \cref{thm:super-regular} and the embedding of stars can be made spread using the robust bipartite matchings discussed in \cref{sec:Kr-short}. The details of the proof are in \cref{sec:bounded-degree-trees}.

We remark that \cref{thm:spanning-trees} is not the first robust version of the results of Koml\'os, S\'ark\"ozy, and Szemer\'edi; Balogh, Csaba, and Samotij \cite{BCS11} consider the resilience of $\mb{G}(n,p)$ with respect to containing \textit{almost} spanning trees. Specifically, for every $\Delta \in \mb{N}$ and $\varepsilon,\eta>0$, in $\mb{G}(n,p)$ with $p = \Omega_{\eta,\epsilon,\Delta}(1/n)$, w.h.p.~any subgraph $G$ with at least $(1/2+\eta)$-fraction of the edges at each vertex contains any tree on $(1-\epsilon)n$ vertices with degree bounded by $\Delta$. Related work has also considered such universality with respect to containment of almost spanning structures of limited bandwidth in $\mb{G}(n,p)$ and resilient subgraphs; we refer the reader to \cite{BKT13} and references therein.

\subsection{Organization}\label{sec:organization}
In \cref{sec:prelim} we collect basic preliminaries which will be used throughout the paper. In \cref{sec:pm-hyper} we give the proof of \cref{thm:pm}. In \cref{sec:Kr-short} we give the first proof of \cref{thm:super-regular} and in \cref{sec:Kr-IA} we give the second proof of \cref{thm:super-regular}. In \cref{sec:robust-counting} we prove the necessary stability analysis to deduce \cref{thm:robust-ch-clique} and thus \cref{thm:robust-corradi-hajnal}. Finally in \cref{sec:bounded-degree-trees} we prove \cref{thm:spanning-trees}.

\subsection{Notation}
We write $[n]=\{1,\ldots,n\}$. We write $f=O(g)$ to mean that $f\le Cg$ for some absolute constant $C$, and $g=\Omega(f)$ to mean the same. We put a subscript, say $O_\eps$, to mean that the constant may depend on some outside parameter. We write $f=o(g)$ if for all $c > 0$ we have $f\le cg$ once the implicit growing parameter (usually $n$) is large enough, and $g=\omega(f)$ means the same. For parameters $\alpha,\beta$, we write $\alpha\ll\beta$ to mean that $\alpha>0$ is less than some suitably chosen function of $\beta$. Finally, we may implicitly round large real numbers to integers if they are counting objects or performing a similar role, and the exact number is not required to be precise.

As discussed earlier, for a hypergraph $H$ and subset $S\subseteq V(H)$ we write $\deg_H(S)$ for the number of hyperedges of $H$ containing $S$. If $S, A \subseteq V(H)$ we also write $\deg_H(S,A)$ for the number of hyperedges $h$ in $H$ such that $S \subseteq h$ and $h \setminus S \subseteq A$.

For a graph $G$ and two sets $X_1,X_2 \subseteq V(G)$ we write $d_G(X_1,X_2) = e_G(X_1,X_2) / (|X_1||X_2|)$ and $d_G(X_1) = e_G(X_1) / \binom{|X_1|}{2}$.

\subsection*{Acknowledgements}
The first author would like to thank David Conlon and Jacob Fox for helpful comments and suggestions. The second and third authors thank Matthew Kwan and Vishesh Jain for helpful discussions. Part of this research was conducted while the third author was visiting IST Austria, Tel Aviv University, and Cambridge University and they would like to thank these institutions for their hospitality. Finally we thank the authors of \cite{KKKOP22} for pointing out the extension to the non-divisible case suggested in \cref{rmk:pm}. 

We thank Aliaksei Semchankau for identifying an error in the proof of the $\ell>1$ case of \cref{thm:pm} that appeared in a previous version of this paper (cf.\ \cref{rmk:pm,rmk:ell=1 regular complex proof}).

\section{Preliminaries}\label{sec:prelim}
We will repeatedly use the Chernoff bound for binomial and hypergeometric distributions (see for example \cite[Theorems~2.1,~2.10]{JLR00}) without further comment.
\begin{lemma}[Chernoff bound]\label{lem:chernoff}
Let $X$ be either:
\begin{itemize}
    \item a sum of independent random variables, each of which take values in $[0,1]$, or
    \item hypergeometrically distributed (with any parameters).
\end{itemize}
Then for any $\delta>0$ we have
\[\mb{P}[X\le (1-\delta)\mb{E}X]\le\exp(-\delta^2\mb{E}X/2),\qquad\mb{P}[X\ge (1+\delta)\mb{E}X]\le\exp(-\delta^2\mb{E}X/(2+\delta)).\]
\end{lemma}
We also record a lemma for comparing a sequence of random variables to an independent sum of Bernoulli random variables.
\begin{lemma}[{\cite[Lemma~8]{Ram90}}]\label{lem:binomial-comparison}
Let $X_1,\ldots, X_n$ be $\{0,1\}$-valued random variables such that for all $i\in [n]$, we have that $\mb{P}[X_i = 1|X_1,\ldots,X_{i-1}]\le p$ then $\mb{P}[\sum_{i=1}^nX_i\ge t]\le \mb{P}[\mr{Bin}(n,p)\ge t]$ for all $t\ge 0$.
\end{lemma}

Finally, we will use McDiarmid's inequality, which follows from \cite[Lemma 1.2]{McDia89}.

\begin{lemma}\label{lem:McDiarmid}
	Let $X_1,\ldots,X_n$ be independent random variables, with each $X_i$ taking values in a finite set $\Lambda_i$. Let $f:\prod_{i=1}^n \Lambda_i \to \mb{R}$ be a function satisfying: for some $L>0$ if $\vec{x},\vec{y} \in \prod_{i=1}^n \Lambda_i$ differ by at most one coordinate then $|f(\vec{x})-f(\vec{y})| \leq L$. Then, for every $t>0$ there holds
	\[
	\mathbb{P} \left[ |f(X_1,\ldots,X_n) - \mb{E}[f(X_1,\ldots,X_n)]| \geq t \right] \leq 2\exp \left( -2t^2 / (2nL^2) \right).
	\]
\end{lemma}

Next we require a version of the R\"odl nibble which gives a spread distribution over nearly complete hypergraph matchings in a dense situation. We sketch a short proof based on applying the R\"odl nibble given in \cite[Theorem~4.7.1]{AS16} to a harshly subsampled random hypergraph.

\begin{lemma}\label{lem:nibble}
Fix $\gamma,\eta > 0$. There exists some $\delta = \delta_{\ref{lem:nibble}}(\eta)$ such that the following holds. Suppose that $n$ is sufficiently large (in terms of $\gamma$ and $\eta$). Given a $k$-uniform hypergraph $\mc{H}$ on $n$ vertices with at least $\gamma\binom{n}{k}$ hyperedges such that each vertex is in at most $(1+\delta) \gamma \binom{n}{k-1}$ hyperedges of $\mc{H}$, there exists an $O_{k,\gamma,\eta}(1/n^{k-1})$-spread distribution $\mf{p}$ on matchings in $\mc H$ that cover at least $(1-\eta)n$ vertices.
\end{lemma}

\begin{proof}
Let $\mc L \coloneqq \mc H((C / \gamma)n^{1-k})$ for a large constant $C$ (which will be taken to depend only on $\eta$). Note that as long as $n$ is sufficiently large in terms of $C$ then for every distinct $v_1,v_2 \in V(\mc H)$ the codegree $\mr{codeg}_{\mc{L}}(v_1,v_2)$ is stochastically dominated by $\on{Bin}(n^{k-2},n^{-k+5/4})$. Hence,
\[\Pr[ \mr{codeg}_{\mc{L}}(v_1,v_2) \geq 3] \leq (n^{k-2})^3 (n^{5/4-k})^3 = o(n^{-2}).\]
Since there are $O(n^2)$ choices of the vertices $v_1,v_2$, a union bound implies that w.h.p.\ all codegrees in $\mc{L}$ are at most $2$.

Next notice that at least $(1-\delta^{1/3})n$ vertices have degree at least $(1-\delta^{1/2}) \gamma \binom{n}{k-1}$, else $\mc{H}$ would have fewer than the specified number of hyperedges. Hence, as long as $\delta$ is sufficiently small in terms of $C$, the Chernoff and Markov inequalities imply that with probability at least $7/8$ we have that $(1-e^{-\Omega(C^{1/3})})n$ vertices have $\mc{L}$-degree in $[C\pm C^{3/4}]$.

Let $\mc L' \subseteq \mc L$ be the set of edges in $\mc L$ where every vertex has $\mc{L}$-degree at most $C+C^{3/4}$. We will show that for most vertices their $\mc L'$-degree is equal to their $\mc L$-degree. We will do this by bounding the number of edges $e \in \mc H$ such that $e \in \mc L$ and, for some vertex $v \in e$, there holds $\deg_{\mc L}(v) > C+C^{3/4}$.

Let $e \in \mc H$ and a vertex $v \in e$. By definition, $\Pr[e \in \mc L] = (C/\gamma n^{k-1})$. Observe that $\deg_{\mc L}(v) > C+C^{3/4}$ only if $\deg_{\mc L \setminus \{e\}}(v) > C+\frac{1}{2} C^{3/4}$. Applying Chernoff's inequality we obtain $\Pr[\deg_{\mc L \setminus \{e\}}(v) > C+\frac{1}{2} C^{3/4}] = \exp(-\Omega(C^{1/3}))$. Therefore, using the fact that the hypergraph $\mc L \setminus \{e\}$ is independent of the event $e \in \mc L$, we conclude
\begin{align*}
\Pr[e \in \mc L \land \deg_{\mc L}(v) > C+C^{3/4}] &\leq \Pr[e \in \mc L] \times \Pr[\deg_{\mc L \setminus \{e\}}(v) > C+\frac{1}{2} C^{3/4}]\\
&\leq \frac{C}{\gamma n^{k-1}} \exp (-\Omega(C^{1/3})) = \frac{\exp(-\Omega(C^{1/3}))}{n^{k-1}}.
\end{align*}
Since there are fewer than $n^k$ choices of $e \in \mc H$ and $v\in e$, using Markov's inequality we conclude that with probability at least (say) $7/8$, there are fewer than $e^{-\Omega(C^{1/3})} n$ vertices $v$ that are incident in $\mc L$ to a vertex $u$ with $\mc L$-degree less than $C+C^{3/4}$.

Note that for a vertex to have smaller degree in $\mc L'$ than in $\mc L$ it must either have $\mc L$-degree larger than $C+C^{3/4}$ or be incident in $\mc L$ to a vertex with degree larger than $C+C^{3/4}$. Applying a union bound we conclude that with probability at least $3/4$, all but $e^{-\Omega(C^{1/3})}n$ vertices have $\mc L'$-degree $C \pm C^{3/4}$.

Let $\mc G$ be the (good) event that all but $e^{-C^{1/4}}n$ vertices have $\mc L'$-degree $C \pm C^{3/4}$, that the maximal degree in $\mc L'$ is at most $C+C^{3/4}$ (which holds by construction), and that the maximal codegree in $\mc L'$ is at most $2$. Applying a union bound, we conclude that $\Pr[\mc G] \geq 1/2$.

%Fix a vertex $v$ and assume that $\max_{v'\neq v}\mr{codeg}_{\mc L}(v,v')\le 2$ and that its degree is in $[C \pm C^{3/4}]$. Having revealed the outcome of edges containing $v$, we need to prove that it is highly likely that the remaining neighbors in all these edges have degree at most $C+C^{3/4}$, so that the $\mc{L}'$-degree of $v$ is the same as the $\mc{L}$-degree and hence is in $[C\pm C^{3/4}]$. This occurs for $v$ with probability $1-e^{-\Omega(C^{1/3})}$ by Chernoff, having revealed every edge containing $v$ (using that the original degrees in $\mc{H}$ are tightly controlled).

%Therefore we find from Markov that with probability at least $3/4$, $\mc{L}'$ has at least $(1-e^{-\Omega(C^{1/3})})n$ vertices with degree within $C\pm C^{3/4}$, no vertices with degree above $C+C^{3/4}$, and all codegrees bounded by $3$.

Assuming that $C$ is sufficiently large in terms of $\eta$, if $\mc G$ occurs then applying \cite[Theorem~4.7.1]{AS16} to the (at least $(1-e^{-C^{1/4}}n$) non-isolated vertices of $\mc{L}'$ implies the existence of a vertex cover $\mc C' \subseteq \mc L'$ with at most $(1+\eta/(10k))n/k$ edges. Let $\mc C \subseteq \mc C'$ be the set of edges that are disjoint from all others. Provided that $C$ is sufficiently large in terms of $\eta$, a simple calculation shows that $\mc C$ covers at least $(1-\eta)n$ vertices.

We define the measure $\mf p$ by conditioning on $\mc G$ and choosing $\mc C$ as above. To see that it is spread, let $S \subseteq \mc H$. There holds
\[
\Pr [S \subseteq \mc C | \mc G] = \frac{ \Pr [S \subseteq \mc C \land \mc G]}{\Pr[\mc G]} \leq \frac{ \Pr [S \subseteq \mc L]}{1/2} \leq ((2C/\gamma)n^{1-k})^{|S|}.
\]
\end{proof}

Finally we will also require a number of basic definitions regarding regular and super-regular pairs. Recall that have already defined the notion of $(d,\epsilon, \delta)$-super-regular as well as $(d^{+},\epsilon)$-super-regular in \cref{def:super-regular}. We will require the following lemma which allows one to transfer between these notions at the cost of passing to a suitable subgraph; this appears as \cite[Lemma~2.12]{ABCDJMRS22}.
\begin{lemma}\label{lem:plus-to-super-regular}
For every $\epsilon > 0$ and $n = n_{\ref{lem:plus-to-super-regular}}(\epsilon)$ such that the following holds. Consider a bipartite graph on $V_1,V_2$ with parts of size $n$ and $G$ which is $(\epsilon^{2},d^{+})$-super-regular for $d$ such that $4\epsilon\le d\le 1$ and $dn^2\in\mb{N}$. Then there is a spanning subgraph $G'$ of $G$ so that $(V_1,V_2)$ is $(4\epsilon, d)$-super-regular in $G'$.
\end{lemma}

We will next require the counting lemma which counts embeddings of a subgraph $H$ into fixed parts of a collection of regular pairs. This follows immediately from the standard proof of the counting lemma.

\begin{lemma}[Counting lemma]\label{lem:counting-lemma}
Fix $\eps > 0$ and a graph $H$ on vertices $v_1,\ldots,v_k$. There exists an absolute constant $C_H = C_{H,\ref{lem:counting-lemma}}$ such that the following holds. Fix a graph $k$-partite $G = (V,E)$ where $V = \bigcup_{i=1}^{k}A_i$ and suppose that for each $(i,j)\in E(H)$ that $G[A_i,A_j]$ is $(d_{i,j},\eps)$-regular with $d_{i,j}\ge \eps$. Fix any subsets $X_i\subseteq A_i$ with $|X_i|\ge \eps |A_i|$. Then the number of homomorphisms from $H$ to $G$ with $v_i$ mapping into $X_i$ is  
\[\prod_{(i,j)\in E(H)}d_{i,j}\prod_{i=1}^{k}|X_i| \pm C_H\eps \prod_{i=1}^{k}|A_i|.\]
\end{lemma}

We will also require that any $(d,\eps)$ super-regular bipartite graph has a subset of edges such that both it and its complement are super-regular.

\begin{lemma}\label{lem:inner-regularity}
Suppose that $1/n\ll \eps \ll d\le 2/3$. Then given a graph $G = (A_1\cup A_2, E)$ such that $|A_i| = n$ and $(A_1,A_2)$ is $(d,\eps)$-super-regular we have a spanning subgraph $G'\subseteq G$ which is $(d\pm \eps^{1/3}, \eps^{1/3})$-super-regular and such that the bipartite complement $K_{A_1,A_2}\setminus G'$ is $(1-d\pm \eps^{1/3}, \eps^{1/3})$-super-regular.
\end{lemma}
\begin{proof}
Let $A_1' = \{ v \in A_1 \colon \deg_G(v,A_2)\ge (d+2\eps)n \}$ and define $A_2'$ similarly. By $(d,\eps)$-regularity applied to the sets $A_1'$ and $A_2$ (and symmetric) we have that $|A_1'|, |A_2'|\le \eps n$. By removing all edges between $A_1'$ and $A_2'$, note that the degrees in $A_1\setminus A_1'$ and $A_2\setminus A_2'$ are unchanged and we still have a lower bound of $(d-\epsilon)n$ for all such vertices. Now for each vertex $v \in A_i'$ choose a set of approximately $dn$ edges to keep and remove the rest; this is possible to do on a vertex-by-vertex basis as all remaining edges have at most $1$ endpoint in $A_i'$. Furthermore note that in this procedure all vertices outside $A_1'\cup A_2'$ have degrees adjusted by at most $\epsilon n$. Finally note that since at most $2\epsilon n^2$ edges have been modified the desired result follows immediately by the definition of regularity. 
\end{proof}
A similar proof can in fact be used to show the following; this combined with \cref{lem:plus-to-super-regular,lem:inner-regularity} show that any $(d,\eps,\delta)$-super-regular graph will have very super-regular spanning subgraphs.
\begin{lemma}\label{lem:super-boost}
Suppose that $1/n\ll\eps\ll\delta\ll d\le 2/3$. Then given a graph $G = (A_1\cup A_2, E)$ such that $|A_i| = n$ and $(A_1,A_2)$ is $(d,\eps,\delta)$-super-regular we have a spanning subgraph $G'\subseteq G$ which is $(\delta^+,\eps^{1/3})$-super-regular.
\end{lemma}
\begin{proof}[Proof sketch]
Using $(d,\eps)$-regularity, we can show that the vast majority of vertices have degree $(d\pm2\eps)n$, say. Call the exceptional vertices $A_1',A_2'$ as before. We take random $(\delta/d)$-samples of the edges $G[A_1\setminus A_1',A_2\setminus A_2']$, delete the edges $G[A_1',A_2']$, and then for each $v\in A_1'\cup A_2'$ independently choose $\delta n$ edges to keep.
\end{proof}

Furthermore we will need the influential blowup lemma of  Koml\'os, S\'ark\"ozy, and Szemer\'edi \cite[Theorem~1]{KSS97}. 

\begin{theorem}[{\cite[Theorem~1]{KSS97}}]\label{lem:blowup}
Given a graph $R$ of order $r$ and parameters $\delta, \Delta$, there exists $\eps = \eps(\delta, \Delta, r)$ such that the following holds. Suppose that one replaces the vertices of $R$ with sets of size $n_1,\ldots,n_r$ and define a graph $G_1$ where each edge of $R$ is blown up to a $(\delta,\eps)$-super-regular pair and a graph $G_2$ where each edge of $R$ is blown up to a complete graph. Then if a graph $H$ with maximum degree $\Delta$ embeds into $G_2$ it embeds into $G_1$.  
\end{theorem}

Finally we will also require a version of the regularity lemma which respects the minimum degree of the input graph. The version stated appears as \cite[Lemma~2.6]{ABCDJMRS22}.

\begin{definition}\label{def:regularity-setup}
We say a partition $V(G) = V_0 \cup V_1\cup \cdots\cup V_t$ is an $\epsilon$-regular partition if $|V_0|\le \epsilon |V(G)|$, $|V_1| = \cdots = |V_t|$, and all but $\epsilon t^2$ pairs $(V_i,V_j)$ are $\epsilon$-regular. Given an $\epsilon$-regular partition and $d\in[0,1]$, we say that $R$ is the $(\epsilon, d)$-reduced graph with respect to the partition if $V(R) = [t]$ and $(i,j)\in E(R)$ if and only if $(V_i,V_j)$ is a $(d^+,\epsilon)$-regular pair.
\end{definition}

\begin{lemma}[{\cite[Lemma~2.6]{ABCDJMRS22}}]\label{lem:regularity}
For all $\epsilon > 0$ and $m_0\in \mb{N}$, there is $M_0 = M_{\ref{lem:regularity}}(m_0,\epsilon)$ such that the following holds. For all $0<d<\gamma<1$, $n>M_0$, and graphs $G$ with $\delta(G)\ge \gamma n$, there exists an $\epsilon$-regular partition $V_0\cup V_1\cup \cdots\cup V_m$ with $m_0\le m\le M_0$ such that the $(\epsilon,d)$-reduced graph $R$ has $\delta(R)\ge(\gamma-d-2\epsilon)m$. 
\end{lemma}

\section{Robust conditions for perfect matchings in random hypergraphs}\label{sec:pm-hyper}

The first lemma will allow us to construct the vortex; a crucial feature of our analysis is that the randomness in the choice of the vortex is taken into account when calculating the spread. For this reason the lemma guarantees a \textit{distribution} over vortices, rather than the existence of any particular one.

\begin{lemma}[Vortex]\label{lem:vortex}
For every $\alpha > 0$ and $\eps \in (0,1/10)$ there exists some $C=C_{\ref{lem:vortex}}(\eps)$ such that if $\mc{H}$ is a $k$-uniform hypergraph on $n\ge C$ vertices satisfying $\delta_\ell(H) \ge (\alpha + \eps) \binom{n}{k-\ell}$ then there exists a distribution on set sequences $V(\mc{H}) = V_0 \supseteq V_1 \supseteq \cdots \supseteq V_N = X$ with the following properties:
\begin{enumerate}[{\bfseries{V\arabic{enumi}}}]
	\item\label{itm:vortex-set-sizes} For every $0 \le i < N$ there holds $|V_{i+1}| = (1 \pm \eps^2/N^2)\eps^2|V_i|$;
	\item\label{itm:vortex-final-set-size} $|X| \in [n^{1/(k+2)},n^{1/(k+1)}]$;
	\item\label{itm:vortex-high-degrees} For every $S \in \binom{V(\mc{H})}{\ell}$ and every $0 \le i \le N$, there holds $\deg_\mc{H}(S,V_i) \ge (\alpha + \eps/2) \binom{|V_i|}{k-\ell}$;
    \item\label{itm:vortex-prob-guarantee} For every vertex set $\{v_1,\ldots,v_m\} \subseteq V(\mc{H})$ and every vector $\vec{x} \in \{0,\ldots,N\}^m$ there holds
	\[\mb{P}\left[ \bigwedge_{i=1}^m (v_i \in V_{x_i}) \right] \le \prod_{i=1}^m \frac{2|V_{x_i}|}{n}.\]
	\end{enumerate}
\end{lemma}

\begin{proof}
First, consider the distribution on set sequences $V(\mc{H}) = U_0 \supseteq \cdots \supseteq U_N$ obtained as follows: Set $U_0 = V(\mc{H})$. For as long as $|U_i| > n^{1/(k+1)}$, let $U_{i+1}$ be a binomial random subset of $U_i$ of density $\eps^2$. Let $\mathcal E$ be the event that properties \cref{itm:vortex-set-sizes,itm:vortex-final-set-size,itm:vortex-high-degrees} hold. McDiarmid's inequality (\cref{lem:McDiarmid}) and a union bound imply that $\mathcal E$ holds w.h.p.
	
Next, let $\{v_1,\ldots,v_m\} \subseteq V(\mc{H})$ and $\vec{x}\in\{0,\ldots,N\}^m$. Clearly:
\[\mb{P}\bigg[\bigwedge_{i=1}^m(v_i\in U_{x_i})\bigg] = \prod_{i=1}^m \eps^{2x_i}.\]
	
Let $V_0\supseteq\cdots\supseteq V_N$ be the distribution obtained by conditioning $U_0\supseteq\cdots\supseteq U_N$ on the occurrence of $\mathcal E$. By definition, $V_0\supseteq\cdots\supseteq V_N$ satisfies properties \cref{itm:vortex-set-sizes,itm:vortex-final-set-size,itm:vortex-high-degrees}. Furthermore, for every nonempty $\{v_1,\ldots,v_m\}\subseteq V(\mc{H})$ and $\vec{x} \in \{0,\ldots,N\}^m$:
\[\mb{P}\bigg[\bigwedge_{i=1}^m(v_i\in V_{x_i})\bigg] = \mb{P}\bigg[\bigwedge_{i=1}^m(v_i\in U_{x_i})\bigg|\mc{E}\bigg]\le\frac{\mb{P}\big[\bigwedge_{i=1}^m(v_i\in U_{x_i})\big]}{\mb{P}[\mc{E}]}\le (1+o(1))\prod_{i=1}^m\eps^{2x_i}\le \prod_{i=1}^m \frac{2|V_{x_i}|}{n},\]
as desired. The last inequality comes from applying \cref{itm:vortex-set-sizes} iteratively at most $N$ times.
\end{proof}

Next, we show that if $\mc{H}$ satisfies the conditions in \cref{thm:pm} then it contains a large nearly regular subgraph. We note that \cref{lem:regular-subgraph} plays the role of a ``regularity--boosting'' lemma which has various previous applications in iterative absorption.

\begin{lemma}\label{lem:regular-subgraph}
For every $\eps,\gamma>0$ there exists $\gamma' > 0$ such that if $\mc{H}$ is a $k$-uniform hypergraph on $n \ge C_{\ref{lem:regular-subgraph}}(\eps,\gamma)$ vertices and $\delta_\ell(\mc H) \ge (\delta_{\ell,k}^+ + \eps) \binom{n}{k-\ell}$ then there exists a subgraph $\mc{H}'$ with at least $\gamma'\binom{n}{k}$ edges and maximum degree at most $(1+\gamma)\gamma'\binom{n}{k-1}$.
\end{lemma}

\begin{remark}\label{rmk:ell=1 regular complex proof}
In the $\ell=1$ case, \cref{lem:regular-subgraph} has a particularly simple proof. First, let $W \subseteq V(\mc H)$ be a set of at most $k-1$ vertices such that $n-|W|$ is a multiple of $k$. Now, consider $\mc H'' \coloneqq \mc H[V(\mc H) \setminus W]$. The minimum degree assumption implies that $\mc H''$ contains a perfect matching. In fact, for $\gamma'$ sufficiently small, we can iteratively remove $(1+\gamma)\gamma' \binom{n}{k-1}$ such matchings. The union of these matchings satisfies the requirements for $\mc H'$.

Unfortunately, this simple argument does not generalize to $\ell>1$ (an attempt to do so is the error alluded to in \cref{rmk:pm}). This is because we can only iteratively remove $O(n^{k-\ell})$ perfect matchings while guaranteeing the minimum degree condition. We circumvent this difficulty by constructing $\Omega(n^{k-1})$ \textit{approximate} matchings.
\end{remark}

We require the following lemma which follows immediately from work of Ferber and Kwan \cite{FK22}.

\begin{lemma}\label{lem:partial-large}
Fix $\gamma,\epsilon>0$. There exists $\delta>0$ such that the following holds. If $\mc{H}$ is a $k$-uniform hypergraph on $n \ge C_{\ref{lem:partial-large}}(\gamma,\epsilon)$ vertices and all but $\delta \binom{n}{\ell}$ sets of size $\ell$ are contained in at least $(\delta_{\ell,k}^+ + \epsilon/2)\binom{n}{k-\ell}$ edges of $\mc{H}$, then $\mc{H}$ contains a matching with at least $(1-\gamma)n$ vertices.
\end{lemma}

\begin{proof}
Let $Q = Q(k,\gamma,\epsilon)$ divisible by $k$. We assume that $\delta$ is small with respect to $\gamma$ and that $Q$ is large with respect to $\gamma,\varepsilon$ and $k$. By \cite[Lemma~3.4]{FK22} a uniformly random set $S \subseteq V(\mc H)$ of size $Q$ satisfies $\delta_{\ell}(H[Q]) \ge (\delta_{\ell,k}^+ + \epsilon/3)Q$ with probability at least $1-\binom{Q}{k}(\delta + e^{-c_k\epsilon^2 Q})\ge 1-\gamma/2$. By definition of $\delta_{\ell,k}^+$, this minimum degree condition implies that $H[S]$ contains a perfect matching.

Now, randomly partition $V(\mc H)$ into $\lfloor n/Q \rfloor$ sets of size $Q$ and one ``leftover'' set of size at most $Q$. Markov's inequality implies that with probability at least $1/3$ at least $1-3\gamma/4$ fraction of these sets have a perfect matching. The union of these perfect matchings covers at least $(1-\gamma)n$ vertices. This implies the result.
\end{proof}

\begin{proof}[Proof of \cref{lem:regular-subgraph}]
We may assume that $\gamma$ is sufficiently small in terms of $\varepsilon$. We iteratively apply \cref{lem:partial-large} by removing partial matchings of size at least $(1-\gamma/2)n/k$. Notice that if we have removed $\eta\binom{n}{k-1}$ such matchings then we have removed at most $\eta \binom{n}{k-1} \frac{n}{k} \leq \eta n^k$ total edges. Hence, after removing these edges, at most $\eta n^k \binom{k}{\ell} / (\varepsilon n^{k-\ell}\!/2) =  O_k((\eta/\epsilon)\binom{n}{\ell})$ sets of size $\ell$ are contained in fewer than $(\delta_{\ell,k}^+ + \epsilon/2)\binom{n}{k-\ell}$ hyperedges. Therefore, for $\gamma'>0$ sufficiently small in terms of $\gamma,\varepsilon$, and $k$, we can iterate \cref{lem:partial-large} and remove $(1+\gamma)\gamma'\binom{n}{k-1}$ matchings of size at least $(1-\gamma/2)n/k$. We complete the construction by taking $\mc H'$ as the union of these matchings.
\end{proof}

We are ready to state the cover-down lemma.

\begin{lemma}[Cover-down lemma]\label{lem:cover-down}
For every $\eps \in (0,1/(10k!))$ there exists some $C=C_{\ref{lem:cover-down}}(\eps) > 0$ such that the following holds. Let $\mc{H}$ be a $k$-uniform hypergraph with $|V(\mc{H})| \ge C$. Let $U \subseteq V(\mc{H})$ satisfy $|U| = (1 \pm \eps)\eps^2 |V(\mc{H})|$. Suppose that $\mc{H}$ satisfies $\delta_\ell(\mc{H}) \ge (\delta_{\ell,k}^+ + \eps/3)\binom{|V(\mc{H})|}{k-\ell}$ and, for every $S \in \binom{V(\mc{H})}{\ell}$, there holds $\deg_\mc{H}(S,U) \ge (\delta_{\ell,k}^+ + \eps/3) \binom{|U|}{k-\ell}$. Then there exists a $C / |V(\mc{H})|^{k-1}$-spread distribution on matchings $M \subseteq \mc{H}$ that satisfy:
\begin{enumerate}[{\bfseries{C\arabic{enumi}}}]
	\item\label{itm:M-covers-V-U} $M$ covers every vertex in $V(\mc{H}) \setminus U$; and
	\item\label{itm:small-cover-in-U} $M$ covers at most $\eps^2 |U|$ vertices in $U$.
\end{enumerate}
\end{lemma}

\begin{proof}
For notational conciseness we set $V = V(\mc{H})$ and fix a constant $\delta$, sufficiently small, to be chosen later. We first construct the random matching $M$.

\textit{Step 1: Finding an approximately regular subgraph.} We will find a regular subgraph of $\mc{H}[V \setminus U]$ by applying \cref{lem:regular-subgraph}. Let $V' \coloneqq V \setminus U$. Let $\mc{H}' \coloneqq \mc{H}[V']$. Observe that $\delta_\ell(\mc{H}') \ge \delta_\ell(\mc{H}) - |U|^{k-\ell} \ge (\delta_{\ell,k}^+ + 0.3\eps)\binom{|V'|}{k-\ell}$. Assuming that $|V|$ is sufficiently large, we can apply \cref{lem:regular-subgraph} to find $\wt{\mc H}$ with at least $\gamma'\binom{n}{k}$ edges and maximum degree at most $(1+\delta)\gamma'\binom{n}{k-1}$ (for some $\gamma'>0$ as a function of $\epsilon$ and $\delta$).

\textit{Step 2: Finding a spread approximate matching.} We apply \cref{lem:nibble} to $\wt{\mc{H}}$ to find an $O_\eps(1/|V|^{k-1})$-spread random matching $\wt{M} \subseteq \wt{\mc{H}}$ covering all but at most $\eps^6|V(\wt{H})|$ vertices; this is possible given that $\delta$ is sufficiently small with respect to $\epsilon$.

\textit{Step 3: Covering remaining vertices in $V \setminus U$.} Conditioning on $\wt{M}$, let $v_1,\ldots,v_m$ be an enumeration of the uncovered vertices $V \setminus (V(\wt{M}) \cup U)$, noting $m\le\eps^6|V(\wt{\mc H})|$. We extend $\wt{M}$ to a matching $M \subseteq \mc{H}$ using a random greedy algorithm: Iterating through $i=1,\ldots,m$, for each $v_i$ choose, uniformly at random, a hyperedge $T_i \in H$ containing $v_i$ and $k-1$ vertices in $U$ that is vertex-disjoint from $\wt{M}$ and all $T_j$ for $j<i$.

We note that this procedure is sure to be successful. Indeed, before choosing any hyperedge $T_i$, every vertex $v \in V$ satisfies $\deg_H(v,U) \ge \binom{|U|}{\ell-1} (\delta_{\ell,k}^+ + \eps/3) \binom{|U|}{k-\ell} \binom{k-1}{\ell-1}^{-1} \ge (\delta_{\ell,k}^+ + \eps/4) \binom{|U|}{k-1}$. Furthermore, since $\wt{M}$ is contained entirely in $V \setminus U$, none of these hyperedges intersect $\wt{M}$. Thus, there are at least $(\delta_{\ell,k}^+ + \eps/4) \binom{|U|}{k-1}$ choices for $T_i$. Additionally, every hyperedge $T_j$ intersects at most $(k-1)|U|^{k-2}$ possible choices for $T_i$. Since $m (k-1)|U|^{k-2} \le \eps^6|V(\wt{\mc{H}})| (k-1)|U|^{k-2} \le \eps^2|U|^{k-1}$, there are always at least $(\delta_{\ell,k}^+ + \eps/8) \binom{|U|}{k-1}$ choices available for $T_i$.

For the final matching we take $M \coloneqq \wt{M} \cup \{T_1,\ldots,T_m \}$. Clearly, $M$ covers all vertices in $V \setminus U$, proving \cref{itm:M-covers-V-U}. Moreover it covers $(k-1)m \le \eps^2|U|$ vertices in $U$, proving \cref{itm:small-cover-in-U}.

It remains to show that $M$ is $O(1/|V|^{k-1})$-spread. Let $S \subseteq \mc{H}$ be a set of hyperedges. We need to show that $P_S \coloneqq \mb{P}[S \subseteq M] = (O(1/|V|))^{(k-1)|S|}$. First, we may assume that $S$ is a matching. Furthermore, if $S \subseteq M$ then every hyperedge in $S$ is either included in $\wt{M}$ (in which case it has all $k$ vertices in $V \setminus U$) or it is one of the hyperedges $T_1,\ldots,T_m$ (in which case it has exactly one vertex in $V \setminus U$). So we may assume that every hyperedge in $S$ has either one or $k$ vertices in $V \setminus U$. Let $S_k$ be those hyperedges in $S$ with all vertices in $V \setminus U$, and let $S_1 = S \setminus S_k$ be those hyperedges in $S$ with only one vertex in $V \setminus U$. We now have:
\[P_S = \mb{P} \big[ S_k \subseteq \wt{M} \big] \mb{P} \big[ S_1 \subseteq M \setminus \wt{M} | S_k \subseteq \wt{M} \big].\]

By construction, $\wt{M}$ is $O(1/|V|^{k-1})$-spread, so $\mb{P}[ S_k \subseteq \wt{M}] = (O(1/|V|^{k-1}))^{|S_k|}$. Next, we observe that after conditioning on any outcome of $\wt{M}$, it holds that $S_1 \subseteq M \setminus \wt{M}$ only if for every hyperedge $T \in S_1$, the hyperedge chosen to match the (unique) vertex in $T \setminus U$ was $T$. Since every such choice is made uniformly from at least $(\delta_{\ell,k}^+ + \eps/8) \binom{|U|}{k-1} = \Omega(|V|^{k-1})$ possibilities, it follows that $\mb{P}[ S_1 \subseteq M \setminus \wt{M} | S_k \subseteq \wt{M} ] = (O(1/|V|^{k-1}))^{|S_1|}$. Thus $P_S = (O(1/|V|))^{(k-1)|S|}$, as desired.
\end{proof}

We are now in position to complete the proof of \cref{thm:pm}.

\begin{proof}[Proof of \cref{thm:pm}]
We assume, without loss of generality, that $\eps < 1/(100k!)$.

Using \cref{lem:vortex}, let $V(\mc{H}) = V_0 \supseteq V_1 \supseteq \cdots \supseteq V_N = X$ be a random sequence of sets satisfying properties \cref{itm:vortex-set-sizes,itm:vortex-final-set-size,itm:vortex-high-degrees,itm:vortex-prob-guarantee} in \cref{lem:vortex}.

We will inductively construct (random) matchings $\emptyset=M_0 \subseteq M_1 \subseteq \cdots \subseteq M_N$, satisfying the following properties for every $0 \le i \le N$. For notational convenience we set $V_{N+1} = \emptyset$.
\begin{enumerate}
	\item $M_i$ is $O(1/|V_i|^{k-1})$-spread;
	\item $M_i$ covers all vertices in $V(H) \setminus V_i$;
	\item $|V(M_i) \cap V_i| \le 2\eps^2|V_i|$; and
	\item $V(M_i) \cap V_{i+1} = \emptyset$.
\end{enumerate}
	
We begin by taking $M_0 = \emptyset$. Now, suppose that for $0\le i < N$ we have constructed $M_i$ with the properties above. Let $V_i' = V_i \setminus (V(M_i) \cup V_{i+2})$ and let $\mc{H}_i = \mc{H}[V_i']$. Note $V_{i+1}\subseteq V_i'$. Observe that $\delta_\ell(\mc{H}_i) \ge (\delta_{\ell,k}^+ + \eps/3) \binom{|V_i'|}{k-\ell}$. Indeed, for every $S \in \binom{V'}{\ell}$, by \cref{itm:vortex-high-degrees} there holds $\deg_\mc{H}(S,V_i) \ge (\delta_{\ell,k}^+ + \eps/2) \binom{|V_i|}{k-\ell}$. At most $(|V(M_i) \cap V_i|+|V_{i+2}|) \binom{|V_i|}{k-\ell-1} \le 3\eps^2 |V_i|\binom{|V_i|}{k-\ell-1} \le \frac{\eps}{10} \binom{|V_i'|}{k-\ell}$ of these hyperedges are not contained in $V'$. Therefore $\deg_{\mc{H}_i}(S) \ge (\delta_{\ell,k} + \eps/3)\binom{|V_i'|}{k-\ell}$, as desired. By applying \cref{lem:cover-down} to $\mc{H}_i$ with $U = V_i' \cap V_{i+1}$ we obtain an $O_\eps(1/|V_i|^{k-1})$-spread matching $M_i'$ covering all vertices in $V_i' \setminus V_{i+1}$, at most $2\eps^2 |V_{i+1}|$ vertices in $V_{i+1}$, and no vertices in $V_{i+2}$. (We can apply this after checking a similar minimum degree condition with respect to $\deg_{\mc{H}_i}(S,U)$.) By taking $M_{i+1} = M_i \cup M_i'$ we complete the inductive step.

Finally, to obtain a perfect matching, note that if $M_N$ satisfies the properties above then $\delta_\ell(\mc{H}[V(\mc{H}) \setminus V(M_N)]) \ge (\delta_{\ell,k}^+ + \eps/3)\binom{|X|}{k-\ell}$. Furthermore, $|V(\mc{H}) \setminus V(M_N)|$ is a multiple of $k$ since it was obtained from $V(\mc{H})$ by removing a matching. Therefore, by \cref{rmk:large-degree-implies-matching}, there exists a perfect matching $\wt{M} \subseteq \mc{H}[V(\mc{H}) \setminus V(M_N)]$. Take $M = M_N \cup \wt{M}$.

It remains to prove that $M$ is $O_\eps(1/n^{k-1})$-spread. Let $S \subseteq \mc{H}$ be a set of hyperedges, which we may assume are vertex-disjoint. We need to show that $P_S \coloneqq \mb{P} \left[ S \subseteq M \right] = (O_\eps(1/n^{k-1}))^{|S|}$. Let $T_1,\ldots,T_m$ be an enumeration of the hyperedges in $S$. For a vector $\vec{x} \in [N+1]^{m}$, let $P(\vec{x})$ be the probability that for every $j \in [m]$, the hyperedge $T_j$ is in $M_{x_j} \setminus M_{x_j-1}$ if $x_j \le N$, and $T_j \in \wt{M}$ if $x_j = N+1$. We will show that
\begin{equation}\label{eq:P_S}
P(\vec{x}) = \left( \prod_{i=1}^N \left( O_\eps\left( \frac{|V_{i-1}|}{n^k} \right) \right)^{|\{ j\colon x_j=i \}|} \right) \left( \frac{|V_N|}{n} \right)^{k|\{ j\colon x_j=N+1 \}|}.
\end{equation}
This will suffice, since then
\begin{align*}
P_S = \sum_{ \vec{x} \in [N+1]^m } P(\vec{x}) & = \sum_{ \vec{x} \in [N+1]^m } \left( \prod_{i=1}^N \left( O_\eps\left( \frac{|V_{i-1}|}{n^k} \right) \right)^{|\{ j\colon x_j=i \}|} \right) \left( \frac{|V_N|}{n} \right)^{k|\{ j\colon x_j=N+1 \}|}\\
& = \left( \frac{O_\eps(1)}{n^{k-1}} \right)^m \sum_{ \vec{x} \in [N+1]^m } \left( \prod_{i=1}^{N}(\eps^{2i})^{|\{ j\colon x_j=i \}|} \right) \left(\frac{|V_N|^k}{n}\right)^{|\{ j\colon x_j=N+1 \}|}\\
& = \left( O_\eps\left( \frac{1}{n^{k-1}}\right) \right)^m \prod_{j=1}^m \left( \sum_{i=1}^N \eps^{2i} + \frac{n^{k/(k+1)}}{n} \right) = \left( O_\eps\left( \frac{1}{n^{k-1}}\right) \right)^m.
\end{align*}

We now prove \eqref{eq:P_S}. For $1 \le i \le N+1$, let $C_i$ be the event that $\{ T_j\colon x_j = i \} \subseteq \mc{H}[V_{i-1}]$ and let $D_i$ be the event that $\{ T_j\colon x_j = i \} \subseteq M_i \setminus M_{i-1}$ if $i \le N$, and $\{ T_j\colon x_j = i \} \subseteq \wt{M}$ if $i = N+1$. We then have
\[P(\vec{x}) \le \mb{P} \left[ \bigcap_{i=1}^{N+1} C_i \right] \prod_{i=1}^{N+1} \mb{P} \left[ D_i \bigg| \bigcap_{i=1}^{N+1} C_i, D_1 \cap \cdots \cap D_{i-1} \right].\]
By the randomness guarantee in the vortex construction (\cref{itm:vortex-prob-guarantee} in \cref{lem:vortex}), we have:
\[\mb{P} \left[ \bigcap_{i=1}^{N+1} C_i \right] = \prod_{i=1}^{N+1} \left( O \left( \frac{|V_{i-1}|}{n} \right) \right)^{k |\{ j\colon x_j=i \}| }.\]
Next, we note that conditioned on any outcome of $M_{i-1}$ and the vortex, the matching $M_i\setminus M_{i-1}$ is $O_\eps(1/|V_i|^{k-1})$-spread. Thus, for every $i \le N$:
\[\mb{P} \left[ D_i \bigg| \bigcap_{i=1}^{N+1} C_i, D_1 \cap \cdots \cap D_{i-1} \right] = \left(O_\eps\left( \frac{1}{|V_i|^{k-1}} \right)\right)^{|\{ j\colon x_j = i \}|}.\]
Finally, we use the trivial bound $\mb{P} \left[ D_{N+1} \bigg| \bigcap_{i=1}^{N+1} C_i, D_1 \cap \cdots \cap D_{\ell} \right] \le 1$ to obtain:
\[P(\vec{x}) \le \left( \prod_{i=1}^{N+1} \left( O \left( \frac{|V_{i-1}|}{n} \right) \right)^{k |\{ j\colon x_j=i \}| } \right) \left( \prod_{i=1}^N\left(O_\eps\left( \frac{1}{|V_i|^{k-1}} \right)\right)^{|\{ j\colon x_j = i \}|} \right),\]
which implies \eqref{eq:P_S} upon using $|V_{i-1}|=O_\eps(|V_i|)$.
\end{proof}

\section{\texorpdfstring{$K_r$}{Kr}-factors in \texorpdfstring{$r$}{r}-partite super-regular systems via robust-perfect matchings}\label{sec:Kr-short}
We now give the first of two proofs of \cref{thm:super-regular}. A key ingredient is the case where $r=2$, and in particular proving that in a bipartite graph $G=(A,B,E)$ with $|A|=|B|=n$ one can find a spread perfect matching. This tool will again be crucial in the proof of \cref{thm:spanning-trees}. To do so, we consider the following subgraph
of $G$. For each vertex $v$ of $G$, choose a uniform and independent
random set of $C$ neighbors of $v$ (with repetitions). Let $H$
be the graph containing all of the edges chosen by either vertex. 
\begin{lemma}\label{lem:bip-pm}
Let $G$ be $(d,\delta)$-super-regular for $\delta\ll d$. Then, with probability at least $3/4$, for $C$
sufficiently large depending only on $d$, the subgraph $H$ contains
a perfect matching. 
\end{lemma}

\begin{proof}
We will prove that $H$ satisfies Hall's condition with high probability.
In particular, we want to show that there is no subset $T$ of $B$
of size $k$ for which there is a subset $S$ of $A$ of size $k+1$
and the neighborhood of any vertex in $S$ is contained in $T$; and
similarly for $T$ a subset of $A$ of size $k$ and $S$ a subset
of $B$ of size $k+1$. Note that if all vertices in $S$ have their
neighborhood contained in $T$, then all vertices in $B\setminus T$
have their neighborhood contained in $A\setminus S$. Thus, by symmetry,
we only need to consider the case $k\le n/2$, since for $k>n/2$
we have $|A\setminus S|<n/2$. In the following, let $T$ be a subset
of $B$ of size $k$, and $S$ a subset of $A$ of size $k+1$. We
bound the probability that $N_{H}(v)\subseteq T$ for all $v\in S$.
Let $\eta=4\delta^{1/3}$. 

\noindent\textbf{Case 1:} $k\in(\eta n,n/2]$. In this case, by the
assumption that $G$ is $(d,\delta)$-regular, for $\epsilon=\delta^{1/3}/d$,
the number of vertices $v$ with $|N_{G}(v)\cap T|>(k/n+\epsilon)dn$
is at most $\delta^{1/3}n$. Thus there are at least $k-\delta^{1/3}n$ vertices
$v$ in $S$ with $|N_{G}(v)\cap T|\le(k/n+\epsilon)dn$. For each
such $v$, the chance that $N_{H}(v)\subseteq T$ is at most $(k/n+\epsilon^{1/2})^{C}$.
Hence, the probability that $N_{H}(S)\subseteq T$ is at most $(k/n+\epsilon^{1/2})^{C(k-\delta n)}$.
By the union bound, the probability that there exists $S$ and $T$
with $N_{H}(S)\subseteq T$ is at most 
\begin{align*}
\binom{n}{k}\binom{n}{k+1}(k/n+\epsilon^{1/2})^{C(k-\delta^{1/3}n)} & \le\left(\frac{e^{2}n^{2}}{k^{2}}\right)^{k+1}\left(\frac{k}{n}+\epsilon^{1/2}\right)^{C(k-\delta^{1/3}n)}\\
&\le\left(\frac{e^{2}n^{2}}{k^{2}}\right)^{2k}\left(\frac{k}{n}+\epsilon^{1/2}\right)^{Ck/2}
=\left(\frac{e^{4}n^{4}}{k^{4}}\cdot\left(\frac{k}{n}+\epsilon^{1/2}\right)^{C/2}\right)^{k}.
\end{align*}
Note that $\frac{k}{n}+\epsilon^{1/2}\le\min(2/3,2(k/n)/d)$, and hence
\[
\left(\frac{e^{4}n^{4}}{k^{4}}\cdot\left(\frac{k}{n}+\epsilon^{1/2}\right)^{C/2}\right)^{k}\le\left(\frac{e^{4}n^{4}}{k^{4}}\left(\frac{k}{n}+\epsilon^{1/2}\right)^4\cdot(2/3)^{C/2-2}\right)^{k}\le\left(\frac{(2e)^{4}}{d^{4}}\cdot(2/3)^{C/2-2}\right)^{k}<2^{-k},
\]
assuming that $C$ is sufficiently large in $d$. 

\noindent\textbf{Case 2:} $k\le\eta n$. In this case, for each $v\in S$,
$|N_{G}(v)\cap T|\le k\le\eta n$. Hence, the chance that $N_{H}(v)\subseteq T$
is at most $(2k/(dn))^{C}$. Hence, the probability that $N_{H}(S)\subseteq T$
is at most $(2k/(dn))^{Ck}$. By the union bound, the probability that
there exists $S$ and $T$ with $N_{H}(S)\subseteq T$ is at most
\[
\binom{n}{k}\binom{n}{k+1}(2k/(dn))^{Ck}\le\left(\frac{e(2k)^{C-4}}{d^{C}n^{C-4}}\right)^{k}<\left(\frac{e(2k/n)^{C-4}}{d^{C}}\right)^{k}.
\]
Combining the cases, by the union bound, the probability that $H$
does not satisfy Hall's condition is at most 
\[
2\sum_{k\le\eta n}\left(\frac{e(2k/n)^{C-2}}{d^{C}}\right)^{k}+2^{-\eta n+1}=o(1).\qedhere
\]
\end{proof}

Using \cref{lem:bip-pm}, we can give a general procedure for finding spread matchings in super-regular bipartite graphs. 
\begin{theorem}\label{thm:bip-spread}
Let $d>0$ and $\delta\ll d$. Let $G$ be a $(d,\delta)$-super-regular bipartite graph with parts of size $n$. There exists a distribution $\mu$ on perfect matchings in $G$ which is
$O_{d}(1/n)$-spread. 
\end{theorem}
\begin{proof}
From $G$ pick a subgraph $H$ as in \cref{lem:bip-pm}, which
has a perfect matching with probability at least $3/4$. Condition on this event and pick and output an arbitrary perfect matching $W$ of $H$. This induces a distribution $\mu$
on perfect matchings of $G$. We show that $\mu$ is $O_{d}(1/n)$-spread.
Indeed, given any subset $S$ of edges of $G$, if $S$ is not a matching,
then $\mu(W\supseteq S)=0$. If $S$ is a matching, $W$ can only
contain $S$ if for each edge $e=\{x,y\}\in S$, either $x$ or $y$
picks the other vertex as one of the $C$ neighbors, which happens
with probability at most $2C/n$. Furthermore, the above events are
independent across different edges of the matching $S$. Hence, $\mu(W\supseteq S)\le(4/3)(2C/n)^{|S|}$.
Thus, $\mu$ is $(4C/n)$-spread. 
\end{proof}

One can now prove \cref{thm:super-regular} via an inductive argument on $r$.

\begin{lemma}\label{lem:matching-reg}
Let $G=(V_1,V_2,V_3,E)$ be a $(d,\delta)$-super-regular tripartite graph with $|V_1|=|V_2|=|V_3|$. Assume that for each edge $\{v_2,v_3\}\in E(G)$, there are at least $(d^2-\delta^{1/2})n$ vertices $v_1\in V_1$ which are adjacent to both $v_2$ and $v_3$. Let $\mu_1$ be the distribution on perfect matchings $M_1$ between $V_2$ and $V_3$ given by \cref{thm:bip-spread}, which is $p$-spread for some $p\le C/n$. Construct a graph $\Gamma_{M_1}$ where the vertices are the edges in $M_1$ and vertices in $V_1$, and an edge $e=\{v_2,v_3\}$ of $M_1$ is connected to a vertex $v_1$ of $V_1$ if and only if $v_1$ is adjacent to both $v_2$ and $v_3$. There exists $d',c'$ depending only on $d$ and $C$ such that for $\delta$ sufficiently small in terms of $d,C$, with probability at least $1-\exp(-c'n)$, we have that $\Gamma_{M_1}$ contains a $(d',16C/\log(1/\delta)^{1/4})$-super-regular subgraph. 
\end{lemma}

\begin{proof}
Let $\eta=(\log\delta^{-1})^{-1/2}$. We first prove that $\Gamma_{M_{1}}$
is $4\eta^{1/2}$-regular with high probability. Indeed, for each
subset $S_{1}$ of $V_{1}$ of size $\rho n$ and $\rho\ge4\eta^{1/2}$,
we say that $e\in E(G_{1})$ is bad if the number of common neighbors
of the endpoints of $e$ in $S_{1}$ is not in $(d^{2}\pm\eta)\rho n$,
and say that $v_{i}\in V_{i}$ is bad for $i\in\{2,3\}$ if the number
of neighbors of $v_{i}$ in $S_{1}$ is not in $(d\pm\eta)\rho n$.
The number of bad $v_{i}$ is at most $\delta n$. For each $v_{i}$
which is not bad, the number of bad edges $e$ adjacent to $v_{i}$ is at
most $\delta n$. We say that a perfect matching $M_{1}$ of $G[V_2,V_3]$ is bad if it contains
at least $\eta n$ bad edges which are adjacent only to good vertices. The number
of bad perfect matchings is at most 
\[
\binom{n}{\eta n}(\delta n)^{\eta n}\cdot n^{n-\eta n}.
\]
For each such perfect matching, the probability (under $\mu_{1}$)
that it is realized is at most $(C/n)^{n}$. Hence, the probability that there are at least $\eta n$
bad edges adjacent to good vertices selected in $M_{1}$ is at most
\begin{align*}
\binom{n}{\eta n}(\delta n)^{\eta n}\cdot(n)^{n-\eta n}\cdot(C/n)^{n} & \le\exp(-(\log\delta^{-1})\eta n)C^{n}\exp(\eta n\log(\delta^{-1})/2)\\
 & \le\exp(-(\log\delta^{-1})^{1/2}n/4).
\end{align*}
Note that if there are at most $\eta n$ bad edges adjacent to good
vertices in $M_{1}$, then the number of vertices in $E(M_{1})$ whose
number of edges to $S_{1}$ is not $(d^{2}\pm\eta)\rho n$ is at
most $(2\delta+\eta)n$. In that case, for any subset $T$ of
$E(M_{1})$ of size at least $\rho n$, the number of edges between $T$ and $S_{1}$ is $(d^{2}\pm2(\eta+2\delta)\rho^{-1})\rho n|T|$.
Hence, by the union bound over $S_{1}$, we obtain that $\Gamma_{M_{1}}$
is $4\eta^{1/2}$-regular with probability at least 
\[
2^{n}\exp(-(\log\delta^{-1})^{1/2}n/4)<2^{-n}.
\]

By assumption, the minimum degree of each $e\in M_1$ in $\Gamma_{M_1}$ is at least $(d^2-\delta^{1/2})n$. For each vertex $v_1\in V_1$, let $E(v_1)$ be the set of edges in $E(G_1)$ whose endpoints are both adjacent to $v_1$. Then each vertex in $V_2\cup V_3$ is adjacent to at least $(d^2-\delta^{1/2})n$ edges in $E(v_1)$. By the remark following \cref{thm:bip-spread}, the probability that the degree of $v_1$ in $\Gamma_{M_1}$ is at most $(d^2/(2e^2))^{C}n$ is at most $\exp(-(d^2/(4e^2))^Cn)$. Hence, by the union bound, with probability at least $1-n\exp(-(d^2/(4e^2))^Cn)$, the minimum degree of $\Gamma_{M_1}$ is at least $(d^2/(2e^2))^Cn$. The conclusion of the lemma then follows from \cref{lem:super-boost,lem:plus-to-super-regular}.
\end{proof}

\begin{proof}[Proof of \cref{thm:super-regular}]
We prove the result by induction on $r$. The case $r=2$ is shown in \cref{thm:bip-spread}. 

Let $G_{r-1,r}$ be the graph induced on vertex sets $V_{r-1}$ and $V_r$. By \cref{thm:bip-spread}, there is a distribution $\mu_1$ on perfect matchings $M_1$ of $G_{r-1,r}$ which is $(C_d/n)$-spread for $C_d$ depending only on $d$. By \cref{lem:matching-reg}, if for $i<r-1$, we construct the graph $\Gamma_{i,M_1}$ which has as vertices edges $e=\{v_{r-1},v_r\}$ in $M_1$ and vertices $v_i$ in $V_i$ for which $v_i$ is adjacent to both $v_{r-1}$ and $v_r$, then with probability at least $1-\exp(-c'n)$, $\Gamma_{i,M_1}$ has a subgraph which is $(d',16C/(\log \delta^{-1})^{1/4})$-super-regular. By the union bound, with high probability, this property holds for all $i<r-1$. Now we have an $(r-1)$-partite graph $G'$ where $G'[V_i,V_j]=G[V_i,V_j]$ for $i,j<r-1$ and $G'[V_i,V_{r-1}]=\Gamma_{i,M_1}$ for $i<r-1$, for which each pair of parts is $(d^{\prime +},16C/(\log \delta^{-1})^{1/4})$-super-regular. By the inductive hypothesis, we have a distribution $\mu$ on perfect matchings $\wt{M}$ of $G'$ which is $O_{d}(1/n^{r-2})$-spread. The perfect matching $\wt{M}$ the corresponds to a perfect matching $M$ of $G$. We now verify that $M$ is $O_{d}(1/n^{r-1})$-spread.

Fix a subset $S$ of hyperedges. As before, we can assume that $S$
is a matching. Let $S_{1}$ be the matching of $G_{r-1,r}$ induced by
$S$. If $M\supseteq S$, then $M_{1}\supseteq S_{1}$, which holds
with probability at most $(C/n)^{|S|}$ for some $C$ depending
only on $d$. Furthermore, conditioned on a consistent realization
of $M_{1}$, we need the matching $\wt{M}$ to contain a corresponding
set of edges of $E(G')$ of size $|S|$, which holds with
probability at most $(O_d(1/n^{r-2}))^{|S|}$. Hence, the probability that $M$
contains $S$ is at most $(C'/n^{r-1})^{|S|}$. Thus, the distribution $\mu$ is $O_{d}(1/n^{r-1})$-spread. 
\end{proof}

\section{\texorpdfstring{$K_r$}{Kr}-factors in \texorpdfstring{$r$}{r}-partite super-regular systems via iterative absorption}\label{sec:Kr-IA}

In order to prove \cref{thm:super-regular}, we will require a regularity boosting lemma for the clique complex above a set of super-regular pairs. Regularity boosting plays a crucial role in the applications of iterative absorption; see e.g. \cite[Lemma~4.2]{BGKLMO20}. While our regularity boost generally follows a similar strategy of using local gadgets to adjust the initial uniform weighting on the hypergraph to give a good fractional weighting, generally partite instances requires a great deal more care (see e.g.~work of Montgomery \cite{Mon17}); however in our case a substantially simpler proof suffices.

\begin{lemma}[Fractional matching]\label{lem:fractional-matching}
Fix $r\ge 2$ and suppose $1/n\ll\eps\ll d,1/r\le 2/3$. Let $G=(V,E)$ be an $r$-partite graph on $V=\bigcup_{i=1}^rA_i$ where $|A_i|=n$ for all $i\in[r]$. Suppose $G[A_i,A_j]$ is $(d_{i,j},\eps)$-super-regular with $d_{i,j}\ge d$ for all $i\neq j$. Let $\mc{H}$ be the $r$-uniform hypergraph where edges in $\mc{H}$ correspond to $r$-partite cliques of $G$.
Then there exists a weighting $\omega\colon\mc{H}\to[0,1]$ such that for all $v\in \bigcup_{i=1}^rA_i$ we have that 
\[\sum_{\substack{H\ni v\\H\in\mc{H}}}\omega(H) = \frac{1}{2} n^{r-1}\prod_{1\le i<j\le r}d_{i,j}.\]
\end{lemma}

By sampling cliques according to $\omega$ and applying the Chernoff bound we have the following immediate corollary. 
\begin{corollary}\label{cor:fractional-matching-sampled}
Fix $r\ge 2$ and suppose $1/n\ll\eps\ll d,1/r\le 2/3$. Let $G=(V,E)$ be a $r$-partite graph on $V=\bigcup_{i=1}^rA_i$ where $|A_i|=n$ for all $i\in[r]$. Suppose $G[A_i,A_j]$ is $(d_{i,j},\eps)$-super-regular with $d_{i,j}\ge d$ for all $i\neq j$. Let $\mc{H}$ be the $r$-uniform hypergraph where edges in $\mc{H}$ correspond to $r$-partite cliques of $G$.
Then there exists a weighting $\omega\colon\mc{H}\to\{0,1\}$ such that for all $v\in\bigcup_{i=1}^{r}A_i$ we have that
\[\sum_{\substack{H\ni v\\H\in\mc{H}}}\omega(H) = \frac{1}{2} n^{r-1}\prod_{1\le i<j\le r}d_{i,j}\pm n^{r-4/3}.\]
\end{corollary}

We now prove \cref{lem:fractional-matching}. 
\begin{proof}[{Proof of \cref{lem:fractional-matching}}]
By applying \cref{lem:inner-regularity}, we may assume that each pair $(A_i,A_j)$ is $((d_{i,j}-4\eps)^+, \eps^{1/3})$-super-regular and that the complement is $((1-d_{i,j}+4\eps)^+, \eps^{1/3})$-super-regular. The weight function $\omega$ will be a perturbation of the function which is uniformly $1/2$ on $\mc{H}$.

Fix a vertex $v\in A_i$. By the degree lower bounds for $v$ to $A_j$ with $j\neq i$ and by \cref{lem:counting-lemma} applied to $|N(v)\cap A_j|$ for $j\neq i$ to count copies of $K_{r-1}$, we have that
\[\deg_\mc{H}(v)=\sum_{\substack{H\in \mc{H}\\H\ni v}} 1 = n^{r-1}\prod_{1\le i<j\le r} d_{i,j} \pm \eps^{1/4}n^{r-1}.\]
Note that this implies that $|E(\mc{H})| = n^r\prod_{1\le i<j\le r} d_{i,j} \pm \eps^{1/4}n^{r}$. We define the \textit{defect} of a vertex $v$ as 
\[D_v = \deg_\mc{H}(v) - \frac{|E(\mc{H})|}{n}.\]
Note that $|D_v|\le 2\eps^{1/4}n^{r-1}$.

We will now define weight-shifting gadgets. For a pair of vertices $v_1,v_2 \in A_r$, let $\mc{R}_{v_1,v_2}$ be the set of $(2r-1)$-tuples of distinct vertices $(v_1',\ldots,v_{2r-1}')$ where $v_{2r-1}'\in A_r$, $v_i',v_{i+r-1}'\in A_{i}$ for $i\in [r-1]$, and $(v_1, v_1',\ldots, v_{r-1}')$, $(v_{2r-1}',v_1',\ldots, v_{r-1}')$, $(v_2, v_{r}',\ldots, v_{2r-2}')$, and $(v_{2r-1}', v_{r}',\ldots, v_{2r-2}')$ are all in $\mc{H}$. By applying the counting lemma (\cref{lem:counting-lemma}) we have that 
\[|\mc{R}_{v_1,v_2}| = n^{2r-1}\prod_{1\le i<j\le r-1}d_{i,j}^2\prod_{1\le i\le r-1}d_{i,r}^4\pm\eps^{1/4}n^{2r-1}.\]
For each $R\in\mc{R}_{v_1,v_2}$, define $f_{v_1,v_2,R}^{r}\colon\mc{H}\to\mb{R}$ by assigning $0$ to everything outside the four distinguished $r$-cliques of $R$, assigning $(D_{v_1}-D_{v_2})/(2n|\mc{R}_{v_1,v_2}|)$ to $(v_{2r-1}',v_1',\ldots, v_{r-1}')$ and $(v_2, v_{r}',\ldots, v_{2r-2}')$, and assigning $-(D_{v_1}-D_{v_2})/(2n|\mc{R}_{v_1,v_2}|)$ to $(v_1, v_1',\ldots, v_{r-1}')$ and $(v_{2r-1}', v_{r}',\ldots, v_{2r-2}')$. Define $f_{v_1,v_2}^r\colon\mc{H}\to\mb{R}$ as the sum of all $f_{v_1,v_2,R}^r$ over $R\in\mc{R}_{v_1,v_2}$. One can define the analogous construction for each pair of vertices in the same part for different $j\in[r]$, not just $A_r$. The modified function $\omega$ will simply be 
\[\omega(H) = \frac{2^{-1}n^{r-1}\prod_{1\le i<j\le r}d_{i,j}}{|E(\mc{H})|n^{-1}}\bigg(1+\sum_{\substack{1\le i\le r\\v_1,v_2\in A_i}} f_{v_1,v_2}^{i}(H)\bigg).\]

Notice that the functions $f^{i}_{v_1,v_2}$ are mean zero upon averaging over all $H\in\mc{H}$ by construction and therefore to prove the desired result it suffices to prove that for each pair of vertices $v_1,v_2\in A_i$ there holds
\[
\sum_{\substack{H\in\mc{H}\\H\ni v_1}} \omega(H) = \sum_{\substack{H\in\mc{H}\\H\ni v_2}} \omega(H)
\]
and that $\omega(H)\in [0,1]$ for every $H \in \mc{H}$. For the first claim notice that if $v\in A_\ell$ then by construction
\begin{align*}
\bigg(\frac{2^{-1}n^{r-1}\prod_{1\le i<j\le r}d_{i,j}}{|E(\mc{H})|n^{-1}}\bigg)^{-1}\sum_{\substack{H\in\mc{H}\\H\ni v}} \omega(H) &=  \sum_{H\ni v}\bigg(1+\sum_{\substack{1\le i\le r\\v_1,v_2\in A_i}} f_{v_1,v_2}^{i}(H)\bigg)\\
&= \sum_{H\ni v}\bigg(1+2\sum_{v'\in A_\ell} f_{v,v'}^{\ell}(H)\bigg) = \sum_{H\ni v}1 + \sum_{v'\in A_\ell}\frac{(D_{v'}-D_v)}{n} \\
&= \deg_\mc{H}(v) - D_v + \frac{1}{n} \sum_{v'\in A_{\ell}}D_v'= \frac{|E(\mc{H})|}{n}.
\end{align*}
Thus it remains to prove that $\omega(H)\in[0,1]$ for every $H \in \mc{H}$. We saw above that for every $v_1,v_2 \in A_i$ there holds $|\mc{R}_{v_1,v_2}^i| \geq d^{O(1)}n^{2r-1}$. Furthermore each hyperedge is given non-zero weight by at most $4rn^{r+1}$ gadgets $f_{v_1,v_2}^i$. Finally, there holds
\[
|f_{v_1,v_2}^i(H)| \leq (|D_{v_1}|+|D_{v_2}|)/(2n|\mc{R}_{v_1,v_2}^i|) \leq 4\varepsilon^{1/4}n^{r-1} / (2n^{2r} d^{O(1)})
\]
for every $H \in \mc{H}$, and every $v_1,v_2 \in A_i$. This implies that 
\begin{align*}
\bigg|\sum_{\substack{1\le i\le r\\v_1,v_2\in A_i}} f_{v_1,v_2}^{i}(H)\bigg|&\le 4rn^{r+1} \times 4\varepsilon^{1/4} n^{r-1} / (2n^{2r}d^{O(1)}) = \varepsilon^{1/4}d^{-O(1)}.
\end{align*}
As $\eps\ll d$ and $\frac{2^{-1}\prod_{1\le i<j\le r}d_{i,j} n^{r-1}}{|E(\mc{H})|n^{-1}}\in [1/3,2/3]$ by the counting lemma (\cref{lem:counting-lemma}), the desired result follows immediately. 
\end{proof}

We are now in position to prove \cref{thm:super-regular}. Given \cref{cor:fractional-matching-sampled} and \cref{lem:blowup}, the proof (via iterative absorption) is analogous to that of \cref{thm:pm}. First we prove a vortex lemma similar to \cref{lem:vortex}, but starting with the setup in \cref{thm:super-regular}.

\begin{lemma}[Vortex]\label{lem:super-vortex}
Fix $r\ge 2$ and suppose $1/n\ll\eta\ll\eps\ll d$. Let $G=(V,E)$ be an $r$-partite graph on partition $V=\bigcup_{j=1}^rA_j$ where $|A_j|=n$ for all $j\in[r]$. Suppose $G[A_j,A_k]$ is $(d_{j,k},\eps)$-super-regular with $d_{j,k}\ge d$ for all $j\neq k$. Then there exists a distribution on set sequences $V(G) = V_0 \supseteq V_1 \supseteq \cdots \supseteq V_N = X$ with the following properties:
\begin{enumerate}[{\bfseries{V\arabic{enumi}}}]
    \item\label{itm:super-vortex-balanced} For every $0\le i\le N$ there holds $|V_i\cap A_j|=|V_i|/r$ for all $j\in[r]$;
	\item\label{itm:super-vortex-set-sizes} For every $0\le i<N$ there holds $|V_{i+1}| = (1\pm\eta/N^2)\eta|V_i|$;
	\item\label{itm:super-vortex-final-set-size} $|X| \in [n^{1/(r+2)},n^{1/(r+1)}]$;
	\item\label{itm:super-vortex-super-regular} For every $0\le i\le N$ and $j\neq k$ there holds $G[V_i\cap A_j,V_i\cap A_k]$ is $(d_{j,k}^{(i)},\eps^{1/20})$-super-regular for some $d_{j,k}^{(i)} = d_{j,k}\pm2\eps$;
	\item\label{itm:super-vortex-high-degrees} For every $0\le i<N$, $j\neq k$, and $v\in V_i\cap A_j$ there holds $\deg_G(v,V_{i+1}\cap A_k)\ge(d_{i,j}-2\eps)|V_{i+1}|/r$;
    \item\label{itm:super-vortex-prob-guarantee} For every vertex set $\{v_1,\ldots,v_m\} \subseteq V(H)$ and every vector $\vec{x} \in \{0,\ldots,N\}^m$ there holds
	\[\mb{P}\bigg[ \bigwedge_{i=1}^m (v_i \in V_{x_i}) \bigg] \le \prod_{i=1}^m \frac{2|V_{x_i}|}{n}.\]
	\end{enumerate}
\end{lemma}

\begin{proof}
First, consider the distribution of set sequences $V(G)=U_0\supseteq\cdots\supseteq U_N$ obtained as follows: Set $U_0=V(G)$. For as long as $|U_i|>n^{1/(r+1)}$, for each $j\in[r]$ let $U_{i+1}\cap V_j$ be a uniformly random subset of $U_i\cap V_r$ of size exactly $\lceil\eta|U_{i+1}|/r\rceil$. Observe that \ref{itm:super-vortex-balanced}, \ref{itm:super-vortex-final-set-size}, and \ref{itm:super-vortex-final-set-size} hold by definition. Chernoff's inequality for hypergeometric distributions (\cref{lem:chernoff}) and a union bound imply that \ref{itm:super-vortex-high-degrees} holds w.h.p.

Showing that \ref{itm:super-vortex-super-regular} holds w.h.p.\ requires a little more work. Let $i \in [N]$ and $j,k \in [r]$ be distinct. We use regularity and Chernoff's inequality to show that with all but exponentially small probability there holds $d(U_i \cap A_j, U_i \cap A_k) = d_{j,k} \pm 2\varepsilon$. In order to show that $G[U_i \cap A_j,U_i \cap A_k]$ is $\varepsilon^{1/20}$-regular we use the well-known equivalence between the sums of codegrees and regularity \cite{CGW89}. That is, we first use McDiarmid's inequality (\cref{lem:McDiarmid}) to prove that with all but exponentially small probability there holds
\[
\sum_{u,v \in U_i \cap A_j} |\{w \in U_i \cap A_k : uw,vw \in E(G)\}| = (1\pm o(1)) \left( \frac{|U_i|}{|V|} \right)^3 \sum_{u,v \in A_j} |\{w \in A_k : uw,vw \in E(G)\}|.
\]
Since $(A_j,A_k)$ is $(d_{j,k},\eps)$-regular the sum on the right is equal to $n^3d_{j,k}^2 \pm \varepsilon^{1/4}n^3$. Hence the sum on the left is equal to $|U_i \cap A_j|^3 d_{j,k}^2 \pm \varepsilon^{1/5}|U_i|^3$. But this, in turn, implies that $G[U_i \cap A_j,U_i \cap A_k]$ is $\varepsilon^{1/20}$-regular. In order to show super-regularity it remains to verify the minimum degree condition; this follows immediately from \ref{itm:super-vortex-high-degrees}.

%Let $\mc{E}$ be the event that properties \cref{itm:super-vortex-balanced,itm:super-vortex-set-sizes,itm:super-vortex-final-set-size,itm:super-vortex-super-regular,itm:super-vortex-high-degrees} hold for this set sequence. McDiarmid's inequality and a union bound imply that $\mc{E}$ holds w.h.p.

Next, let $\{v_1,\ldots,v_m\}\subseteq V(G)$ and $\vec{x}\in\{0,\ldots,N\}^m$. Clearly
\[
\mb{P}\bigg[\bigwedge_{i=1}^m(v_i\in U_{x_i})\bigg]\le\prod_{i=1}^m \frac{3|U_{x_i}|}{2n},
\]
since each next set $U_{i+1}$ of the binomial set process can be coupled inside a $\eta(1+1/N^2)$-binomial random subset, say, of the current set $U_i$ (and multiplying over at most $N$ steps).

Let $\mc{E}$ be the event that properties \cref{itm:super-vortex-balanced,itm:super-vortex-set-sizes,itm:super-vortex-final-set-size,itm:super-vortex-super-regular,itm:super-vortex-high-degrees} hold for $U_0,\ldots,U_N$. Let $V_0\supseteq\cdots\supseteq V_N$ be the distribution obtained by conditioning $U_0\supseteq\cdots\supseteq U_N$ on the occurrence of $\mc{E}$. By definition, $V_0\supseteq\cdots\supseteq V_N$ satisfies properties \cref{itm:super-vortex-balanced,itm:super-vortex-set-sizes,itm:super-vortex-final-set-size,itm:super-vortex-super-regular,itm:super-vortex-high-degrees}. Furthermore, for every nonempty $\{v_1,\ldots,v_m\}\subseteq V(G)$ and $\vec{x}\in\{0,\ldots,N\}^m$:
\[\mb{P}\bigg[\bigwedge_{i=1}^m(v_i\in V_{x_i})\bigg] = \mb{P}\bigg[\bigwedge_{i=1}^m(v_i\in U_{x_i})\bigg|\mc{E}\bigg]\le\frac{\mb{P}\big[\bigwedge_{i=1}^m(v_i\in U_{x_i})\big]}{\mb{P}[\mc{E}]}\le \prod_{i=1}^m\frac{2|V_{x_i}|}{n},\]
as desired. The last inequality comes from applying \cref{itm:super-vortex-set-sizes} iteratively at most $N$ times.
\end{proof}

Next we prove a cover-down lemma similar to \cref{lem:cover-down}.
\begin{lemma}[Cover-down lemma]\label{lem:super-cover-down}
Fix $r\ge 2$ and suppose $1/m\ll\eta\ll\eps\ll d$. Let $G=(V,E)$ be a $r$-partite graph on partition $V=\bigcup_{j=1}^rA_j$ where $|A_j|=m$ for all $j\in[r]$. Suppose $G[A_j,A_k]$ is $(d_{j,k},\eps)$-super-regular with $1/2\ge d_{j,k}\ge d$ for all $j\neq k$. Let $U\subseteq V(G)$ satisfy $|U\cap A_j|=|U|/r$ for all $j\in[r]$ and $|U|=(1\pm\eta)\eta(rm)$. Suppose that for all $j\neq k$ and $v\in A_j$ we have $\deg_G(v,U\cap A_k)\ge(d_{i,j}-\eps)|U|/r$, and for $j\neq k$ we have that $G[U\cap A_j,U\cap A_k]$ is $(d_{j,k}',\eps)$-super-regular for some $1/2\ge d_{j,k}'\ge d$.
Then there exists a $C_\eta/m^{r-1}$-spread distribution on partial $K_r$-factors $M$ of $G$ that satisfies:
\begin{enumerate}[{\bfseries{C\arabic{enumi}}}]
	\item\label{itm:super-M-covers-V-U} $M$ covers every vertex in $V(G)\setminus U$;
	\item\label{itm:super-small-cover-in-U} $M$ covers at most $\eta|U|$ vertices in $U$.
\end{enumerate}
\end{lemma}
\begin{proof}
We first construct the random partial $K_r$-factor $M$.
	
\textit{Step 1: Finding a regular clique system.} We will find a regular collection of cliques of $G[V \setminus U]$ by applying \cref{cor:fractional-matching-sampled}. Let $V'\coloneqq V \setminus U$ and $G' \coloneqq G[V']$. Observe that $G'[V'\cap A_j,V'\cap A_k]$ is $(d_{j,k}',2\eps)$-super-regular for some $d_{j,k}'\in[d_{j,k}\pm2\eps]$ since we removed a small (with respect to $\eps$) fraction. Since $|V'\cap A_j|\ge m/2$ is sufficiently large with respect to $\eps$ and the part sizes are equal by the given conditions, we can apply \cref{cor:fractional-matching-sampled} to find a set $\wt{\mc{H}}$ of $r$-cliques of $G'$ so that every $v\in V'$ is contained in
\[\frac{1}{2}n^{r-1}\prod_{1\le i<j\le r}d_{i,j}'\pm n^{r-4/3}\]
many cliques.

\textit{Step 2: Finding a spread approximate matching.} We apply \cref{lem:nibble} to $\wt{\mc{H}}$ to find an $O_{\eta,r}(1/m^{k-1})$-spread matching $\wt{M}\subseteq\wt{\mc{H}}$ covering all but at most $\eta^6|V'|$ vertices of $V'$.

\textit{Step 3: Covering remaining vertices in $V \setminus U$.} Conditioning on $\wt{M}$, let $v_1,\ldots,v_t$ be an enumeration of the uncovered vertices $V \setminus (V(\wt{M}) \cup U)$, noting $t\le\eta^6|V'|$. Note that there are an equal amount in each part $A_j$ for $j\in[r]$. We extend $\wt{M}$ to a partial $K_r$-factor $M$ covering all of $V\setminus U$ (and some of $U$) using a random greedy algorithm: Iterating through $i=1,\ldots,t$, for each $v_i$ choose, uniformly at random, an $r$-clique $T_i$ of $G$ containing $v_i$ and $k-1$ vertices in $U$ that is vertex-disjoint from $\wt{M}$ and all $T_j$ for $j<i$.

We note that this procedure is sure to be successful. Indeed, before choosing any hyperedge $T_i$, every vertex $v\in A_j$ satisfies $\deg_G(v,U\cap A_k)\ge (d_{j,k}-\eps)|U|/r$. Furthermore, since $\wt{M}$ is contained entirely in $V \setminus U$, none of these hyperedges intersect $\wt{M}$. Thus, for some $T_i$ with $v_i\in A_j$, there are at least
\[\prod_{k\neq j}(|U|/r)\bigg(\prod_{k\neq j}(d_{j,k}-\eps)\prod_{\substack{1\le k_1<k_2\le r\\k_1,k_2\neq j}}d_{k_1,k_2}-C\eps\bigg)\]
choices for $T_i$ by \cref{lem:counting-lemma} and the super-regularity of the pairs $(U\cap A_{k_1},U\cap A_{k_2})$.
Additionally, every other hyperedge $T_j$ intersects at most $(r-1)|U|^{r-2}$ possible choices for $T_i$. Since $t(r-1)|U|^{r-2}\le (\eta^6m)(r-1)|U|^{r-2} \le\eta|U|^{r-1}$, there are always at least say
\[(|U|/(2r))^{r-1}\prod_{1\le j<k\le r}d_{j,k}\]
choices available for $T_i$ regardless of the prior choices.

For the final matching we take $M\coloneqq\wt{M}\cup\{T_1,\ldots,T_t\}$. Clearly, $M$ covers all vertices in $V\setminus U$, proving \cref{itm:M-covers-V-U}. Moreover it covers $(r-1)t\le\eta|U|$ vertices in $U$, proving \cref{itm:small-cover-in-U}.

It remains to show that $M$ is $O_\eta(1/m^{r-1})$-spread. Let $S$ be a set of $r$-cliques of $G$. We need to show that $P_S \coloneqq \mb{P}[ S \subseteq M ] = (O(1/m))^{(r-1)|S|}$. First, we may assume that $S$ is a partial $K_r$-factor. Furthermore, if $S \subseteq M$ then every $r$-clique in $S$ is either included in $\wt{M}$ (in which case it has all $k$ vertices in $V\setminus U$) or it is one of the hyperedges $T_1,\ldots,T_t$ (in which case it has exactly one vertex in $V\setminus U$). So we may assume that every hyperedge in $S$ has either one or $k$ vertices in $V\setminus U$. Let $S_k$ be those hyperedges in $S$ with all vertices in $V\setminus U$, and let $S_1 = S \setminus S_k$ be those hyperedges in $S$ with only one vertex in $V\setminus U$. We now have:
\[P_S=\mb{P}\big[S_k\subseteq\wt{M}\big]\mb{P}\big[S_1\subseteq M\setminus \wt{M}|S_k\subseteq\wt{M}\big].\]

By construction, $\wt{M}$ is $O_\eta(1/m^{k-1})$-spread, so $\mb{P}[S_k \subseteq\wt{M}] = (O(1/m^{k-1}))^{|S_k|}$. Next, we observe that after conditioning on any outcome of $\wt{M}$, it holds that $S_1 \subseteq M \setminus \wt{M}$ only if for every hyperedge $T\in S_1$, the hyperedge chosen to match the (unique) vertex in $T\setminus U$ was $T$. Since every such choice is made uniformly from at least $(|U|/(2r))^{r-1}\prod_{1\le j<k\le r}d_{j,k}=\Omega_\eta(m^{k-1})$ possibilities, it follows that $\mb{P}[ S_1 \subseteq M \setminus \wt{M} | S_k \subseteq \wt{M} ] = (O(1/m^{k-1}))^{|S_1|}$. Thus $P_S = (O(1/m))^{(k-1)|S|}$, as desired.
\end{proof}

Now we prove \cref{thm:super-regular}.
\begin{proof}[Proof of \cref{thm:super-regular}]
First, note that we may assume $d_{j,k}\le 1/3$ for all $j\neq k$ by considering a random $1/3$-sample of the edges of the graph $G$ and adjusting parameters appropriately. Additionally, we then choose parameters so that
\[1/n\ll\eta\ll\eps\ll d\le 1/3,\]
treating our choice of $\eta$ ultimately as a function of $d,\eps$. By renaming parameters we may assume that all the pairs in $G$ are $\varepsilon^{20}$-super-regular.

Using \cref{lem:super-vortex}, let $V(G) = V_0 \supseteq V_1 \supseteq \cdots \supseteq V_N = X$ be a random sequence of sets satisfying properties \cref{itm:super-vortex-balanced,itm:super-vortex-set-sizes,itm:super-vortex-final-set-size,itm:super-vortex-super-regular,itm:super-vortex-high-degrees,itm:super-vortex-prob-guarantee} in \cref{lem:super-vortex}.

We will inductively construct (random) partial $K_r$-factors $M_0 \subseteq M_1 \subseteq \cdots \subseteq M_N$, satisfying the following properties for every $0 \le i \le N$. For notational convenience we set $V_{N+1} = \emptyset$.
\begin{enumerate}
	\item $M_i$ is $O(1/|V_i|^{k-1})$-spread;
	\item $M_i$ covers all vertices in $V(G) \setminus V_i$;
	\item $|V(M_i) \cap V_i| \le 2\eta|V_i|$; and
	\item $V(M_i) \cap V_{i+1} = \emptyset$.
\end{enumerate}
	
We begin by taking $M_0 = \emptyset$. Now, suppose that for $i < N$ we have constructed $M_i$ with the properties above. Let $V_i' = V_i \setminus (V(M_i) \cup V_{i+2})$ and let $G_i = G[V_i']$. Note $V_{i+1}\subseteq V_i'$. Observe that $G_i$ and $U=V_{i+1}\setminus V_{i+2}\subseteq V(G_i)$ satisfy the hypotheses of \cref{lem:super-cover-down} for $m=|V_i'|$ and slightly modified parameters $d_{j,k},d_{j,k}',\eps$, using \cref{itm:super-vortex-balanced,itm:super-vortex-set-sizes} for the set sizes, \cref{itm:super-vortex-high-degrees} for the degree condition, and \cref{itm:super-vortex-super-regular} for $i,i+1$ for the super-regularity conditions (plus $V_{i+1}\subseteq V_i'$ as well as the fact that removing $V_{i+2}$ is removing a negligible $\eta$-fraction of $V_{i+1}$ and thus does not affect these conditions severely). By applying \cref{lem:super-cover-down} to $G_i,U=V_{i+1}\setminus V_{i+2}$ in this situation we obtain an $O_\eta(1/|V_i|^{k-1})$-spread partial $K_r$-factor $M_i'$ covering all vertices in $V_i'\setminus V_{i+1}$, at most $2\eta|V_{i+1}|$ vertices in $V_{i+1}$, and no vertices in $V_{i+2}$. By taking $M_{i+1} = M_i\cup M_i'$ we complete the inductive step.

Finally, to obtain a perfect matching, note that if $M_N$ satisfies the properties above then \cref{itm:super-vortex-super-regular} for $i=N$ and the fact that we only delete a small fraction of $V_N$ means that we have a super-regular remainder, to which \cref{lem:blowup} applies (with $H$ being a disjoint collection of $r$-cliques and $R=K_r$ underlying the $r$-partite structure of $G_1=G[V(G)\setminus V(M_N)]$ induced by $A_1,\ldots,A_r$). That is, we have a $K_r$-factor $\wt{M}$ of $G[V(G)\setminus V(M_N)]$. Take $M = M_N \cup \wt{M}$.

It remains to prove that $M$ is $O(1/n^{r-1})$-spread. Let $S$ be a set of $r$-cliques. We need to show that $P_S\coloneqq\mb{P}[S\subseteq M] = (O_\eta(1/n^{r-1}))^{|S|}$. Let $T_1,\ldots,T_m$ be an enumeration of the $r$-cliques in $S$. For a vector $\vec{x}\in[N+1]^m$, let $P(\vec{x})$ be the probability that for every $j\in[m]$, the clique $T_j$ is in $M_{x_j} \setminus M_{x_j-1}$ if $x_j \le N$, and $T_j \in \wt{M}$ if $x_j = N+1$. We can, essentially identically to the proof of \cref{eq:P_S} within the proof of \cref{thm:pm}, show that
\begin{equation}\label{eq:P_S-2}
P(\vec{x}) = \bigg( \prod_{i=1}^N \bigg( O_\eta\bigg( \frac{|V_{i-1}|}{n^r} \bigg) \bigg)^{|\{ j\colon x_j=i \}|} \bigg) \bigg( \frac{|V_N|}{n} \bigg)^{r|\{ j\colon x_j=N+1 \}|}.
\end{equation}
We truncate the details: we use the randomness guarantee \cref{itm:super-vortex-prob-guarantee} in \cref{lem:super-vortex} and the $O_\eta(1/|V_i|^{k-1})$-spreadness of the $M_i'=M_{i+1}\setminus M_i$. Then, similarly, we deduce
\begin{align*}
P_S = \sum_{ \vec{x} \in [N+1]^m } P(\vec{x}) & = \sum_{ \vec{x} \in [N+1]^m } \left( \prod_{i=1}^N \left( O_\eta\left( \frac{|V_{i-1}|}{n^k} \right) \right)^{|\{ j\colon x_j=i \}|} \right) \left( \frac{|V_N|}{n} \right)^{r|\{ j\colon x_j=N+1 \}|}\\
& = \left( \frac{O_\eta(1)}{n^{r-1}} \right)^m \sum_{ \vec{x} \in [N+1]^m } \left( \prod_{i=1}^{N} (\eta^i)^{|\{ j\colon x_j=i \}|} \right) \left(\frac{|V_N|^r}{n}\right)^{|\{ j\colon x_j=N+1 \}|}\\
& = \left( O_\eta\left( \frac{1}{n^{r-1}}\right) \right)^m \prod_{j=1}^m \left( \sum_{i=1}^N \eta^i + \frac{n^{r/(r+1)}}{n} \right) = \left( O_\eta\left( \frac{1}{n^{r-1}}\right) \right)^m.\qedhere
\end{align*}
\end{proof}

\section{Robust Hajnal--Szemer\'edi Theorem and counting \texorpdfstring{$K_r$}{Kr}-factors}\label{sec:robust-counting}

We are now in position to use our analysis in order to prove that a version of the robust Hajnal--Szemer\'edi theorem, \cref{thm:robust-corradi-hajnal}, holds. We first show that \cref{thm:robust-corradi-hajnal} follows from \cref{thm:robust-ch-clique}, which will be the main focus for the rest of the section.
\begin{proof}[Proof of \cref{thm:robust-corradi-hajnal} from \cref{thm:robust-ch-clique}]
We rely on theorems of Riordan (\cite[Theorem 1]{Rio18} for $r \geq 4$ and \cite[Theorem 16]{Rio18} for $r=3$) that couple the random binomial hypergraph with the clique complex of the random binomial graph. The result we need is that there exists a constant $a>0$ such that for every (fixed) $r \in \mb N$ and $p \leq \log^2(n)/n^{r-1}$, denoting $q = ap^{1/\binom{r}{2}}$, there exists a coupling of $\mb{G}^{(r)}(n,p)$ and $\mb{G}(n,q)$ such that w.h.p.\ $H' \sim \mb{G}^{(r)}(n,p)$ is contained in the $r$-clique complex of $G' \sim \mb{G}(n,q)$. (We note that Riordan proves a substantially more precise result when $r\ge 4$ obtaining the optimal constant $a$ for the complete hypergraph; an analogous result for $r=3$ was proven by \cite{Hec21}.)
 
Denote the $r$-clique complex of $G$ by $\mc{H}$. \cref{thm:robust-ch-clique} and \cref{thm:FKNP} imply that for some $C_1>0$, denoting $p = C_1 \log n / n^{r-1}$, w.h.p.\ $\mc{H}(p)$ (i.e., the random binomial subgraph of $\mc H$ with rate $p$) contains a perfect matching. Riordan's theorems imply that there exists a coupling of $G(a p^{1/\binom{r}{2}})$ and $\mc{H}(p)$ such that for $G' \sim G(a p^{1/\binom{r}{2}})$ and $H' \sim \mc{H}(p)$ w.h.p.\ $H'$ is contained in the clique complex of $G'$. In particular, w.h.p.\ both $H'$ both contains a perfect matching and is contained in the clique complex of $G'$, which together imply that $G'$ contains a $K_r$-factor.

For the counting result, note that we have an $O(1/n^{r-1})$-spread measure on the set of $K_r$-factors, each of which is composed of $n/r$ many $r$-cliques. Therefore, for some $C>0$, each factor occurs with probability at most $(C/n^{r-1})^{n/r}$ by the spread condition, so there are at least $(n^{r-1}/C^{r-1})^{n/r}=(n/C)^{(r-1)n/r}$ total factors.
\end{proof}

We define the key property of being somewhat near an extremal structure (complete balanced $r$-partite graph) in a specific sense.

\begin{definition}\label{def:extremal-conditions}
We say that graph $G$ is $(r,\alpha)$-sparse if there is $A\subseteq G$ with $|A|=\lfloor|V(G)|/r\rfloor$ such that $d_G(A)\le\alpha$. We say it is $\alpha$-disconnected if there is a partition $V(G)=A\cup B$ with $|A|=\lfloor|V(G)|/2\rfloor$ and $d_G(A,B)\le\alpha$.
\end{definition}

\subsection{Non-sparse setting}\label{sub:robust-non-sparse}
We next show the result when $G$ is not $(r,\alpha)$-sparse for appropriate $\alpha$.
\begin{lemma}\label{lem:robust-non-sparse}
Let $r|n$ and $\alpha<\alpha_{\ref{lem:robust-non-sparse}}(r)$ and $\theta=\theta_{\ref{lem:robust-non-sparse}}(r,\alpha,\alpha')$. Let $G$ be an $n$-vertex graph. If $\delta(G)\ge(r-1)n/r-\theta n$ and $G$ is not $(r,\alpha)$-sparse, and furthermore if $r=2$ then $G$ is not $\alpha$-disconnected, then there is a $C_{\ref{lem:robust-non-sparse}}(r,\alpha)/n^{r-1}$-spread distribution on the set of $K_r$-factors of $G$.
\end{lemma}

The proof is a slight simplification of the proof for triangles presented in \cite[Lemma~9.1]{ABCDJMRS22}; as the details are not as delicate in the non-extremal case we will be brief. The main task of the algorithm is noting that given a regularity partition, one can find a $K_r$-factor of the reduced graph covering almost all vertices, and the small remainder (and vertices within the regularity partition of exceptional degree) can be handled in a spread manner. 

We will first require a version of the Hajnal--Szemer\'edi Theorem itself \cite{HS70}.
\begin{theorem}\label{thm:Haj-Sz}
Let $n,k\ge 2$ be integers and let $0\le x< 1$. Suppose that $G$ is an $n$-vertex graph with $\delta(G)\ge\big(\frac{k-1}{k}-x\big)n$. Then $G$ contains a $K_k$-matching of size at least $(1-k(k-1)x)\lfloor \frac{n}{k}\rfloor$. 
\end{theorem}
 
The crucial input from \cite{ABCDJMRS22} is a robust fractional version of the Hajnal--Szemer\'edi Theorem. 
\begin{theorem}[{\cite[Theorem~7.4]{ABCDJMRS22}}]\label{thm:stab-fractional}
Fix $k\ge 2$ and $\eta > 0$. There exists $\gamma = \gamma(k,\eta)>0$ such that the following holds for all $m$. Let $G$ be a connected graph on $m$ vertices with $\delta(G)\ge ((k-1)/k-\gamma)m$ and $\alpha(G)<(1/k-\eta)m$. Let $\lambda\colon V(G)\to \mb{N}$ be a weight function such that $\lambda(u) = (1\pm \frac{\gamma}{2})(\frac{1}{m}\sum_{v\in V(G)}\lambda(v))$ and $\lambda(u)\ge m^{2k}$ for all $u\in V(G)$, and $k$ divides $\sum_{v\in V(G)}\lambda(v)$. Then there exists a weight function $\omega\colon K_k(G)\to \mb{N}\cup \{0\}$ such that $\sum_{\substack{K\in K_k(G)\\K\ni u}}\omega(K) = \lambda(u)$ for all $u\in V(G)$.
\end{theorem}

We now prove \cref{lem:robust-non-sparse}.
\begin{proof}[Proof of \cref{lem:robust-non-sparse}]
Fix a sequence of constants $0<\frac{1}{m_0}\ll \theta \ll \epsilon \ll d \ll \alpha\ll 1$ satisfying various constraints throughout the proof. Given a graph $G$, apply \cref{lem:regularity} and consider the $(\epsilon,d)$-reduced graph $R$ which is returned and note that $R$ has $m\in [m_0, M_0]$ vertices. Let the underlying $\epsilon$-regular partition of $V(G)$ be $V_0\cup V_1\cup \cdots\cup V_m$ and note by \cref{lem:regularity} that $\delta(R)\ge ((r-1)/r - 4d) m$. 

We first claim that the independence number of $R$ (denoted $\alpha(R)$) is suitably large. This is the analogue of \cite[Claim~9.4]{ABCDJMRS22}.
\begin{claim}\label{clm:small-independence}
We have that $\alpha(R)<(\frac{1}{r}-\alpha^{2})m$.
\end{claim}
\begin{proof}
Suppose that $R$ has an independent set $S$ of size $(1/r-\alpha^2)m$. Let $S' = \bigcup_{j\in S}V_j\subseteq V(G)$. By the definition of the $(\epsilon,d)$-reduced graph, we have that $e_G(S')\le (2\epsilon + d) n^2$ and $|S'|\ge (1-\epsilon)(1/r-\alpha^2)n\ge (1/r-2\alpha^2)n$. Adding an arbitrary set of $n/r - |S'|$ many additional vertices to $S'$, gives a set of size exactly $n/r$ with at most $4\alpha^2n^2$ edges in $G$. This contradicts the fact that $G$ is not $(r,\alpha)$-sparse.
\end{proof}

We next require that $R$ is a connected graph. This is necessary to verify the connectedness assumption which appears within \cref{thm:stab-fractional}; this is the unique place where the assumption that $G$ is not $\alpha$-disconnected is required. 
\begin{claim}\label{clm:connected}
The graph $R$ is connected.
\end{claim}
\begin{proof}
For $r\ge 3$, the claim is immediate as $\delta(R)\ge 3m/5$. For $r = 2$, note that we have that $\delta(R) \ge (1/2 - 4d)m$. Therefore, if $R$ is not connected, then there are at most $2$ connected components each of size at least $(1/2-4d)m$. Let the connected components of $R$ be $S_1$ and $S_2$ and define $S_i' = \bigcup_{j\in S_i}V_j$ for $i\in [2]$. By the definition of an $(\epsilon,d)$-reduced graph, we have that $e_G(S_1',S_2')\le(2\epsilon + d)n^2$ and that $|S_1'|, |S_2'|\ge (1-\epsilon)(1/2-4d)n\ge (1/2 - 5d)n$. This immediately implies, using $|V_0|\le \epsilon n$, that $e_G(V_0 \cup S_1',S_2')\le 5dn^2$ and that $||(V_0 \cup S_1')|-|S_2'||\le 11dn$. Rebalancing $(V_0 \cup S_1')$ and $S_2'$ to give an equipartition of $V(G)$, we obtain a contradiction to the fact that $G$ is not $\alpha$-disconnected.
\end{proof}

We now consider the reduced graph $R$ and the induced partition on the vertex set $V_0\cup V_1\cup\cdots\cup V_m$. By applying \cref{clm:small-independence} and \cref{thm:Haj-Sz}, there exists a partial $K_r$-factor $\mc{T}_r$ of the reduced graph $R$ which cover all but $O_{r}(dm)$ vertices in $R$. Let $\mc{T}^{\ast} = V(\mc{T}_r)$ and for each clique $\{i_1,\ldots,i_r\}$ in $\mc{T}_r$, we can pass to a subset $V_{i_1}^{\ast},\ldots, V_{i_r}^{\ast}$ with $|V_{i_r}^{\ast}| = (1-r\epsilon)|V_{1}|$ and $(V_{i_j}^{\ast},V_{i_\ell}^{\ast})$ being $(2\epsilon, (d-\epsilon)^{+},d-k\epsilon)$-super-regular (see \cite[Lemma~2.9]{ABCDJMRS22}).

We define $X = V_0 \cup\bigcup_{j\notin \mc{T}^\ast} V_j \cup\bigcup_{i\in [m]}(V_i\setminus V_{i}^{\ast})$. Note that $|X|\lesssim_{r} dn$. We will construct a (random) clique factor in $G$ by first matching the vertices into cliques with a random greedy algorithm and then (after rebalancing the partition) applying \cref{thm:super-regular} to find a clique factor within each clique of $\mc{T}_r$.

We now proceed with the following algorithm.
\begin{itemize}
    \item Order the vertices in $X$ as $\{v_1,\ldots,v_{|X|}\}$ arbitrarily. Define $G_0= V(G)$; $G_i$ will correspond to the vertex set after the vertex $v_i$ has been matched. 
    \item For each vertex $v_{\ell}$ in $X$, choose a uniformly random clique extending $v_{\ell}$ within $G_{\ell-1}$ which does not contain an additional vertex of $X$. Update $G_\ell$ to be the vertices in $G_{\ell-1}$ minus the set of vertices in the chosen clique. 
\end{itemize}

We now prove a number of basic properties of the algorithm
\begin{claim}\label{clm:cleanup-phase-1}
The algorithm satisfies the following properties:
\begin{itemize}
    \item The algorithm always runs to completion;
    \item The random set of cliques created by the algorithm is $O_r(1/n^{r-1})$-spread;
    \item For any subset of vertices $S\subseteq V(G)\setminus X$ we have 
    \[\mb{P}[|S\cap (G_0\setminus G_{|X|})|\ge \sqrt{d}|S|]\le \exp(-\Omega_r(\sqrt{d}|S|)).\]
\end{itemize}
\end{claim}
\begin{proof}
For the first part, we consider the degree of $v_{\ell}$ in the remaining graph. Notice that the degree of $v_{\ell}$ is always at least $((r-1)/r)n - \theta n - r|X|\ge ((r-1)/r - C_rd)n$ for an appropriate constant $C_r$. Note that any subset of size $T$ of size $(r-1)n/r$ in $G$ has at least
\begin{align*}
e(G[T]) &= \frac{1}{2}\sum_{v\in T}\deg_G(v) - \frac{1}{2}|G[T,T^{c}]| \ge |T|/2 \cdot ((r-1)n/r - \theta n) - \frac{1}{2}|T|(n/r)\\
&\ge |T|/2\cdot ((r-2)n/r - d^2n) \ge |T|^2/2 \cdot ((r-2)/(r-1)-d)
\end{align*}
edges. This implies that the density of the edge set of the neighbors of $X$ is at least $(r-2)/(r-1) - C_r'd$ and therefore by supersaturation for Turan's theorem there are at least $\Omega_r(n^{r-1})$ possible cliques at each stage. (Notice that the Turan threshold for finding a $K_{r-1}$ is $(r-3)/(r-2)$ which is strictly below the density specified.) This implies the first two statements in the claim. 

For the third claim, notice that there are at most $|S|n^{r-1}$ cliques of size $r$ containing an element of $S$ and a fixed vertex in $X$. Thus the result follows noting that $|X|\le \sqrt{d}n/r$, that we can remove at most $r$ vertices in $S$ at each stage, and applying the binomial domination lemma (\cref{lem:binomial-comparison}).
\end{proof}

By applying \cref{clm:cleanup-phase-1} we have with probability at least $1/2$ that 
\begin{itemize}
    \item $|V_{i}^{\ast}\cap G_{|X|}|\ge (1-\sqrt{d})|V_{i}^{\ast}|$;
    \item for each edge in the clique factor $\mc{T}_r$ we have that the corresponding pair of parts formed by $V_i^{\ast}\cap G_{|X|}$ and $V_j^{\ast}\cap G_{|X|}$ are still $(4\epsilon, (d/2)^{+},d/4)$-super-regular by considering the number of vertices deleted in each part, and the number of neighbors of a given vertex which are deleted, controlled via the third bullet of \cref{clm:cleanup-phase-1}; and
    \item For a pair $(V_i^{\ast},V_j^{\ast})$ where each $(i,j)$ appears in $R$, we have that $(V_i^{\ast}\cap G_{|X|},V_j^{\ast}\cap G_{|X|})$ is $(4\epsilon, (d/2)^{+})$-regular. This is immediate by considering the number of vertices deleted from $V_i$ to obtain $V_i^{\ast}\cap G_{|X|}$.
\end{itemize}

For the sake of clarity define $V_i' = V_i^{\ast}\cap G_{|X|}$ for $i\in V(\mc{T}_r)$. At this point we face the difficulty that the parts corresponding to each clique of the $K_r$-factor in $\mc{T}_r$, while relatively close in size, are not balanced appropriately. We will use \cref{thm:stab-fractional} in order to remove a certain number of triangles and ``rebalance'' the part sizes in $V_i'$. For conciseness let $R' \coloneqq R[V(\mc{T}_r)]$ denote the subgraph of $R$ induced by $V(\mc{T}_r)$.

In order to apply \cref{thm:stab-fractional}, we define $\lambda(i) = |V_i'| - \lceil n/m (1-d^{1/3})\rceil$. Notice that $\sum_{i\in V(\mc{T}_r)}|V_i'|$ is the number of remaining vertices and that this number is divisible by $r$ (since these are the vertices remaining after deleting a set of $r$-cliques from $G$). Additionally the number of remaining parts is $r|\mc{T}_r|$. Hence $\sum_{i\in V(\mc{T}_r)} \lambda(i)$ is divisible by $r$. Furthermore $\lambda(i)\in (d^{1/3} \pm d^{1/2})n/m$, $\delta(R')\ge ((r-1)/r-C_rd)m$, $\alpha(R')<(1/r-\alpha^2)m$, and $r|\mc{T}_r| = m (1\pm O_r(d))$. As $d\ll \alpha$, \cref{thm:stab-fractional} implies that there exists a weight function $\omega \colon K_r(R')\to \mb{N}$ such that for all $i\in V(\mc{T}_r)$ we have that
\[\sum_{\substack{K\in K_r(R')\\K\ni i}} \omega(K) = \lambda(i).\]
(Note that strictly speaking one needs to check that the graph on $R'$ is connected, but the proof given in \cref{clm:connected} is robust to perturbations of the vertex set of $R$ of size $O(dm)$.)

Notice that if we can remove $\omega(K)$ vertex-disjoint cliques from the corresponding $r$-partite set of regular pairs then we will be left with $|V_i'|-\lambda(i) = \lceil n/m (1-d^{1/3})\rceil$ vertices in each part. This gives the desired result, provided the parts left are sufficiently well behaved. The difficulty is guaranteeing the necessary super-regularity at the end of the algorithm; to do so we split $V_i' = V_i'' \cup V_i'''$ where each vertex of $V_i'$ is placed with probability $1/2$ independently at random in $V_i''$ or $V_i'''$. Consider the clique $K$ which contains $i$ in $\mc{T}_r$; applying super-regularity, Chernoff's inequality, and a union bound we have that for all remaining $j\in K$ and for each vertex $v\in V_j'$ there holds $\deg(v,V_i'')\ge (d/8)|V_i''|$. Furthermore, we have that $|V_i''| = (1/2 \pm 1/4)|V_i'|$ with high probability.

We now proceed with the following algorithm.
\begin{itemize}
    \item Order the cliques in $K_r(R')$ in an arbitrary manner.
    \item For each successive clique $K \in K_r(R')$, or $\omega(K)$ steps, iteratively remove a random clique with one vertex in each part $V_i'''$ such that $i \in K$ and which does not intersect previously chosen cliques. 
\end{itemize}

We now note that the above algorithm indeed runs to completion and is appropriately spread. To see that the algorithm runs to completion note that we remove $\lambda(i) \leq 2d^{1/3}n/m$ vertices from each part $V_i'''$ and hence we have at least (say) $|V_i'''|/2$ choices for the vertex in each part. Hence, the counting-lemma guarantees that there are $\Omega(n^{r})$ choices in each stage. Furthermore the algorithm is trivially sufficiently spread as there are only $O(n)$ rounds and the choice of each clique is uniformly random among a set of size $\Omega(n^{r})$ at each step (conditional on the previous choices). Finally, pairs for each $K \in \mc{T}_r$ are still suitably super-regular as the minimum degree condition is preserved by looking at edges in $V_i''$ and regularity is preserved as we are left with a constant fraction of each vertex set $V_i'$. We then take each pair of parts appearing in a clique in $\mc{T}_r$, apply \cref{lem:plus-to-super-regular}, and then apply \cref{thm:super-regular} in order to give an $O_{d,\eps}(1/n^{r-1})$-spread factor covering the remaining vertices. This completes the proof.
\end{proof}

\subsection{Reduction to non-sparse setting}\label{sub:robust-reduction}
We now prove \cref{thm:robust-ch-clique} by reducing to applications of \cref{lem:robust-non-sparse} and \cref{thm:super-regular}.
\begin{proof}[Proof of \cref{thm:robust-ch-clique}]
We are given $G$ on $n \in r \mb{N}$ vertices $V=V(G)$ with $\delta(G)\ge(r-1)n/r$ and wish to construct a sufficiently well-spread distribution on the perfect matchings in $\mc{H}$, the set of $r$-cliques in $G$. We may assume without loss of generality that in $G$, any two vertices of degree strictly greater than $(r-1)n/r$ are not connected via an edge, else we can reduce the situation by deleting the edge. The argument is similar to \cite[Section~6]{KSS01}. Given $r$, consider parameters
\[
1/n\ll\alpha_1\ll\alpha_2\ll\cdots\ll\alpha_r\ll1/r
\]
with appropriate space between each pair.

We iterate over $i\in\{1,\ldots,r\}$ and at each step, if possible, find $A_i\subseteq V\setminus(A_1\cup\cdots\cup A_{i-1})$ such that $|A_i|=n/r$ and $d_G(A_i)\le\alpha_i$. Let $m\in\{0,\ldots,r\}$ be the number of steps that successfully go through. Let $A_{\ge 1}=A_1\cup\cdots\cup A_m$ and $B=V(G)\setminus A_{\ge 1}$. We write $A_0=B$ for convenience, and let the resulting partition be denoted $\mc{A}$. Note that $m=0$ corresponds directly to \cref{lem:robust-non-sparse}, except in the case $r=2$ where there is a slight difference; in general we will apply \cref{lem:robust-non-sparse} with $r$ replaced by $r-m$ to create a $K_{r-m}$-factor on a set roughly similar to $B$, then find a $K_{m+1}$-factor of a corresponding reduced near-partite graph using \cref{thm:super-regular}.

Note that $m\neq r-1$: if $m=r-1$, then we can check that $A_1,\ldots,A_{r-1}$ and $A_r:=A_{\ge 1}=B$ have the property that they have at least a $1-O(\alpha_{r-1})$ fraction of edges between each pair. On the other hand, the condition that vertices of degree greater than $(r-1)n/r$ are not connected implies that every edge in $G[A_r]$ has at least one endpoint of degree $(r-1)n/r$, and thus corresponds to an edge between that vertex and $A_1\cup\cdots\cup A_{r-1}$ that is ``missing''. The claim easily follows, using $\alpha_{r-1}\ll\alpha_r$.

Let $\alpha=\alpha_m$ and $\beta=\alpha_{m+1}$, note $\alpha\ll\beta$, and note that $G[A_1],\ldots,G[A_m]$ have densities bounded by $\alpha$ while $G[B]$ is such that every set of size $n/r$ has density at least $\beta$. Consider $\eta$ satisfying $\beta\ll\eta\ll1/r$ with appropriate space between the pairs. We define a notion of vertices that do not look like they respect the partition, in the sense that their degrees are not what they should be if $G[A_1],\ldots,G[A_m]$ were empty.

\begin{definition}[Bad and exceptional vertices]\label{def:bad}
Given a partition $\mc{P}=P_1\cup\cdots\cup P_m\cup P_0$ of $V$ with $P_{\ge 1}=P_1\cup\cdots\cup P_m$ and $|P_i|=n/r$ for $i\in[m]$, we define \emph{bad} and \emph{exceptional} vertices as follows. For $i\in[m]$, let $v\in V$ be called \emph{$(i,\eta^\ast)$-bad wrt $\mc{P}$} if $\deg_G(v,P_i)\ge\eta^\ast|P_i|$. We say that $v$ is \emph{$(0,\eta^\ast)$-bad wrt $\mc{P}$} if $\deg_G(v,P_{\ge 1})\le(1-\eta^\ast / m)|P_{\ge 1}|$. For $j\in[m]$ we say that a vertex $v\in V$ is \emph{$j$-exceptional wrt $\mc{P}$} if $\deg_G(v,P_j)\le\eta^\ast|P_j|/2$. For $i\in[m]$ we say that $v\in V$ is \emph{$(0,\eta^\ast)$-exceptional wrt $\mc{P}$} if $\deg_G(v,P_0)\le(r-m-1+\eta^\ast/2)n/r$. We will often drop $\eta^\ast$ and specification of $\mc{P}$ from the notation where it is clear.
\end{definition}

For the purpose of the argument we have defined these notions over all of $V(G)$, but we will be most interested in $i$-bad vertices contained in $P_i$ as well as $j$-exceptional vertices that are not contained in $P_j$.

Notice that if $v$ is $(j,1/8)$-exceptional with respect to $\mc{P}$ then almost all of its ``degree deficit'' (which is at most $n/r-1$) is used up by the non-edges between $v$ and $P_j$, and hence $v$ is nearly complete to $P_i$ for all $i\neq j$ (and is hence $(i,1/8)$-bad for all $i\notin\{0,j\}$). Furthermore, if $i=0$ and $j\neq i$ then by inspection we see $v$ is still $(i,1/8)$-bad. Thus, if $v$ is $(j,1/8)$-exceptional then it is $(i,1/8)$-bad for all $i\neq j$.

Finally, given a partition $\mc{P}$ of some vertex set $P_1\cup\cdots\cup P_m\cup P_0$, we say it is \emph{balanced} if $|P_1|=\cdots=|P_m|$ and $|P_0|=(r-m)|P_1|$. We say a clique $K_r$ is balanced with respect to $\mc{P}$ if it has $1$ vertex in each $P_i$ for $i\in[m]$, and $r-m$ vertices in $P_0$.

\subsubsection{Connected case}\label{sub:connected-subcase}
We first assume that either $m\neq r-2$ or, if $m=r-2$, that $G[B]$ is not $\alpha_{m+1}$-disconnected. For each $i \in [m]$, the fact that $d_G(A_i)\le\alpha$ implies that for every $\eta \geq \alpha^{1/8}$ the set $A_i$ contains at most $\alpha^{2/3}n$ many $(i,\eta^2)$-bad vertices wrt $\mc{A}$. A similar argument shows that $B$ contains at most $\alpha^{2/3}n$ many $(0,\eta^2)$-bad vertices wrt $\mc{A}$. Indeed, by the minimum degree condition and $d_G(A_i)\le\alpha$ we see that the number of edges $e_G(A_i,B)$ is at least $(r-m)n^2/r^2-2\alpha n^2$ hence
\[e_G(A,B)\ge m(r-m)n^2/r^2-2m\alpha n^2.\]
But $|B|=(r-m)n/r$ and $\deg_G(v,A)\le mn/r$ for every $v \in B$, so Markov applied to the quantity $mn/r-\deg_G(v,A)$ over $v\in B$ yields the desired result.

\textbf{Step 1: Cleaning the partition.}
We modify $\mc{A}$ into a partition $\mc{A}'$ slightly in the following manner:
\begin{itemize}
    \item Initialize with $\mc{A}'=\mc{A}$.
    \item If there are distinct $i,j\in\{0,\ldots,m\}$ and $v\in A_i'$ which is $(i,\eta^2)$-bad wrt $\mc{A}$ and $w\in A_j'$ which is $(i,\eta^{1/2})$-exceptional wrt $\mc{A}$, we swap the positions of $v,w$ in $A_i',A_j'$.
    \item Continue this operation until there are no possible choices, then terminate.
\end{itemize}
We emphasize that in this algorithm badness and exceptionality are always wrt the original partition $\mc{A}$, rather than the evolving partition $\mc{A}'$.

We claim that the process terminates in at most $\alpha^{1/2}n$ steps. Indeed, note that at every step the following quantity, which by the analysis above is initially less than $\alpha^{3/5}n$, strictly decreases: the number of $(i,v)\in\{0,\ldots,m\}\times V(G)$ such that $v\in A_i'$ and $v$ is $(i,\eta^2)$-bad wrt $\mc{A}$. This is because at the start of a step both $v$ and $w$ contribute to this condition but at the end of it only $v$ can contribute (and everything else is left unchanged since we are measuring wrt the original partition $\mc{A}$). Here we are using the fact that $(i,1/8)$-exceptional vertices are $(j,1/8)$-bad for all $j\neq i$.

For each $i \in \{0,\ldots,m\}$, let $A_i'$ denote the part of $\mc{A}'$ that corresponds to $A_i$. Since $\mc{A}'$ differs from $\mc{A}$ by few vertices, we actually find that for each $i\in\{0,\ldots,m\}$, either there is no $v\in A_i'$ which is $(i,\eta)$-bad wrt $\mc{A}'$ or there are no $(i,\eta)$-exceptional vertices wrt $\mc{A}'$ not in $A_i'$ (or both). Furthermore, $d_G(A_i')\le\alpha^{1/2}$ for all $i\in[m]$, say, and the part sizes of $\mc{A}'$ are balanced in the same way as before. The goal at this point is to cover vertices with ``balanced'' cliques, i.e., ones where there are $r-m$ vertices in $A_0'$ and $1$ in each $A_i'$ for $i\neq 0$. However, the exceptional vertices will require a different treatment.

\textbf{Step 2: Covering exceptional vertices.}
Our first goal is to cover all exceptional vertices by $r$-cliques in a spread manner, leaving a balanced partition behind. The idea is that an $(i,\eta)$-exceptional vertex, even though it is in some $A_j'$, is better utilized if we ``swap'' it to $A_i'$. To counterbalance this, something must be ``swapped'' back. Hence, we will cover $i$-exceptional vertices in $A_j'$ by $r$-cliques which contain one ``extra'' vertex in $A_j'$ and one fewer vertex in $A_i$ (compared to a balanced clique), and this will be accompanied by a clique with one extra vertex in $A_i'$ and one fewer vertex in $A_j$.

Suppose there are $x_i\neq 0$ many $(i,\eta)$-exceptional vertices wrt $\mc{A}'$ that are not in $A_i'$. Since there exist such vertices, the earlier cleaning step shows that there are no $v\in A_i'$ that are $(i,\eta)$-bad wrt $\mc{A}'$. Now we double-count $e_G(A_i',V\setminus A_i')$. If $i\neq 0$ then every exceptional vertex contributes at most $\eta n/2$ edges, and the rest contribute at most $n/r$ edges, so that
\[e_G(A_i',V\setminus A_i')\le (r-1)n^2/r^2-x_i(n/r-\eta n/2).\]
On the other hand, each of the $n/r$ vertices in $A_i'$ has degree at least $(r-1)n/r$, so the above demonstrates that $e_G(A_i')\ge x_in/(3r)$. A similar argument works for $i=0$, except that the bound is replaced by $e_G(A_0',V\setminus A_0')\le m(r-m)n^2/r^2-x_i(n/r-\eta n/2)$ and thus in fact
\[e_G(A_0')\ge \frac{1}{2}(r-m)(r-m-1)n^2/r^2+x_in/(3r).\]

In the first case $e_G(A_i') \geq x_in/(3r)$, and since there are no $(i,\eta)$-bad vertices in $A_i'$, the induced graph $G[A_i']$ has maximum degree at most $\eta n$. This graph also has $n/r$ vertices since $i\neq 0$. In the case where $i=0$, there are $(r-m)n/r$ vertices and at least $n^2/(2r^2)$ edges, since $r\neq m-1$ (and there is no $i=0$ case when $r=m$). We now sample a size-$x_i$ matching $M \subseteq G[A_i']$ as follows: Let $M' \subseteq G[A_i']$ be a binomial random set of edges where each edge appears with probability $6r/n$. Next, let $M'' \subseteq M'$ be the set of edges that are vertex disjoint from all other edges in $M'$. Let $M$ be the distribution obtained by taking a uniformly random size-$x_i$ subset of $M''$, conditioned on $|M''| \geq x_i$.
%First, sample each edge with probability $6r/n$ and delete edges whose endpoints are included in more than $1$ sampled edge, then condition on having at least $x_i$ remaining edges. Then uniformly at random choose $x_i$ such edges.

Clearly, $M''$ is $O(1/n)$-spread. We claim that $\mb{P}[|M''| \geq x_i] \geq 1/100$. Indeed, when $i \neq 0$, this follows from the maximum degree condition. When $i=0$ we have
\[
\mb{E}[|M''|] \geq e_G(A_0')(6r/n)(1-6r/n)^{2n} \geq n^2/(2r^2)(6r/n)e^{-4r} \geq \alpha^{3/5}n \geq 2x_i.
\]
A second moment calculation implies that $\mb{P}[|M''| \geq x_i] \geq 1/100$ in this case as well. In both cases it follows that $M$ is $O(1/n)$-spread.
%The maximum degree condition in $A_i'$ demonstrates that the desired event we condition on occurs with probability at least $1/100$ when $i\neq 0$. When $i=0$, the fact that $x_i$ is small, say less than $r\alpha^{3/5}n$, implies the claim. Either way, we additionally see that the choice of edges is $O(1/n)$-spread. Having done this, for all $i\in\{0,\ldots,m\}$

We now cover all $(i,\eta)$-exceptional vertices not in $A_i'$ and all edges in $M$ by $r$-cliques in a sufficiently spread manner. Notice that the total amount of vertices and edges to cover is small --- say at most $\alpha^{2/5}n$. Furthermore, every $(i,\eta)$-exceptional vertex not in $A_i'$ has used up most of its ``degree deficit'' on $A_i'$, which means that it is nearly complete to $V\setminus A_i'$. Additionally, edge in $M$ has both endpoints not $(i,\eta)$-bad (since there are no $i$-bad vertices when there are $i$-exceptional vertices outside of $A_i'$) hence they similarly are nearly complete to $V\setminus A_i'$. We therefore easily find that the $(i,\eta)$-exceptional vertices in $A_j'$ (with $j\neq i$) are contained in $\Omega(n^{r-1})$ many $r$-cliques with one ``extra'' vertex in $A_j'$ and one fewer in $A_i'$; similarly, the constructed edges within $A_i'$ are contained in $\Omega(n^{r-2})$ many $r$-cliques with one ``extra'' vertex within $A_i'$ and one fewer in some prescribed $A_j'$. (Here we are using the robust counting version of Tur\'an's theorem.)

Choosing such $r$-cliques iteratively uniformly at random such that they do not overlap the previously chosen cliques, we can easily argue that the resulting collection of $r$-cliques is $O(1/n^{r-1})$-spread. Furthermore, since we chose $x_i$ edges within each $A_i'$, the removal of all these $r$-cliques leaves a resulting (random) partition $\mc{A}''$, $V'=A_1''\cup\cdots\cup A_m''\cup A_0''$ which is still balanced. We additionally remark that so far we have removed at most $\alpha^{1/3}n$ vertices.

\textbf{Step 3: Covering $0$-bad vertices.}
Our second goal is to cover all $(0,\eta)$-bad vertices (wrt $\mc{A}'$) that remain in $A_0''$ (in a spread manner and leaving a balanced partition behind to which we can begin to apply \cref{thm:super-regular}). Again, there are very few of them --- at most $\alpha^{1/3}n$. At this point, since we have removed all $(i,\eta)$-exceptional vertices (wrt $\mc{A}'$) outside $A_i'$ for each $i$, every such $v\in A_0''$ has degree at least $\eta n/(2r)$ to each $A_1',\ldots,A_m'$ and hence degree at least $\eta n/(3r)$ to each $A_1'',\ldots,A_m''$.

We uniformly at random choose vertices in each such neighborhood (that have not yet been chosen), conditional on forming an $(m+1)$-clique with $v$ and conditional on each such vertex in say $A_j''$ being $(j,\eta)$-good wrt $\mc{A}'$ (there are at most say $\alpha^{1/3}n$ such bad vertices, but at least $\eta n/(3r)$ choices, and $\alpha\ll\eta$, so this conditioning does not distort the randomness too much). Then the common neighborhood of these $m$ extra vertices within $A_0''$ is nearly all of $A_0''$, and we thus easily argue that there are $\Omega_\eta(n^{r-1})$ total choices of balanced cliques $K_r$ containing $v$. (We are again using the robust counting version of Tur\'an's theorem.) We iteratively choose such a clique uniformly at random for each $(0,\eta)$-bad vertex in $A_0''$ (so that each clique is vertex-disjoint from its predecessors). Similar to before, this random collection of cliques is $O_\eta(1/n^{r-1})$-spread (conditional on the prior randomness).

We have removed a total of at most $O(\alpha^{1/3}n)$ vertices, and we have a remaining partition $\mc{A}'''$, $V''=A_1'''\cup\cdots\cup A_m'''\cup A_0'''$ which is still balanced. Furthermore, we now have no $(0,\eta)$-bad vertices wrt $\mc{A}'$ in $A_0'''$ and also no $(i,\eta)$-exceptional vertices wrt $\mc{A}'$ in $V''\setminus A_i'''$.

\textbf{Step 4: Covering the remainder with \cref{thm:super-regular}.}
Note that $A_0'''$ and $A_0$ differ by only $O_\alpha(n)$ vertices. Thus, since $\alpha\ll\beta$, $G[A_0''']$ satisfies the hypotheses of \cref{lem:robust-non-sparse} with $r$ replaced by $r-m$ and $\alpha$ replaced by $\beta/2$, using also that $\mc{A}'''$ is balanced and hence $|A_0'''|/(r-m)\in\mb{N}$.

Now we perform the following process. By \cref{lem:robust-non-sparse} applied to $G[A_0''']$ in the manner above, there is a $O_\beta(1/n^{r-m-1})$-spread distribution on the set of $K_{r-m}$-factors of $G[A_0''']$. Sample from this distribution. Now, since every vertex here is not $(0,\eta)$-bad wrt $\mc{A}'$, every vertex in these $(r-m)$-cliques has degree at least $(1-\eta/m)|A_{\ge 1}'|$ to $A_{\ge 1}'$. Hence each of these $(r-m)$-cliques have common degree at least $(1-\eta^{1/2})|A_{\ge 1}'''|$ to $A_{\ge 1}'''$. Create an auxiliary graph $G'$ which is $G[A_{\ge 1}''']$ along with a collection of vertices $C$ corresponding to these $(r-m)$-cliques, each clique connected to its common neighbors in $A_{\ge 1}'''$ within $G$. We see that each pair of parts among the $A_i'''$ have the same size and have density at least $1-\alpha^{1/4}$. Furthermore, $A_i'''$ and $C$ have the same size and the density of edges between them is also at least $1-\alpha^{1/4}$ (indeed, in expectation, a random vertex $v$ in $A_i'''$ is missing at most $\alpha^{1/3}$-fraction of the crossing edges to $A_0'''$ in $G$, and therefore the fraction of $(r-m)$-cliques that do not fully connect to $v$ is on average at most $r\alpha^{1/3}$). Hence, these bipartite graphs are so dense that they are automatically regular pairs with error parameter depending on $\alpha$. Additionally, the covering that we have done so far has ensured that between every pair, the minimum degree is at least say $\eta n/(4r)$, and $\alpha\ll\eta$.

Thus we can apply \cref{lem:super-boost,lem:plus-to-super-regular} to obtain super-regular pairs between each part and then use \cref{thm:super-regular} to obtain a $O_\alpha(1/n^m)$-spread distribution on $K_{m+1}$-factors in $G'$. This corresponds to a $K_r$-factor of $G[\bigcup\mc{A}''']$, so we have constructed a full $K_r$-factor. Furthermore, we can easily argue due to the spread of the $K_{r-m}$ factor and then spread of the $K_{m+1}$-factor in $G'$ that the factor produced at this stage, conditional on the randomness in the previous steps, is $O_\alpha(1/n^{r-1})$-spread. This completes the connected case.

\subsubsection{Disconnected case}\label{sub:disconnected-subcase}
We now consider the only remaining case, that $m=r-2$ and $G[B]$ is $\beta$-disconnected (recall that $\beta=\alpha_{m+1}$). The argumentation is very similar to the connected case in \cref{sub:connected-subcase}, with two differences. First, since $G[B]$ is $\beta$-disconnected, we can essentially break $G[B]$ into two nearly-complete parts which are mostly disconnected. Therefore almost all the balanced cliques we use will have the additional property that they have $2$ vertices in one of these parts. However, in the case that these nearly-complete parts both have odd size, we will need to use $1$ clique with the property that it has $1$ vertex in each of the parts. Beyond this, the pruning of exceptional and $0$-bad vertices is essentially the same. We will therefore not repeat said argumentation in detail, merely stating the essential points.

\textbf{Step 1: Cleaning the partition.}
We modify $\mc{A}$ into $\mc{A}'$ similar to the connected case (\cref{sub:connected-subcase}). $\mc{A}'$ is still balanced, differs from $\mc{A}$ by few vertices (say $\alpha^{1/3}n$), and for each $i\in\{0,\ldots,m\}$ either there is no $v\in A_i'$ which is $(i,\eta)$-bad wrt $\mc{A}'$ or there are no $(i,\eta)$-exceptional vertices wrt $\mc{A}'$ not in $A_i'$ (or both). Since $\mc{A}$ is close to $\mc{A}'$ and $\alpha\ll\beta$, we know that $B'=A_{\ge 1}'$ has the property that $G[B']$ is $2\beta$-disconnected, say.

\textbf{Step 2: Partitioning $B'=A_0'$.}
Consider a partition $B'=B_1'\cup B_2'$ constructed as follows: start with the equipartition $B'=B_1^{(0)\prime}\cup B_2^{(0)\prime}$ guaranteed by the $2\beta$-disconnectedness of $G[B']$ (\cref{def:extremal-conditions}). Then, for each time $t\ge 1$, we define $B_1^{(t)\prime}\cup B_2^{(t)\prime}$ by swapping a vertex $v\in B_i^{(t-1)\prime}$ for some $i\in\{1,2\}$ such that $\deg_G(v,B_i^{(t-1)\prime})\le n/(4r)$. If this is no longer possible we terminate and set $(B_1',B_2')=(B_1^{(t-1)\prime},B_2^{(t-1)\prime})$. Each non-terminating step clearly decreases the cut size by at least $n/(2r)$ each time, and the initial cut size is $e(G[B_1^{(0)\prime},B_2^{(0)\prime}])\le\beta(2n/r)^2/2$. Therefore there are at most $\beta n$ total steps, meaning that the process terminates and furthermore the part sizes are $\beta n$-close to the original. In particular, we find $e(G[B_1',B_2'])\le\beta n^2$, say. Additionally, we see that $\delta(G[B_i'])\ge n/(5r)$ for $i\in\{1,2\}$. We deduce that $e(G[B_i'])\ge(1-O(\beta))(n/r)^2/2$. So the induced graph in each part $B_i'$ is almost-complete, and we have a reasonable minimum degree condition within each part.

\textbf{Step 3: Fixing parities.}
Notice that $|B_1'|+|B_2'|=2(n/r)$ is even. If $B_1',B_2'$ have even size then there is no need to do anything additional and we can move on to the next step. However, if they are both odd size then we wish to make them even in some way. To do this, we attempt to find (in a spread way) a single $K_r$ which is balanced and has $1$ vertex in each $B_i'$. Let us assume $|B_1'|\ge|B_2'|$.

Since $|B_1'|\ge|B_2'|$ we have $|B_2'|\le n/r$ so $\deg_G(v,B_1')\ge 1$ for each $v\in B_2'$ by the minimum degree condition on $G$. Choose a uniformly random edge in $G[B_1',B_2']$ such that its vertex in $B_2'$ is not $(i,\eta)$-exceptional for any $i\in[r-2]$ and is not $(0,\eta/(4r))$-bad. Notice that by Markov, most of the vertices of $B_2'$ satisfy this property and combining with $\deg_G(v,B_1')\ge 1$ for all $v\in B_2'$ shows that the are at least $n/(4r)$ choices for this edge, which is thus $O(1/n)$-spread. Call the edge $e_0=(u_0,v_0)$ where $u_0\in B_1'$ and $v_0\in B_2'$. In the case $r > 2$, if $u_0$ is not $(i,\eta)$-exceptional for any $i\neq 0$, then the common neighbors of $u_0,v_0$ among $A_i$ number at least $\eta n/4$. We easily find many balanced cliques containing this edge in this case, using that almost all vertices (up to error depending only on $\alpha$) are non-exceptional and good in the relevant senses. We then choose a uniformly random such clique which is $O(1/n^{r-1})$-spread. Now after removing this one $K_r$, the remaining partition $\mc{A}^\ast$ is balanced and the new $B_1^\ast,B_2^\ast$ are even in size. If $r=2$ then just choosing the spread edge is enough as we do not need to extend it further, and no condition on $0$-badness is needed or even meaningful (and we now move to the next step).

When $r > 2$ and the endpoint $u_0$ is $(i,\eta)$-exceptional for some $i\neq 0$, then $\deg_G(v,A_{\ge 1}')\le(r-3+\eta/2)n/r$ so $\deg_G(v,B')\ge(2-\eta)n/r$. We can then move $v$ to $B_{3-i}'$ to create $B'=B_1^\ast\cup B_2^\ast$. Now both sizes are even in size again, and the minimum degree in both parts is still at least $n/(4r)$. We continue to the next step.

\textbf{Step 4: Exceptional and $0$-bad vertices.}
At this stage, we cover the exceptional and $0$-bad vertices with basically the same arguments as in \cref{sub:connected-subcase}. When covering $(i,\eta)$-exceptional vertices not in $A_i'$ when $i\neq 0$, we can run the same procedure and easily ensure that the two (slightly unbalanced) $r$-cliques used are such that all vertices within $B'$ appear in the same part $B_1'$ or $B_2'$.

We can do the same for $(0,\eta)$-exceptional vertices, but we must be a bit careful: when we choose the $r$-cliques with one fewer vertex in $B'$, we cannot necessarily choose where the $1$ vertex in $B'$ is; and thus we first choose all of those cliques and then choose the edges within $A_0'$ (which are then extended to cliques with $1$ extra vertex in $B'$). We therefore need to be able to choose up to $x_0$ (spread) edges not just within $B'=A_0'$ but specifically within $B_1'$ and $B_2'$ as necessary to balance the parity. (They must also have non-bad endpoints.) Here we use that $G[B_1'],G[B_2']$ are almost-complete graphs instead of the argumentation used in \textbf{Step 2} of \cref{sub:connected-subcase}, and then choose extension to cliques with $1$ extra vertex in $A_0'$ where the third vertex in $A_0'$ is in the same part $B_i'$. It is not hard to see that this is possible in a spread manner.

When covering $0$-bad vertices (now that exceptional vertices are removed), we create an $(r-1)$-clique formed with one vertex from each $A_i'$ for $i\in[r-2]$ and our $0$-bad vertex $v$ with various goodness properties on the new vertices, and then extend to a balanced $r$-clique; we can easily guarantee that the final added vertex (which is in $B'$) is in the same part of $B'$ as $v$ using the $n/(4r)$ minimum degree condition within the parts.

\textbf{Step 5: Covering the remainder with \cref{thm:super-regular}.}
After doing all this covering, which only deletes around $O_\alpha(n)$ vertices, we can again use \cref{thm:super-regular} (here \cref{lem:robust-non-sparse} is not really needed). This time, however, we apply it two times, once to a partition where we take appropriately regular subsets of $A_1'',\ldots,A_{r-2}''$ of size $|B_1''|/2$ as well as an appropriately regular bipartition of $B_1''$, and a second time where we apply to the complements of these subsets and an appropriately regular bipartition of $B_2''$. All these partitions and subsets can be found randomly, as the necessary host sets already satisfy appropriately super-regularity. The minimum degree conditions follow since $\beta\ll\eta$ and $B_1'',B_2''$ differ from $n/r$ by an amount depending only on $\beta$.
\end{proof}

\section{A spread distribution on bounded-degree trees}\label{sec:bounded-degree-trees}

Fix $\Delta \in \mb N$ and $\varepsilon>0$. Let $G$ be an $n$-vertex graph satisfying $\delta(G) \ge (1/2+\varepsilon)n$ and let $T$ be an $n$-vertex tree with maximal degree at most $\Delta$. To prove \cref{thm:spanning-trees} it suffices to exhibit an $O(1/n)$-spread distribution on the copies of $T$ in $G$ and then apply \cref{thm:FKNP}. To do so we closely follow the original proof of Koml\'os, S\'ark\"ozy, and Szemer\'edi \cite{KSS95} that $G$ contains at least one copy of $T$. We show that if this algorithm is appropriately randomized then the resulting distribution on copies of $T$ is $O_{\Delta,\varepsilon}(1/n)$-spread.

Most of the embedding algorithm consists of embedding vertices of $T$ into $G$ ``random greedily'' (i.e., the image of each successive vertex is chosen uniformly at random from a set of suitable choices). For this type of algorithm it is natural to analyze ``vertex spread'', which we now define.

\begin{definition}
    Let $X$ and $Y$ be finite vertex sets and let $\mu$ be a distribution over injections $\varphi\colon X \to Y$. For $q \in [0,1]$, we say that $\mu$ is \textit{$q$-vertex-spread} if for every two sequences of distinct vertices $x_1,\ldots,x_k \in X$ and $y_1,\ldots,y_k \in Y$:
    \[
    \mb{P} \left[ \bigwedge_{i=1}^k \varphi(x_i) = y_i \right] \le q^k.
    \]
\end{definition}

We will prove that a randomized version of the Koml\'os, S\'ark\"ozy, and Szemer\'edi tree embedding algorithm is $O(1/n)$-vertex-spread. The next lemma implies that this is sufficient.

\begin{lemma}
    Let $G$ be an $n$-vertex graph and let $T$ be an $n$-vertex tree with maximal degree at most $\Delta \in \mb N$. Suppose that there exists an $O_{\Delta,\varepsilon}(1/n)$-vertex-spread distribution on graph embeddings $\varphi$ of $T$ into $G$. Then $\varphi$ is an $O_{\Delta,\varepsilon}(1/n)$-spread distribution on copies of $T$ in $G$.
\end{lemma}

\begin{proof}
    Let $C\ge 1$ be a constant such that $\varphi$ is $(C/n)$-vertex-spread. Let $\varphi(E(T))$ denote the (random) set of edges in the embedding of $T$ into $G$. We will show that for every edge set $S \subseteq E(G)$:
    \[
    \mb{P} \left[ S \subseteq \varphi(E(T)) \right] \le \left( \frac{\Delta C^2}{n} \right)^{|S|}.
    \]
    This will imply the lemma.

    Let $S \subseteq E(G)$. We may assume that $S$ is a non-empty forest. Denote the number of connected components in $S$ by $\ell$. Let $V(S)$ be the set of vertices incident to $S$. We observe that $|V(S)| = |S|+\ell$. We claim that there are at most $n^\ell \Delta^{|S|}$ embeddings of of $S$ into $T$. Indeed, let $v_1,v_2,\ldots,v_{|S|+\ell}$ be an ordering of $V(S)$, where $v_1,\ldots,v_\ell$ are in distinct connected components and each of the remaining vertices is incident to a vertex that appeared previously (for instance, one may take a breadth-first ordering of $V(S)$). We will count the number of ways to embed $V(S)$ into $T$ one vertex at a time. There are fewer than $n^\ell$ ways to embed $v_1,\ldots,v_\ell$ into $T$. Then, each $v_i$ is incident to a previously embedded vertex. Since $\Delta(T) \le \Delta$ there are at most $\Delta$ choices to embed each $v_i$. Hence, the number of embeddings is at most $n^\ell \Delta^{|S|}$, as claimed.

    Now, for a given embedding $\psi\colon V(S) \to V(T)$, by the vertex-spread assumption for $\varphi$, we have
    \[
    \mb{P} \left[ \bigwedge_{i=1}^{|S|+\ell} \varphi(\psi(v_i)) = v_i \right] \le \left( \frac{C}{n} \right)^{|S|+\ell}.
    \]
    Applying a union bound over the at most $n^\ell \Delta^{|S|}$ choices of $\psi$ we conclude that
    \[
    \mb{P} \left[ S \subseteq \varphi(E(T)) \right] \le n^\ell \Delta^{|S|} \left( \frac{C}{n} \right)^{|S|+\ell} \le \left( \frac{\Delta C^2}{n} \right)^{|S|},
    \]
    where the second inequality follows from the fact that $\ell \le |S|$ and the assumption that $C \ge 1$. This completes the proof.
\end{proof}

The remainder of this section is devoted to proving the next lemma. Together with the previous claim and \cref{thm:FKNP} it implies \cref{thm:spanning-trees}.

\begin{lemma}\label{lem:spread-tree-distribution}
	For every $\Delta \in \mb N$ and $\delta > 0$ there exists some $n_{\ref{lem:spread-tree-distribution}} = n_{\ref{lem:spread-tree-distribution}}(\Delta,\delta) > 0$ and $C_{\ref{lem:spread-tree-distribution}} = C_{\ref{lem:spread-tree-distribution}}(\Delta,\delta)>0$ such that for every graph $G$ on $n \ge n_{\ref{lem:spread-tree-distribution}}$ vertices with $\delta(G) \ge (1/2+\delta)n$ and every tree $T$ on $n$ vertices with $\Delta(T) \le \Delta$ there exists a $(C_{\ref{lem:spread-tree-distribution}}/n)$-vertex-spread distribution on graph embeddings of $T$ into $G$.
\end{lemma}

\subsection{Preliminaries}

In this section we will construct distributions over extensions of graph embeddings in super-regular pairs. We will generally be in the setting where we are attempting to extend a partial embedding of $T$ into $G$ to a larger partial embedding. It is thus helpful to introduce notation for rooted embeddings. Suppose that $H$ and $G$ are graphs, that $R \subseteq V(H)$ is a vertex set, and that $\varphi\colon R \hookrightarrow V(G)$ is an injective partial embedding of $H$. We denote the set of graph embeddings $\wt \varphi\colon H \hookrightarrow G$ that extend $\varphi$ by $X(H,G,R,\varphi)$. We say that a distribution over elements $\wt \varphi \in X(H,G,R,\varphi)$ is \textit{$q$-vertex-spread} if the induced distribution over $\wt{\varphi}|_{V(H) \setminus R}$ is $q$-vertex-spread.

A key fact proved in \cite{KSS95} is that one can embed bounded-degree stars in super-regular pairs. The next lemma shows that this can be done in a spread manner. We need the following definition.

\begin{definition}
    Given a bipartite graph $G=(A,B,E)$ and a vector $\vec{d} = (d_a\colon a \in A) \in \mb N^{A}$ let $\mc S_{\vec d}$ be the graph consisting of the disjoint union of the stars $(S_a\colon a \in A)$, where $S_a = K_{1,d_a}$ for each $a \in A$. Let $R_A$ be the set of the roots of these stars and let $\varphi:R_A \hookrightarrow A$ map the root of $S_a$ to $a$. A \textit{$\vec d$-matching} in $G$ from $A$ to $B$ is an element of $X(\mc S_{\vec d}, G, R_A, \varphi)$.
\end{definition}

The next lemma is a randomized version of \cite[Lemma 2.1]{KSS95}; the proof is similar to that of \cref{lem:bip-pm}.

\begin{lemma}\label{lem:spread-star-embedding}
Let $\delta>0$ and $\Delta \in \mb N$ be fixed. Suppose that $G=(A,B,E)$ is a bipartite graph that is $(\delta,\varepsilon,\delta/2)$-super-regular, with $\varepsilon \leq \delta/(10\Delta)$. Suppose that $\vec d \in \mb N^{A}$ is a vector satisfying $\sum_{a \in A}d_a \le |B|$ and $\max_{a \in A} d_a \le \Delta$. There exists some $C = C(\delta,\Delta)$ and a $(C/|B|)$-vertex-spread distribution over $X(\mc S_{\vec d}, G, R_A, \varphi)$.
\end{lemma}

\begin{proof}
    It suffices to show that there exists an $O(1/|B|)$-spread distribution on subgraphs $H \subseteq G$ that w.h.p.~contain $\vec d$-matchings from $A$ to $B$ (i.e., $X(\mc S_{\vec d}, H, R_A, \varphi) \neq \emptyset$). Indeed, suppose that $\mu$ is such a distribution. Let $\nu$ be the distribution on $X(\mc S_{\vec d}, G, R_A, \varphi)$ obtained by first sampling $H \sim \mu$, conditioned on $H$ containing a $\vec d$-matching from $A$ to $B$, and then choosing an element of $X(\mc S_{\vec d}, H, R_A, \varphi)$ uniformly at random. Let $v_1,\ldots,v_k \in V(\mc S_{\vec d}) \setminus R_A$ be distinct and $u_1,\ldots,u_k \in B$. We wish to show that
    \[
    \mb{P}_{\wt \varphi \sim \nu} \left[ \bigwedge_{i=1}^k \wt \varphi(v_i) = u_i \right] \le (O(1/|B|))^k.
    \]
    Observe that every $v_i$ has a unique neighbor $a_i$ in $\mc S_{\vec d}$, and that $\wt\varphi(v_i) = u_i$ only if $\varphi(a_i)u_i \in E(H)$. Hence
    \begin{align*}
    \mb{P}_{\wt\varphi \sim \nu} \left[ \bigwedge_{i=1}^k \wt \varphi(v_i) = u_i \right] & \le \mb{P}_{H \sim \mu} \left[ \{ \varphi(a_1)u_1,\varphi(a_2)u_2,\ldots,\varphi(a_k)u_k \} \subseteq E(H) | X(\mc S_{\vec d}, H, R_A, \varphi) \neq \emptyset \right]\\
    & \le (1+o(1)) \mb{P}_{H \sim \mu} \left[ \{ \varphi(a_1)u_1,\varphi(a_2)u_2,\ldots,\varphi(a_k),u_k \} \subseteq E(H) \right].
    \end{align*}
    Since $\mu$ is $O(1/|B|)$-spread the last quantity is bounded from above by $(O(1/|B|))^k$, as desired.
    
    We now construct the desired distribution on subgraphs $H \subseteq G$. Let $D$ be a large constant (depending only on $\delta$ and $\Delta$). Let $H' = G(D/|B|)$ (i.e., the binomial subgraph of $G$ with density $D/|B|$). Let $H'' \subseteq G$ be a random graph constructed as follows: For each $v \in A \cup B$ choose, uniformly at random and independently of all other choices, a set of $D$ edges in $G$ incident to $v$ and add them to $H''$. Set $H = H' \cup H''$. Clearly, the distribution on $H$ is $O(D/|B|)$-spread. We will show that w.h.p.~$H$ contains a $\vec d$-matching from $A$ to $B$.
	
    For $X \subseteq A$, let $d_X \coloneqq \sum_{a\in X} d_a$. It suffices to show that w.h.p.~$H$ satisfies the K\"onig--Hall criterion
    \[
    \forall X \subseteq A,  |N_H(X)| \ge d_X.
    \]

    For notational conciseness set $\varepsilon \coloneqq \delta/(10\Delta)$. We first consider $0 < |X| \le \varepsilon |B|$. We will bound the probability that for some $Y \subseteq B$ with $|Y| = d_X$ we have $N_{H''}(X) \subseteq Y$. Indeed, given such $X$ and $Y$, we have
    \[
    \mb{P} \left[ N_{H''}(X) \subseteq Y \right] \le \left( \frac{\binom{|Y|}{D}}{\binom{\delta|B|/2}{D}} \right)^{|X|} \leq \left( \frac{2e|Y|}{\delta |B|} \right)^{D|X|} \leq \left( \frac{2e\Delta|X|}{\delta|B|} \right)^{D|X|}.
    \]
    We now apply a union bound over choices of such $X$ and $Y$ to obtain:
    \[
    \sum_{k=1}^{\varepsilon|B|} \binom{|A|}{k}\binom{|B|}{\Delta k} \left( \frac{2e\Delta k}{\delta |B|} \right)^{D k} \le \sum_{k=1}^{\varepsilon|B|} \left( \frac{2^D e^{1+\Delta+D} \Delta^{D-\Delta}}{\delta^D} \times \frac{k^{D-\Delta-1}}{|B|^{D-\Delta-1}} \right)^k = o(1).
    \]
    
    Next, we consider $|X|$ such that $\varepsilon |B| \le |X|$ and $d_X \le (1-\varepsilon)|B|$. Let $Y \in \binom{B}{d_X}$ (so, in particular, $|Y| \ge |X| \ge \varepsilon |B|$). Since $G$ is $(\delta,\varepsilon)$-regular we have $e_G(X,B \setminus Y) > \delta|X|(|B|-|Y|)/2 \geq \delta\varepsilon^2|B|^2\!/2$. Therefore
    \[
    \mb{P} \left[ e_{H'}(X,Y) = 0 \right] \le \left( 1 - \frac{D}{|B|} \right)^{\delta\varepsilon^2|B|^2\! / 2} \le \exp \left( - \frac{D\delta\varepsilon^2}{2} |B| \right).
    \]
    Since there are fewer than $2^{|B|}$ choices for $X$ and $Y$, applying a union bound, the probability that there exist such $X$ and $Y$ with $e_{H'}(X,B\setminus Y) = 0$ is at most $2^{2|B|}\exp \left( -D\delta\varepsilon^2 |B|/ 2 \right)$, which tends to zero provided $D$ is sufficiently large.
	
    Observe that the two cases above complete the proof when $|B| > 2d_A$. Henceforth, we assume that $|B| \le 2d_A \le 2\Delta|A|$.
	
    It remains to consider $X$ such that $d_X > (1-\varepsilon)|B|$. This implies $|A \setminus X| \le \varepsilon|B|$ (since $|B| \ge d_A = d_X + d_{A \setminus X} \ge (1-\varepsilon)|B| + |A \setminus X|$). Let $Y \subseteq B$ satisfy $|Y| = d_X-1$. Thus $|B \setminus Y| \ge |A \setminus X|$ (since $|B\setminus Y| \geq d_A-|Y|>d_A-d_X=d_{A\setminus X} \geq |A\setminus X|$). Now, if $N_H(X) \subseteq Y$ then $N_H(B \setminus Y) \subseteq A \setminus X$. Hence, there exists a subset of $B$ of size $|B|-d_X+1 \leq 2\varepsilon|B|$, all of whose neighbors are contained in a smaller set. We use a union bound to show that w.h.p.~there is no such set:
    \[
    \sum_{k=1}^{2\varepsilon|B|} \binom{|B|}{k} \binom{|A|}{k} \left( \frac{\binom{k}{D}}{\binom{\delta|A|/2}{D}} \right)^k \leq \sum_{k=1}^{2\varepsilon|B|} \left( \frac{e^2 2^{D+1} \Delta}{\delta^D} \times \frac{k^{D-2}}{|A|^{D-2}} \right)^k = o(1).
    \]
    
    Thus, w.h.p.~$H$ satisfies the K\"onig--Hall condition and contains a $\vec d$-matching from $A$ to $B$.
\end{proof}

A second key claim that is proved in \cite{KSS95} allows the embedding of forests of length-$3$ paths into super-regular pairs. We use the following definitions.

\begin{definition}
    A \textit{four-layer $(d,\varepsilon)$-super-regular graph} is a graph $G=(V,E)$ with a vertex partition $V = V_1 \cup V_2 \cup V_3 \cup V_4$, where all parts are the same size, and for every $i=1,2,3$ the induced bipartite graph $G[V_i,V_{i+1}]$ is $(d^+,\varepsilon)$-super-regular.
    
    For $m \in \mb N$ let $\mc P_m$ be the graph consisting of $m$ vertex disjoint copies of length-$3$ paths. Denote by $V_{\mr{out}}(\mc P_m)$ the set of $2m$ leaves in $\mc P_m$ and denote by $V_{\mr{in}}(\mc P_m)$ the set of $2m$ inner vertices in $\mc P_m$.

    If $\pi\colon V_1 \to V_4$ is a bijection then $\varphi\colon V_{\mr{out}}(\mc P_{|V_1|}) \to V_1 \cup V_4$ is \textit{$\pi$-respecting} if for every $v \in V_1$ both $\varphi^{-1}(v)$ and $\varphi^{-1}(\pi(v))$ belong to the same path in $\mc P_{|V_1|}$.

\end{definition}

The next lemma is a randomized version of \cite[Theorem 2.1]{KSS95}.

\begin{lemma}\label{lem:spread-path-embedding}
    For every $d>0$ there exist $C = C_{\ref{lem:spread-path-embedding}}(d),\varepsilon = \varepsilon_{\ref{lem:spread-path-embedding}}(d),n_{\ref{lem:spread-path-embedding}}(d)>0$ such that the following holds for all $n\ge n_{\ref{lem:spread-path-embedding}}$. If $G$ is a four-layer $(d,\varepsilon)$-super-regular graph on $4n$ vertices then for any bijection $\pi:V_1\to V_4$ and any $\pi$-respecting bijection $\varphi\colon V_{\mr{out}}(\mc P_n) \to V_1 \cup V_4$ there exists a $(C/n)$-vertex-spread distribution on $X(\mc P_n, G, V_{\mr{out}}(\mc P_n), \varphi)$.
\end{lemma}

\begin{proof}
    We prove the lemma by applying \cref{thm:super-regular} to the auxiliary graph $G'$ which is obtained by identifying $V_1$ and $V_4$ according to $\pi$. In other words, the vertex set of $G'$ is $V_1\cup V_2 \cup V_3$ and the edge set consists of $G[V_1,V_2] \cup G[V_2,V_3]$ and $ \{ xy \in V_1 \times V_3\colon \pi(x)y \in G[V_4,V_3] \} $. As long as $\varepsilon$ is sufficiently small and $n$ is sufficiently large then $G'$ satisfies the assumptions of \cref{thm:super-regular}. Hence, for some $C>1$ there exists a $(Cn^{-2})$-spread distribution $\mu$ on perfect matchings in the $3$-clique complex $\mc H$ of $G'$.

    Observe that there is a natural correspondance between elements $X(\mc P_n, G, V_{\mr{out}}(\mc P_n), \varphi)$ and perfect matchings in $\mc H$. Explicitly, the embedding $\wt\varphi \in X(\mc P_n, G, V_{\mr{out}}(\mc P_n), \varphi)$ corresponds to the perfect matching $M \subseteq \mc H$ consisting of all triples $v_1v_2v_3 \in \mc H$ such that $v_1v_2v_3\pi(v_1)$ is the image under $\wt\varphi$ of a path in $\mc P_n$. Hence, $\mu$ induces a distribution $\nu$ on $X(\mc P_n, G, V_{\mr{out}}(\mc P_n), \varphi)$. We will show that $\nu$ is $(C/n)$-vertex-spread.

    Let $v_1,\ldots,v_k \in V_{\mr{in}}(\mc P_{n})$ be a sequence of distinct vertices and let $u_1,\ldots,u_k \in V_2 \cup V_3$. Let $M \sim \mu$ and let $\wt\varphi \in X(\mc P_n, G, V_{\mr{out}}(\mc P_n), \varphi)$ be the corresponding embedding. We wish to show that
    \[
    \mb{P} \left[ \bigwedge_{i=1}^k \wt\varphi(v_i) = u_i \right] \le (C/n)^k.
    \]
    Observe that every vertex $v \in V_{\mr{in}}(\mc P_n)$ has a unique neighbor $v' \in V_{\mr{in}}(\mc P_n)$ and these vertices are connected to a unique pair in $V_{\mr{out}}(\mc P_n)$. Hence,  specifying the image of $v$ and $v'$ is equivalent to prescribing that a specific triangle appear in $M$. Let $a$ be the number of vertices in $v \in \{v_1,\ldots,v_k\}$ such that $v' \in \{v_1,\ldots,v_k\}$. For each of the remaining $k-a$ vertices $v$ there are at most $n$ ways to embed $v'$. Using the spread of $\mu$, for each such choice, the probability that all the corresponding triangles will be in $M$ is at most $(Cn^{-2})^{a/2+k-a}$. Applying a union bound we conclude that
    \[
    \mb{P} \left[ \bigwedge_{i=1}^k \wt\varphi(v_i) = u_i \right] \le n^{k-a} \left( \frac{C}{n^2} \right)^{a/2 + k-a} \le \left( \frac{C}{n} \right)^k.\qedhere
    \]
\end{proof}

Before stating the next lemma we introduce some notation. Suppose that $\varphi$ is a partial embedding of a graph $H$ into a graph $G$, defined on $D \subseteq V(H)$. We call the vertices in $V(G) \setminus \varphi(D)$ \textit{unoccupied by $\varphi$}. For $v \in V(H) \setminus D$ we write $A(\varphi,v,H,G)$ for the set of vertices $u \in V(G)$ that are unoccupied by $\varphi$ and that are adjacent to the images of all embedded neighbors of $v$ (formally, for every $w \in D$ such that $wv \in E(H)$ it holds $u\varphi(w) \in E(G)$). This is the set of \textit{available} locations for $v$. If $H$ and $G$ are bipartite (as will always be the case for us) we additionally restrict vertices so that they respect the bipartition.

The next lemma, which is closely related to the blow-up lemma, allows one to extend an embedding of a graph within a regular pair in a spread manner. It also allows the designation of ``buffers'' which are either target locations for certain vertices or sets that should be left mostly unoccupied.

\begin{lemma}\label{lem:spread-blow-up-precursor}
	Let $\Delta \in \mb N$ and $\alpha, d \in (0,1)$. There exists some $\varepsilon = \varepsilon(\Delta,\alpha,d) > 0$ such that for every $k \in \mb N$ the following holds for every $n \ge n_0 = n_0(k,\Delta,\alpha,d)$. Suppose that $G = (A,B,E)$ is a $(d^+,\varepsilon)$-regular pair with $2n \ge |A|,|B| \ge n$ and that $H = (C,D,E_H)$ is a bipartite graph satisfying
	\begin{itemize}
            \item every connected component of $H$ has size at most $\varepsilon^2 n$,
 
		\item $|C| \le (1-\alpha)|A|$ and $|D| \le (1-\alpha)|B|$, and
		
		\item $\Delta(H) \le \Delta$.
	\end{itemize}
	Suppose further that $F \subseteq C \cup D$ and $B_1,B_2 \subseteq A \cup B$ satisfy
	\begin{itemize}
		\item $|F \cap C| = |B_1 \cap A|$ and $|F \cap D| = |B_1 \cap B|$.
		
		\item $|B_2 \cap A| = |A|-|C|$ and $|B_2 \cap B| = |B|-|D|$.

            \item $B_1 \cap B_2 = \emptyset$.
	\end{itemize}
	Finally, suppose that for a set $S \subseteq V(H)$ of size at most $k$ there is a graph embedding $\varphi_0\colon H[S] \hookrightarrow G$ such that for every $v \in V(H) \setminus S$ it holds $|A(\varphi_0,v,H,G)| \ge (d/2)^\Delta n$. There exists an $O(1/n)$-vertex-spread distribution on graph embeddings $\varphi \in X(H,G,S,\varphi_0)$ such that
	\begin{itemize}
		\item At most $\varepsilon^{1/3} n$ vertices of $B_2$ are occupied by $\varphi$, and
		
		\item At most $\varepsilon^{1/3} n$ vertices of $F$ are not mapped by $\varphi$ to $B_1$.
	\end{itemize}
\end{lemma}

\begin{proof}
    We describe a randomized embedding procedure that embeds the vertices of $V(H) \setminus S$ one by one, where the image of each vertex is chosen uniformly at random from a set of size $\Omega(n)$. This guarantees the vertex-spread. Additionally, we design the algorithm in such a way that vertices from $B_2$ are chosen at most $\varepsilon^{1/3} n$ times, and vertices from $F$ will fail to be mapped to $B_1$ at most $\varepsilon^{1/3} n$ times. The lemma follows by considering the corresponding distribution over embeddings.

    Naively, one might try a simple random greedy algorithm, where each successive vertex is mapped to a uniformly random available vertex, all while avoiding $B_2$. However, special care is needed to handle the buffers. A particular concern is that since (for example) $|A \setminus B_2| = |C|$, if we have embedded nearly all vertices in $C$ then the remaining set of unoccupied vertices in $A \setminus B_2$ might be so small that regularity fails, in which case the embedding algorithm might get stuck.

    To circumvent this we set aside additional buffer zones within $B_2$ which we allow ourselves to use in the embedding procedure. These are large enough to guarantee that regularity never fails but still small enough that almost all vertices in $B_2$ remain unoccupied. We define two buffer zones, as follows. Let $Z_1$ and $Z_2$ be disjoint subsets of $B_2$, each consisting of $\sqrt{\varepsilon} n$ vertices in each of $B_2 \cap A$ and $B_2 \cap B$ (so that $|Z_1 \cap A| = |Z_1 \cap B| = |Z_2 \cap A| = |Z_2 \cap B| = \sqrt{\varepsilon} n$).
    
    For each vertex we will now define a target set where it will be embedded. Let $N(S)$ denote the set of neighbors of $S$. For $v \in N(S)$ let $X(v) \subseteq A(\varphi_0,v,H,G) \setminus (Z_1 \cup Z_2)$ be a set of size $\sqrt{\varepsilon} (d/2)^\Delta n/(2\Delta k)$ such that:
    \begin{itemize}
        \item For distinct $u,v \in N(S)$ the sets $X(u)$ and $X(v)$ are mutually disjoint;

        \item For $u,v \in N(S)$ that lie in different sides of the partition of $H$ the pair $(X(u),X(v))$ is a $((4d/5)^+,\sqrt{\varepsilon})$-regular pair in $G$; and

        \item for every $v \in N(S)$, if $(W_1,W_2)$ is the ordering of $A,B$ such that $X(v) \subseteq W_1$, then $(X(v),W_2)$ is a $((d/2)^+,\varepsilon^{3/4})$-regular pair.
    \end{itemize}
    We note that such a choice of set $\{X(v)\}_{v\in N(S)}$ is possible since the regularity properties are satisfied w.h.p.\ by choosing uniformly random disjoint sets of the appropriate size. We write $Z_3 \coloneqq \bigcup_{v \in N(S)}X(v)$. Observe that $|Z_3| = |N(S)| \sqrt{\varepsilon} (d/2)^\Delta n/(2\Delta k) \leq \sqrt{\varepsilon} n / 2$.

    Now, for $v \in V(H)$ let $W(v) = A$ if $v \in C$ and let $W(v) = B$ if $v \in D$. For $v \in F \setminus (S \cup N(S))$ set its target set as $X(v) \coloneqq W(v) \cap ( (B_1 \setminus Z_3) \cup Z_1)$. Finally, for $v \in V(H) \setminus (S \cup N(S) \cup F)$ set $X(v) \coloneqq W(v) \cap (V(G) \setminus (B_1 \cup B_2 \cup Z_3) \cup Z_2)$.

    The upshot of choosing the sets $X(v)$ in this way that is that if (as will indeed be the case) we succeed in embedding every $v$ into $X(v)$ then vertices in $F$ will fail to be embedded into $B_1$ at most $\varepsilon^{1/3}n$ times. This is because such a vertex is either in $N(S)$ or else, if it is not embedded into $B_1$, it must have been embedded into $Z_1$. Since $|Z_1| + |N(S)| \leq 2\sqrt{\varepsilon}n + k\Delta \leq \varepsilon^{1/3}n$ this is an upper bound on the number of ``wrongly embedded'' vertices in $B_1$. Similarly, vertices in $B_2$ will be used at most $|Z_1|+|Z_2|+|Z_3| \leq 5\sqrt{\varepsilon}n \leq \varepsilon^{1/3}n$ times. Additionally, for as long as we embed vertices only into their target sets, for $v \in V(H) \setminus (S \cup N(S))$ there will always be at least $\sqrt{\varepsilon}n$ unoccupied vertices in its target set. This is ensured by the excess vertices added by the buffers: each vertex in $X(v)$ appears in at most $|X(v)|-\sqrt{\varepsilon}n$ other target sets.

    We observe that for every $v \in V(G) \setminus S$ it holds $|X(v)| \geq \sqrt{\varepsilon} (d/2)^\Delta n/2 = \Omega_{d,k,\Delta}(n)$ and, furthermore, if $v \notin N(S)$ then $|X(v)| \geq \min\{\alpha,1-\alpha\} n$. Additionally, if $u,v \in V(H)$ lie in different sides of the partition of $H$ then the pair $(X(u),X(v))$ is $((3d/5)^+,\varepsilon^{3/5})$-regular.

    The embedding algorithm is as follows.
    \begin{itemize}
	\item Set $\varphi = \varphi_0$.

        \item Let $v_1,v_2,\ldots,v_{|V(H)|-|S|}$ be an ordering of $V(H) \setminus S$ where the vertices lying in each connected component of $H$ appear as a contiguous interval of vertices.
		
	\item For each $i = 1,\ldots,|V(H)|-|S|$:
		
	\begin{itemize}
		\item Let $N_i \subseteq V(H)$ be the set of neighbors of $v_i$ that have not already been embedded.
			
		\item Let $X_i$ be the set of vertices $v \in A(\varphi,v_i,H,G) \cap X(v_i)$ that, for each $u \in N_i$, satisfy $\deg_G(v,A(\varphi,u,H,G) \cap X(u)) \ge \frac{d}{2}|A(\varphi,u,H,G) \cap X(u)|$.
			
		\item Choose some $v \in X_i$ uniformly at random and update $\varphi(v_i) = v$.
	\end{itemize}
    \end{itemize}
	
    We claim that the algorithm is guaranteed to succeed and that the set $X_i$ always has size at least $\varepsilon^2 n = \Omega_{\Delta,d,\alpha}(n)$. Indeed, we note that at step $i$ of the algorithm, if each of the sets $A(\varphi,v_i,H,G)\cap X(v_i)$ and $\{A(\varphi,u,H,G) \cap X(u)\}_{u \in N_i}$ contains at least a $(\Delta+1)\varepsilon^{3/5}$-fraction of the vertices in their respective target sets then by regularity, for every $u \in N_i$, the inequality $\deg_G(v,A(\varphi,u,H,G) \cap X(u)) \ge \frac{d}{2}|A(\varphi,u,H,G) \cap X(u)|$ holds for all but at most $\varepsilon^{3/5} |X(v_i)|$ vertices in $A(\varphi,v_i,H,G)\cap X(v_i)$. Hence, if this is the case, then $|X_i| \geq (\Delta+1)\varepsilon^{3/5}|X(v_i)| - \Delta \varepsilon^{3/5} |X(v_i)| = \varepsilon^{3/5} |X(v_i)|  \geq \varepsilon^2 n$.
 
    We now show that as long as vertex $v$ has not been embedded the set $A(\varphi,v,H,G) \cap X(v)$ indeed contains at least a $(\Delta+1)\varepsilon^{3/5}$-fraction of $X(v)$. We consider two cases. First, if $v \in N(S)$ then $X(v)$ is disjoint from all other target sets. Thus, it decreases only when a neighbor of $v$ is embedded, in which case (by the algorithm's design) it decreases by a factor $f$ with $f \geq d/2$. Hence, since $v$ has at most $\Delta$ neighbors the number of available locations is always at least $(d/2)^\Delta |X(v)| \geq (\Delta+1)\varepsilon^{3/5}|X(v)|$, where the inequality holds provided $\varepsilon$ is sufficiently small.

    In the second case $v \in V(H) \setminus (S \cup N(S))$. In this case the number of locations can decrease in two ways: in the first, a neighbor of $v$ is embedded, in which case the number of available locations can decrease multiplicatively by a factor of $f$ with $f \geq d/2$. As before, this can happen at most $\Delta$ times. In the second, the set of available locations can decrease if some vertex in $X(v)$ becomes occupied, in which case it decreases by $1$. However, as long as no neighbor of $v$ is embedded the number of available locations is precisely the number of unoccupied vertices in $X(v)$ which is of size at least $\sqrt{\varepsilon} n$ (this is the excess space guaranteed by the buffers). Recall that the embedding is done connected component by connected component and that each connected component has size at most $\varepsilon^2n$. Hence, after the first neighbor of $v$ is embedded at most $\varepsilon^2n$ additional vertices in $X(v)$ become occupied before $v$ itself is embedded. Thus the number of available locations for $v$ is at least $(\sqrt{\varepsilon} - \varepsilon^2)n (d/2)^\Delta \geq (\Delta+1)\varepsilon^{3/5}|X(v)|$, where the last inequality holds provided $\varepsilon$ is sufficiently small in terms of $d$ and $\alpha$.
\end{proof}

\subsection{A randomized tree embedding algorithm}

We now describe our adaptation of the tree embedding algorithm in \cite[Section 5.5]{KSS95}. That algorithm has eight steps. Steps 1--6 are ``preprocessing'' of the tree $T$ and the graph $G$, and we make no changes to these steps. The output of these steps is a regularization of $G$ together with an assignment of $V(T)$ to regularized clusters that determine into which cluster each tree vertex will be embedded. We do not give details of this construction; instead, we simply list its salient properties in \cref{clm:preprocessing-output}, below. For additional details we refer the reader to \cite{KSS95}.

The actual embedding is carried out in steps 7 and 8. One approach to prove \cref{lem:spread-tree-distribution} would be to randomize the embedding strategy and note that it is $O(1/n)$-vertex-spread. Indeed, Step 7 consists of a random greedy algorithm, where vertices of $T$ are embedded greedily into $G$ one at a time, where for each vertex there are $\Omega_{\Delta,\varepsilon}(n)$ suitable choices. By choosing the embedding random-greedily it becomes appropriately spread. In Step 8 vertices are embedded by employing derandomized versions of either \cref{lem:spread-star-embedding} or \cref{lem:spread-path-embedding}; if these are replaced by the appropriate randomized counterpart the embedding strategy as a whole regains the necessary spread.

We do not follow the embedding strategy in \cite{KSS95} exactly. This is mostly for organizational purposes, as well as to avoid duplicating large parts of \cite{KSS95} verbatim. Instead, we first embed a small number of vertices that serve as ``bridges'' between super-regular pairs in the regularized graph. We then embed the remainder of the tree into the super-regular pairs.

We begin by describing the outcome of the preprocessing steps in \cite{KSS95}. Following \cite{KSS95}, we will make $T$ rooted by fixing an arbitrary root. Any subgraph of $T$ can then be viewed as a rooted forest. We define a \textit{secondary leaf} in a forest as a non-leaf vertex all of whose children are leaves.

\begin{claim}\label{clm:preprocessing-output}
	In the setting of \cref{lem:spread-tree-distribution}, for every $\varepsilon>0$ there exists some $M = M(\varepsilon)>0$, $\alpha = \alpha(\Delta)>0$, such that if we fix any vertex $r \in V(T)$ as a root then there exists a decomposition of $G$ into clusters $\mc C$ with the following properties.
	\begin{enumerate}
		\item $|\mc C| \le M$.
		
		\item For every $C \in \mc C$ it holds $2n/M \ge |C| \ge n/(2M)$.
		
		\item There exists a perfect matching $\mc M$ of the clusters in $\mc C$ such that every pair in $\mc M$ is $((\delta/2)^+,\varepsilon)$-super-regular. For $C \in \mc C$ we denote its match by $C'$.
	\end{enumerate}
	There also exists an assignment $a\colon V(T) \to \mc C$, a set $S \subseteq V(T)$, and a constant $K=K(\Delta,\delta,\varepsilon)$ with the following properties:
	\begin{enumerate}
		\item $|S| \le K$.
		
		\item For every $C \in \mc C$ it holds $|a^{-1}(C)| = |C|$.
		
		\item For every edge $uv \in E(T)$, the pair $(a(u),a(v))$ is $\varepsilon$-regular with density at least $\delta/2$.
		
		\item For every edge $uv \in E(T)$, if $(a(u),a(v)) \notin \mc M$ then $u,v \in S$.
	\end{enumerate}
	Finally, for every cluster pair $(C,C') \in \mc M$, let $F_{C,C'} = T[a^{-1}(C),a^{-1}(C')]$ be the subforest of $T$ that is spanned by the vertices assigned to $C$ and $C'$. Then each connected component of $F_{C,C'}$ has size at most $\varepsilon^2 n$. Additionally, there exist sets $F_{C,C'}^1,F_{C,C'}^2 \subseteq V(F_{C,C'}) \setminus S$ such that one of the following holds:
	\begin{enumerate}
		\item $F_{C,C'}^2$ consists of $\alpha |C|$ leaves of $F_{C,C'}$, equally divided between $a^{-1}(C)$ and $a^{-1}(C')$, and $F_{C,C'}^1$ consists of the parents (within $F_{C,C'}$) of $F_{C,C'}^2$.
		
		\item $F_{C,C'}^2$ consists of $\alpha |C|$ secondary leaves in $C$ and their children, and $F_{C,C'}^1$ consists of the parents of the secondary leaves (all within $F_{C,C'}$).
		
		\item $F_{C,C'}^2$ consists of $\alpha |C|$ secondary leaves in $C'$ and their children, and $F_{C,C'}^1$ consists of the parents of the secondary leaves (all within $F_{C,C'}$).
		
		\item For a set of $\alpha |C|$ vertex-disjoint length-$3$ paths in $F_{C,C'}$ in which the internal vertices all have degree $2$, $F_{C,C'}^2$ consists of the $2\alpha|C|$ internal vertices and $F_{C,C'}^1$ consists of the $2\alpha|C|$ endpoints.
	\end{enumerate}
\end{claim}

%TODO add a little motivation
We observe that in the setup above, for every $v \in V(T) \setminus S$, if $a(v) = C$ then all neighbors of $v$ are assigned to $C'$.

We are ready to describe the randomized embedding procedure.

\begin{proof}[Proof of \cref{lem:spread-tree-distribution}]
    Take the output of \cref{clm:preprocessing-output} with respect to $\varepsilon$ sufficiently small so that \cref{lem:spread-blow-up-precursor} can be applied to the super-regular pairs (provided $n$ is sufficiently large).

	We initialize a partial embedding $\varphi$ whose domain is the empty set. We call a vertex of $G$ \textit{occupied} if a vertex of $T$ has been assigned to it (otherwise it is \textit{unoccupied}).

	We begin by embedding the vertices of $S$. For this we will use a random greedy algorithm. The remaining vertices of $T$ will be embedded using \cref{lem:spread-blow-up-precursor,lem:spread-star-embedding,lem:spread-path-embedding}.
	
	Let $s_1,\ldots,s_k$ be an ordering of the vertices of $S$ where each $s_i$ is incident (in $T$) to at most one vertex that precedes it in the ordering. (Such an ordering is possible since $T$ is a tree.) Iterating through $i=1,\ldots,k$, let $A_i \subseteq a(s_i)$ be the set of vertices that are adjacent to all $\varphi(u)$ for all neighbors of $s_i$ that precede it in the order (of which there is at most one). Then, let $B_i \subseteq A_i$ be the set of vertices $v \in A_i$ such that $\deg_G(v,a(u)) \ge \delta|a(u)|/3$ for every $u \in V(T)$ that is adjacent to $s_i$. Choose some $v \in B_i$ uniformly at random, and set $\varphi(s_i) = v$.
	
	Using the fact that in this stage we embed only $O(1)$ vertices and the regularity properties of the decomposition, there are always at least $\delta n / (8M)$ choices for each embedding. Hence $\varphi|_{S}$ is $8M / (\delta n) = O(1/n)$-vertex-spread.
	
	We now embed the remaining vertices of $T$. We do this separately for each cluster pair $(C,C') \in \mc M$. Let $(C,C') \in \mc M$. We will use \cref{lem:spread-blow-up-precursor} to embed the vertices in $F_{C,C'}$ besides those in $F_{C,C'}^2$, and then use either \cref{lem:spread-star-embedding} or \cref{lem:spread-path-embedding} to embed $F_{C,C'}^2$. 
	Before applying \cref{lem:spread-blow-up-precursor} we set aside buffer zones in $(C,C')$. For $i=1,2$, let $B_i \subseteq C \setminus S$ have size $|F_{C,C'}^i \cap a^{-1}(C)|$. Similarly, let $B'_i \subseteq C' \setminus S$ have size $|F_{C,C'}^i \cap a^{-1}(C')|$. We choose these sets in such a way that for all $i,j=1,2$, the pairs $(B_i,B'_j)$ are $((\delta/3)^+,2\varepsilon)$-super-regular. We also choose the buffer zones so that they are mutually disjoint. (To see that this is possible note that if appropriately-sized disjoint sets are chosen uniformly at random then w.h.p.~they satisfy the super-regularity.)
	
	We now apply \cref{lem:spread-blow-up-precursor} to extend $\varphi$ to a partial embedding that embeds $F_{C,C'} \setminus F_{C,C'}^2$ such that:
	\begin{enumerate}
		\item All but $10\varepsilon|C|$ vertices of $F_{C,C'}^1$ are embedded to $B_1 \cup B_1'$.
		
		\item The set of unoccupied vertices differs from $B_2 \cup B_2'$ by at most $10\varepsilon|C|$.
	\end{enumerate}
	
    Let $D_1 = \varphi(F_{C,C'}^1) \cap C$ and $D_1' = \varphi(F_{C,C'}^1) \cap C'$. Let $D_2$ and $D_2'$ be the set of unoccupied vertices in $C$ and $C'$ respectively. Observe that for every $i,j \in \{1,2\}$ the pair $(D_i,D_j')$ is $((\delta/10)^+,20\varepsilon)$-super-regular. It remains to extend $\varphi$ so that it embeds $F_{C,C'}^2$ into $D_2 \cup D_2'$. We consider three cases, depending on how $F_{C,C'}^2$ was constructed.
	
	In the first case $F_{C,C'}^2$ is a set of leaves, evenly divided between ``even'' leaves in $a^{-1}(C)$ and ``odd'' leaves in $a^{-1}(C')$. We apply \cref{lem:spread-star-embedding} twice, first to match the even leaves to their parents (which are already embedded) and then to match the odd leaves to their parents.
	
	In the second and third cases $F_{C,C'}^2$ consists of $\alpha|C|$ secondary leaves (all in either $a^{-1}(C)$ or $a^{-1}(C')$) and their children, and $F_{C,C'}^1$ consists of the secondary leaves' parents. We again apply \cref{lem:spread-star-embedding} twice: first to embed the secondary leaves and then to embed their children.
	
	Finally, in the fourth case, $F_{C,C'}^2$ consists of the internal vertices of $\alpha|C|$ vertex-disjoint paths (in which all vertices in $F_{C,C'}^2$ have degree $2$). We apply \cref{lem:spread-path-embedding} to embed the desired length-$3$ paths.
	
	In all cases, \cref{lem:spread-star-embedding,lem:spread-path-embedding} ensure that the embedding is completed in an $O(1/n)$-vertex spread manner, as desired.
\end{proof}

\bibliographystyle{amsplain0.bst}
\bibliography{main.bib}

\end{document}